\documentclass{article}
\usepackage[T1]{fontenc}
\usepackage{amsthm, fullpage} 
\usepackage{amsmath}
\usepackage{amssymb}
\usepackage{amsfonts}
\usepackage{eucal}
\usepackage[all]{xy}
\usepackage[frenchb,english]{babel}
\usepackage{pslatex}
\usepackage[pagebackref]{hyperref}

\newtheorem{prop}{Proposition}[section]
\newtheorem{theo}[prop]{Théor\`eme}
\newtheorem*{theo**}{Théorème}
\newtheorem{coro}[prop]{Corollaire}
\newtheorem*{conj*}{Conjecture}
\newtheorem{lemm}[prop]{Lemme}

\theoremstyle{definition}

\newtheorem{vide}[prop]{}
\newtheorem{defi}[prop]{Définition}
\newtheorem*{defi*}{Définition}

\theoremstyle{remark}
\newtheorem{rema}[prop]{Remarques}

\newtheorem{nota}[prop]{Notations}

\numberwithin{equation}{prop}

\newcommand{\riso}{ \overset{\sim}{\longrightarrow}\, }
\newcommand{\liso}{ \overset{\sim}{\longleftarrow}\, }

\newcommand{\Spec}{\mathrm{Spec}\,}
\newcommand{\Spf}{\mathrm{Spf}\,}

\renewcommand{\sp}{\mathrm{sp}}
\renewcommand{\det}{\mathrm{det}}

\renewcommand{\AA}{{\mathcal{A}}}

\newcommand{\FF}{{\mathcal{F}}}
\newcommand{\B}{{\mathcal{B}}}

\newcommand{\E}{{\mathcal{E}}}
\newcommand{\G}{{\mathcal{G}}}

\newcommand{\D}{{\mathcal{D}}}

\newcommand{\PP}{{\mathcal{P}}}

\renewcommand{\O}{{\mathcal{O}}}

\newcommand{\V}{\mathcal{V}}
\renewcommand{\S}{\mathcal{S}}
\newcommand{\T}{{\mathcal{T}}}
\newcommand{\Y}{\mathcal{Y}}
\newcommand{\ZZ}{\mathcal{Z}}
\newcommand{\X}{\mathfrak{X}}

\newcommand{\U}{\mathfrak{U}}

\newcommand{\A}{\mathbb{A}}
\renewcommand{\P}{\mathbb{P}}
\newcommand{\F}{\mathbb{F}}
\newcommand{\C}{\mathbb{C}}
\newcommand{\DD}{\mathbb{D}}
\renewcommand{\L}{\mathbb{L}}
\newcommand{\R}{\mathbb{R}}
\newcommand{\Q}{\mathbb{Q}}
\newcommand{\Z}{\mathbb{Z}}
\newcommand{\N}{\mathbb{N}}

\newcommand{\hdag}{  \phantom{}{^{\dag} }    }

\begin{document}
\title{$\D$-modules arithmétiques surholonomes\\
(Overholonomic arithmetical $\D$-modules)}
\author{Daniel Caro \footnote{L'auteur a bénéficié du soutien du réseau européen TMR \textit{Arithmetic Algebraic Geometry}
(contrat numéro UE MRTN-CT-2003-504917).}}
\date{}

\maketitle



\begin{abstract}
\selectlanguage{english}
Let $k$ be a perfect field of characteristic $p >0$, $U$ be a variety over $k$ and $F$ be a power of Frobenius.
  We construct the category of overholonomic arithmetical ($F$-)$\D$-modules over $U$ and
  the category of overholonomic ($F$-)complexes over $U$.
  We show that overholonomic complexes over $U$ are stables by direct images,
inverse images, extraordinary inverse images, extraordinary direct images, dual functors.
Moreover, when $U$ is smooth,
we check that unit-root overconvergent $F$-isocrystals on $U$
are overholonomic. This implies that they are holonomic, which proves in part a Berthelot's conjecture.
\end{abstract}

\selectlanguage{frenchb}
\begin{abstract}
  Soient $k$ un corps parfait de caractéristique $p>0$, $U$ une variété sur $k$ et $F$ une puissance de Frobenius.
  Nous construisons la catégorie des ($F$-)$\D$-modules arithmétiques surholonomes sur $U$ et celle des
   ($F$-)$\D$-complexes surholonomes sur $U$. Nous montrons que les complexes surholonomes sur $U$ sont stables par images
   directes, images inverses, images inverses extraordinaires, images directes extraordinaires, foncteurs duaux.
   De plus, dans le cas lisse, nous vérifions que les $F$-isocristaux surconvergents unités sont surholonomes.
   Cela implique que ceux-ci sont holonomes, ce qui prouve en partie une conjecture de Berthelot.
\end{abstract}

\maketitle

\tableofcontents

\begin{center}
{\bf Introduction}
\end{center}

Donnons-nous $\F _q$, un corps fini de caractéristique $p$, et $X$ une variété sur $\F _q$
(i.e. un $\F _q$-schéma séparé et de type fini).
Afin de prouver la rationalité de la fonction zêta de Weil associée à $X$ et surtout de leur donner une
interprétation cohomologique qui fût conjecturée par Weil, l'idée de Grothendieck
fût de généraliser ces fonctions à des coefficients sur $X$ plus généraux, le {\it coefficient constant}
redonnant la fonction zêta de $X$.
Ces coefficients,
les $\Q _l$-faisceaux constructibles (pour la topologie étale),
avec $l$ un nombre premier différent de $p$,
vérifient les trois propriétés fondamentales suivantes :
\begin{enumerate}
\item Lorsque $X$ est réduit à un point, les $\Q _l$-faisceaux constructibles sont
les $\Q _l$-espaces vectoriels de dimension {\it finie} ;
  \item Ils contiennent les coefficients constants ;
  \item Ils sont stables par les {\it six opérations cohomologiques de Grothendieck},
  à savoir $\otimes$,  $\mathcal{H}om$, $f _*$, $f ^*$, $f _!$ et $f^!$.
\end{enumerate}
Les opérations cohomologiques permettent de changer de variétés.
Afin d'établir l'interprétation cohomologique des fonctions, dites $L$,
associées aux $\Q _l$-faisceaux constructibles,
il procède alors, grâce à la stabilité de la constructibilité,
à une récurrence sur la dimension de $X$ (en effectuant des fibrations).

Nous aimerions avoir une analogie $p$-adique de ceci, i.e.,
construire une catégorie de $\Q _p$-faisceaux ou de complexes de $\Q _p$-faisceaux sur $X$ qui vérifie les trois propriétés ci-dessus.
À cette fin, l'idée de Berthelot est,
en s'inspirant de la caractéristique nulle,
d'élaborer une théorie de \textit{$\D$-modules arithmétiques}.

Plus précisément, soient $\V$ un anneau de valuation discrète complet d'inégales caractéristiques
($0$, $p$) de corps résiduel parfait $k$, $\PP$ un $\V$-schéma formel lisse et $P$ sa fibre spéciale.
Berthelot construit alors
le faisceau des opérateurs différentiels d'ordre infini et de niveau fini sur $\PP$,
$\D ^\dag _{\PP, \, \Q}$ (voir \cite{Be1}), que l'on peut voir
comme le complété faible tensorisé par $\Q$ du faisceau des opérateurs différentiels usuels sur $\PP$, $\D _{\PP}$.
D'où la notion de $\D$-modules arithmétiques sur $\PP$, i.e., de $\D ^\dag _{\PP, \, \Q}$-modules
(toujours à gauche par défaut).
Ensuite, après avoir défini l'holonomie en s'inspirant de la caractéristique nulle (\cite[5]{Beintro2}),
après avoir construit leurs opérations cohomologiques (\cite[2,3,4]{Beintro2}),
il donne les conjectures de base (voir \cite[5.3.6]{Beintro2}) sur la stabilité de l'holonomie.
Si celles-ci étaient exactes, nous pourrions définir les foncteurs
$g _+$, $g ^+$, $g _!$ et $g^!$ au niveau des {\it variétés sur $k$}
de façon analogue
à ceux définis ici dans \ref{def-surhol-rele2} pour les complexes surholonomes.
De plus, les $F$-complexes holonomes (le symbole {\og$F$\fg} signifie que les complexes sont munis
d'une action de Frobenius) constitueraient ainsi une catégorie de coefficients stables par
les cinq opérations cohomologiques de Grothendieck
$\otimes$, $g _+$, $g ^+$, $g _!$ et $g^!$ ainsi que par le foncteur dual, le produit tensoriel externe.
Pour l'instant, seule la stabilité de l'holonomie par produits
tensoriels externes (\cite[5.3.5.(v)]{Beintro2}) et foncteurs duaux (\cite[III.4.4]{virrion}) est obtenue
dans le cas général.
Rappelons que,
lorsque $\PP$ est une courbe,
sa stabilité par foncteur cohomologique local a déjà été validée (\cite[2.3.3]{caro_courbe-nouveau}).

Cependant, afin de disposer d'une bonne cohomologie $p$-adique,
nous proposons  ici, une approche légèrement différente : nous remplaçons l'étude de l'holonomie
par celle de la {\it surholonomie}.
Les liens entre ces deux notions sont les suivants :
un $F$-$\D$-module arithmétique surholonome est holonome. De plus, la réciproque est validée
si et seulement si la conjecture sur la stabilité de l'holonomie par foncteur
cohomologique local est exacte (voir \ref{prop-conjB-eq}).

Décrivons maintenant le contenu de cet article.

Nous donnons dans la première partie
quelques rappels sur les $\D$-modules arithmétiques qui nous seront utiles.
Nous formulons dans une deuxième partie
un critère d'holonomie. Il nous permettra d'établir
qu'un $F$-$\D$-module arithmétique surholonome est holonome.

Nous introduisons dans le troisième chapitre la notion de
$F\text{-}\D ^\dag _{\PP,\,\Q}$-modules surholonomes et de $F$-complexes de
$\D ^\dag _{\PP,\,\Q}$-modules surholonomes.
Cette appellation vient du fait qu'un $F\text{-}\D ^\dag _{\PP,\,\Q}$-module surholonome est holonome.
Soient $f$ : $\PP' \rightarrow \PP$ un morphisme de $\V$-schémas formels lisses.
Nous prouvons
la stabilité de la surholonomie
par foncteurs duaux, par foncteurs cohomologiques locaux,
par images inverses extraordinaires et
images inverses par $f$.
De plus, lorsque $f$ est propre,
nous obtenons sa stabilité par les foncteurs images directes et
images directes extraordinaires par $f$.
La stabilité par produits tensoriels internes et externes reste conjecturale.

Nous construisons dans la quatrième partie,
pour toute $k$-variété $U$ plongeable dans un $\V$-schéma formel propre et lisse,
la catégorie des {\it ($F$-)complexes surholonomes sur $U$}
que l'on note $(F\text{-})D ^\mathrm{b} _\mathrm{surhol} (\D _{U })$.
Nous définissons ensuite le foncteur dual sur $U$ noté
$\DD _U$ : $D ^\mathrm{b} _\mathrm{surhol} (\D _{U })\rightarrow
D ^\mathrm{b} _\mathrm{surhol} (\D _{U })$.
Pour tout morphisme $g$ : $U' \rightarrow U$ de $k$-variétés
(on suppose ici qu'il existe des $\V$-schémas formels propres et lisses
dans lesquels $U$ et $U'$ peuvent se plonger respectivement),
nous construisons, sur les catégories de $F$-complexes surholonomes correspondantes,
les images directes, images directes extraordinaires, images inverses extraordinaires,
images inverses par $g$ que l'on désigne respectivement par
$g _+$, $g _!$, $g ^!$, $g ^+$ (voir \ref{def-surhol-rele2}).
De plus, nous prouvons la stabilité de la surholonomie par celles-ci.

Dans la cinquième partie (et dans la quatrième partie pour le cas simple), nous construisons
pour toute variété $U$ sur $k$,
la catégorie des
{\it ($F$-)$\D$-modules arithmétiques surholonomes sur $U$}, que l'on note
($F$-)$\mathfrak{M} ^+ _{U}$.
On obtient la catégorie
$(F\text{-})D ^\mathrm{b} (\mathfrak{M} ^+ _{U})$, ce qui donne un analogue de
$(F\text{-})D ^\mathrm{b} _\mathrm{surhol} (\D _{U })$.
On définit de plus le {\og coefficient constant sur $U$\fg} noté $\O _U$.

Soit $U $ une variété sur $k$ se plongeant dans un
$\V$-schéma formel propre et lisse.
Dans une sixième partie, nous définissons les fonctions $L$ associées aux
$F$-complexes surholonomes sur $U$ et nous en donnons une formule cohomologique.
Cela correspond à une extension de la formule cohomologique d'\'Etesse et Le Stum
des fonctions $L$ associées aux $F$-isocristaux surconvergents sur $U$ (voir \cite{E-LS}).
La preuve se fait par récurrence sur la dimension de $U$ en se ramenant à la situation géométrique
déjà connue (voir \cite{caro_surcoherent}).
Remarquons que l'hypothèse que $U$ se plonge dans un $\V$-schéma formel propre et lisse est sûrement évitable
(il s'agit d'étendre les opérations cohomologiques via des diagrammes de topos \cite[V.3.4.1]{Becris}),
mais nous ne nous sommes pas intéressé à ce point.

Soient $\PP$ un $\V$-schéma formel séparé et lisse,
$X$ un sous-schéma fermé lisse de $P$ et $T$ un diviseur de $P$ tel que
$T _X := T\cap X$ soit un diviseur de $X$.
Dans \cite[4]{caro_courbe-nouveau}, lorsque $\PP$ est une courbe propre et lisse et $P=X$, nous avions établi
l'holonomie des $F$-isocristaux surconvergents sur $X \setminus T _X$, ce qui avait généralisé
les premiers examples {\it non triviaux} de $F\text{-}\D ^\dag _{\PP,\,\Q}$-modules holonomes sur les courbes
de Berthelot de \cite[5]{Be0}
(le cas trivial est celui où $T$ est vide, i.e., celui où les $F$-isocristaux sont {\it convergents}).
Dans la septième partie de cet article,
nous obtenons les premiers exemples (non triviaux) en dimension supérieure
de $F$-$\D$-modules arithmétiques holonomes.
Plus précisément, nous démontrons que les $F$-isocristaux surconvergents unités
sur $X \setminus T _X$ surconvergents le long de $T _X$
sont surholonomes. Cela implique que le coefficient constant $\O _U$ d'une variété $U$ sur $k$ est surholonome.
De plus, nous vérifions que, pour toute variété $U$ lisse sur $k$, le $F$-$\D _U$-module arithmétique associé
à un $F$-isocristal unité surconvergent sur $U$ (voir \ref{eqcat-gen}) est surholonome.
En particulier, on obtient leur holonomie, ce qui valide une conjecture fondamentale de Berthelot
(voir \cite[5.3.6.D]{Beintro2}).
Il faut noter que l'utilisation de la surholonomie (et surtout de sa stabilité)
est l'élément clé de la preuve.
La stabilité de la surholonomie permet ainsi de résoudre des problèmes sur l'holonomie.
Plus généralement, il est raisonnable de penser que l'on puisse établir de manière analogue
(mais il faut au préalable disposer d'un théorème de réduction semi-stable en dimension supérieure :
voir les travaux de Kedlaya de \cite{kedlaya-semistableI} et \cite{kedlaya-semistableII})
la surholonomie des $F$-isocristaux surconvergents sur les variétés lisses.
Il en résulterait leur holonomie, ce qui homologue la conjecture \cite[5.3.6.D]{Beintro2} de Berthelot.
Par dévissage (voir \cite{caro_devissge_surcoh}), on en déduirait aussitôt
la stabilité par produits tensoriels de la surholonomie.
Les travaux de \cite{caro_log-iso-hol} sont une première étape dans cette direction.

Nous prouvons dans la dernière partie que la conjecture \cite[5.3.6.B]{Beintro2} de Berthelot
sur la stabilité de l'holonomie par image inverse extraordinaire implique que les notions
d'holonomie, de surholonomie, de surcohérence sont identiques.
On en déduit que cette conjecture \cite[5.3.6.B]{Beintro2}
entraîne en particulier
(pour une version plus générale voir \ref{B->A})
la stabilité de l'holonomie par image directe par un morphisme propre de $\V$-schémas formels lisses
(ce qui était conjecturé dans \cite[5.3.6.A]{Beintro2}).

\begin{center}
{\bf Notations}
\end{center}

Tout au long de cet article, nous garderons les notations
suivantes :
les schémas formels seront notés par des lettres calligraphiques ou
gothiques et leur fibre spéciale par les lettres romanes
correspondantes. Si $f\,:\, X' \rightarrow X$ est un morphisme de schémas ou de schémas formels, on note
$d _X$ la dimension de $X$ et $d _{X'/X}$ la dimension relative.
De plus, la lettre $\V$ désignera un anneau de valuation discrète complet,
de corps résiduel parfait $k$ de caractéristique $p>0$, de corps de
fractions $K$ de caractéristique $0$.
On fixe $s\geq 1$ un entier naturel et $F$ sera la puissance
$s$-ième de l'endomorphisme de Frobenius.
Si $\E$ est un faisceau abélien, $\E _\Q$ désignera $\E \otimes _\Z \Q$.
Lorsque cela n'est pas précisé,
les modules sont des modules à gauche,
les {\it variétés} sont des variétés sur $k$, i.e.,
des $k$-schémas séparés et de type fini.
Tous les $k$-schémas seront {\it réduits}.
Lorsque nous considérerons des produits de
$\V$-schémas formels (resp. $k$-schémas), nous omettrons parfois d'indiquer $\Spf \V$ (resp. $\Spec k$).
Lorsque le symbole ensemble vide {\og$\emptyset$\fg} apparaît dans une notation, nous ne l'indiquerons pas non plus
(e.g., dans les opérations cohomologiques de la section \ref{section1}
lorsqu'un diviseur est vide ;
e.g., {\og $\D ^\dag  _{\X,\Q}$\fg} à la place de {\og $\D ^\dag   _{\X} ( \emptyset ) _\Q$ \fg}).

Si $\AA$ est un faisceau d'anneaux,
$D ^\mathrm{b}( \AA )$, $D ^+( \AA )$ et $D ^-( \AA )$ désignent respectivement les catégories dérivées
des complexes de $\AA$-modules à cohomologie bornée,
bornée inférieurement et bornée supérieurement.
Lorsque l'on souhaitera préciser entre droite et gauche, on notera alors $D ^* (
\overset{ ^g}{}\AA )$ ou $D ^* ( \AA \overset{ ^d}{})$, où
$*$ est l'un des symboles $\emptyset$, $+$, $-$, ou $\mathrm{b}$.
De plus, les indices {\og $\mathrm{qc}$\fg}, {\og $\mathrm{coh}$\fg}, {\og $\mathrm{surcoh}$\fg}, {\og $\mathrm{parf}$\fg} et {\og $\mathrm{hol}$\fg}
signifient respectivement {\og quasi-cohérent\fg}, {\og  cohérent\fg}, {\og  surcohérent\fg}
{\og  parfait\fg} et {\og holonome\fg} (voir la section \ref{section1}).

\begin{center}
{\bf Remerciements}
\end{center}
La surholonomie répond à des questions de K. Kedlaya, B. Le Stum et N. Tsuzuki concernant
une extension logique de la notion de surcohérence. Leur intérêt à ce sujet et aussi notamment ceux de
A. Abbes et C. Huyghe-Noot ont été une grande source de motivation.
Je remercie B. Le Stum et A. Virrion pour leur attention très stimulante
relative au critère d'holonomie.

\section{Rappels sur les $\D$-modules arithmétiques}
\label{section1}

Nous avons regroupé dans cette section toutes les définitions et propriétés concernant
la théorie des $\D$-modules arithmétiques que nous utiliserons constamment dans ce manuscrit.
Nous nous sommes donc focalisé sur ce qui sera sans cesse nécessaire.
Ainsi, ce survol ne se prétend pas exhaustif. Le lecteur néophyte pourra aussi se reporter à
\cite{Beintro2} ou pour une version plus récente à \cite{kedlaya-padiccohomology}.

Soient $\X$, $\X'$ deux $\V$-schémas formels lisses,
$f _0\, : \, X' \rightarrow X$ un morphisme de $k$-schémas sur les fibres spéciales,
$T$ un diviseur de $X$ tel que $T' :=f _0 ^{-1} (T)$ soit un diviseur de $X'$.
On note $\Y:=\X \setminus T$, $\Y':=\X' \setminus T'$ les ouverts respectifs de $\X$, $\X'$.
Sauf mention explicite du contraire (e.g., pour les définitions de l'image directe et de l'image inverse extraordinaire),
pour simplifier l'exposé, on supposera que $f _0$ se relève en un morphisme
$f \, :\, \X ' \rightarrow \X$ de $\V$-schémas formels lisses.

\begin{vide}
[Opérateurs différentiels de niveau $m$]
  Berthelot a construit pour tout entier $m\geq 0$ donné le {\og faisceau des opérateurs différentiels de niveau $m$ sur $\X$\fg},
qu'il note $\D ^{(m)} _{\X}$ (voir \cite[2]{Be1}).
Comme $p$ n'est pas inversible sur $\O _\X$, on remarque que le faisceau $\D _\X$ des opérateurs différentiels usuels sur $\X$
n'est pas cohérent. Par contre,
les faisceaux $(\D ^{(m)} _{\X}) _{m\in \N}$ donne une filtration croissante exhaustive de $\D _\X$
par des $\O _\X$-algèbres cohérentes
(voir par exemple \cite[2.2.1.7, 2.2.3.1]{Be1} et remarquer que $\O _\X$ n'a pas d'éléments de $p$-torsion).

Il a aussi construit une $\O _\X$-algèbre qu'il note $\B _{\X} ^{(m)} (T)$ (voir \cite[4.2]{Be1}).
Lorsque $f\in \O _\X$ relève une équation définissant $T$ dans $X$,
$\B _{\X} ^{(m)} (T) \riso \O _{\X} [t]/(f ^{p ^{m+1}} t - p)$, cette dernière étant indépendante à isomorphisme canonique près du choix
du relèvement $f$. Il en déduit ainsi la construction de $\B _{\X} ^{(m)} (T)$ par recollement.

On pose enfin
$\smash{\widehat{\D}} _{\X } ^{(m)}(T):=
\B _{\X} ^{(m)} (T) \widehat{\otimes} _{\O_\X}\D _{\X } ^{(m)}$ que l'on pourrait appelé, s'il fallait absolument lui donner un nom,
{\og complété $p$-adique du faisceau des opérateurs différentiels de niveau $m$ sur $\X$ à coefficients surconvergents le long de $T$\fg}.

\end{vide}

\begin{vide}
  [Complexes quasi-cohérents, définitions]
  \label{comp-qc}
La notion de {\og complexes quasi-cohérents sur les schémas formels\fg} a été introduite par Berthelot (voir \cite{Beintro2}).
Soient $\B$ un faisceau de $\O _{\X}$-algèbres (non nécessairement commutatif), $\E \in D ^- ( \B \overset{ ^d}{} )$,
  $\FF \in D ^- (\overset{ ^g}{} \B  )$.
  On pose $\B _i :=\B / \pi ^{i+1} \B$,
  $\E _i := \E \otimes ^\L _{\B} \B _i$, $\FF _i := \B _i\otimes ^\L _{\B} \FF$,
  $\E \widehat{\otimes} ^\L _{\B} \FF:=
  \R \underset{\underset{i}{\longleftarrow}}{\lim} \E _i \otimes ^\L _{\B _i} \FF _i$.
On définit de même $\E  ^{(\bullet)} \widehat{\otimes} ^\L _{\B  ^{(\bullet)}} \FF ^{(\bullet)}$
lorsque $\B  ^{(\bullet)}$ est un système inductif de $\O _\X$-algèbres,
$\E  ^{(\bullet)}$ et $\FF  ^{(\bullet)}$ sont des complexes de $\B  ^{(\bullet)}$-algèbres.

$\bullet$ Le complexe $\E$ (resp. $\FF$) est dit {\og $\B$-quasi-coherent\fg} si
le morphisme canonique
$\E \rightarrow \E \widehat{\otimes} ^\L _{\B} \B$
(resp. $\FF \rightarrow \B \widehat{\otimes} ^\L _{\B} \FF$)
est un isomorphisme.
Pour $* = g$ ou $*=d$, on désigne par $D ^- _\mathrm{qc} (\overset{ ^*}{} \B  )$
(resp. $D ^\mathrm{b} _\mathrm{qc} (\overset{ ^*}{} \B  )$)
la sous catégorie pleine de
$D ^- (\overset{ ^*}{} \B  )$
(resp. $D ^\mathrm{b} (\overset{ ^*}{} \B  )$)
des complexes quasi-cohérents.

$\bullet$ On remarque que
puisque $\smash{\widehat{\D}} _{\X } ^{(m)}(T)$ est un
$\smash{\widehat{\B}} _{\X} ^{(m)}(T)$-module plat (pour les structures droite ou gauche),
un complexe de $\smash{\widehat{\D}} _{\X } ^{(m)}(T)$-modules à gauche ou à droite
est $\smash{\widehat{\D}} _{\X } ^{(m)}(T)$-quasi-cohérent si et seulement
s'il est
$\smash{\widehat{\B}} _{\X} ^{(m)}(T)$-quasi-cohérent.
De plus, on dispose de l'isomorphisme :
$\smash{\widehat{\B}} _{\X} ^{(m)}(T) \otimes ^\L _{\V} \V /\pi ^{i+1}
\riso
\smash{\widehat{\B}} _{\X} ^{(m)}(T) \otimes  _{\V} \V /\pi ^{i+1}$
(cela découle de \cite[4.3.3.(i)]{Be1}).
Ainsi, un complexe de $\smash{\widehat{\B}} _{\X} ^{(m)}(T)$-modules est
$\smash{\widehat{\B}} _{\X} ^{(m)}(T)$-quasi-cohérent si et seulement s'il est
$\O _{\X} $-quasi-cohérent si et seulement s'il est
$\V$-quasi-cohérent.

$\bullet$ On note $\smash{\widehat{\D}} _{\X} ^{(\bullet)}(T) $ le système inductive $(\smash{\widehat{\D}} _{\X} ^{(m)}(T))_{m \in \N}$.
En localisant deux fois $D ^{\mathrm{b}} ( \smash{\widehat{\D}} _{\X} ^{(\bullet)}(T))$ (ces localisations
remplacent respectivement les foncteurs $-\otimes _\Z \Q$ et {\og limite inductive sur le niveau $m$ \fg}),
on obtient
  $\smash{\underset{^{\longrightarrow}}{LD}} ^{\mathrm{b}} _{\Q }
  ( \smash{\widehat{\D}} _{\X} ^{(\bullet)}(T))$
  (voir \cite[4.2.1 et 4.2.2]{Beintro2}).
Soit $\E ^{(\bullet)}=(\E ^{(m)}) _{m \in \N} \in \smash{\underset{^{\longrightarrow}}{LD}} ^{\mathrm{b}} _{\Q}
  ( \smash{\widehat{\D}} _{\X } ^{(\bullet)}(T))$.
D'après \cite[4.2.3]{Beintro2} (voir aussi \cite[1.1.3]{caro_courbe-nouveau}),
 $\E ^{(\bullet)}$ est dit quasi-cohérent si et seulement si, pour tout $m $,
$\E ^{(m)} $ est $\smash{\widehat{\D}} _{\X } ^{(m)}(T)$-quasi-cohérent.
La sous-catégorie pleine des complexes quasi-cohérents est notée
$\smash{\underset{^{\longrightarrow}}{LD}} ^{\mathrm{b}} _{\Q ,\mathrm{qc}}
  ( \smash{\widehat{\D}} _{\X } ^{(\bullet)}(T))$.
\end{vide}

\begin{vide}
  [Quasi-cohérence et variation du diviseur]
  \label{qc-vardiv}
  Soient $T _1 \subset T_2$ deux diviseurs de $X$. Avec les remarques de \ref{comp-qc},
  on dispose du foncteur oubli :
$\smash[b]{\underset{^{\longrightarrow }}{LD }}  ^\mathrm{b} _{\Q, \mathrm{qc}} (
\smash{\widehat{\D}} _{\X} ^{(\bullet)} (T _2))
\rightarrow
\smash[b]{\underset{^{\longrightarrow }}{LD }}  ^\mathrm{b} _{\Q, \mathrm{qc}} (
\smash{\widehat{\D}} _{\X} ^{(\bullet)}(T _1) )$
que l'on note $oub _{T_1, T _2}$.

Réciproquement, on définit un foncteur
$(\hdag T _2,\,T _1)\,:\,\smash[b]{\underset{^{\longrightarrow }}{LD }}  ^\mathrm{b} _{\Q, \mathrm{qc}} (
\smash{\widehat{\D}} _{\X} ^{(\bullet)} (T _1))
\rightarrow
\smash[b]{\underset{^{\longrightarrow }}{LD }}  ^\mathrm{b} _{\Q, \mathrm{qc}} (
\smash{\widehat{\D}} _{\X} ^{(\bullet)}(T _2) )$
en posant, pour tout
$\E  ^{(\bullet)} \in \smash[b]{\underset{^{\longrightarrow }}{LD }}  ^\mathrm{b} _{\Q, \mathrm{qc}}
( \smash{\widehat{\D}} _{\X} ^{(\bullet)} (T _1))$,
$$(\hdag T _2,\,T _1) (\E  ^{(\bullet)} ) :=
\widehat{\D} ^{(\bullet)} _{\X} ( T_2)  \smash{\widehat{\otimes}} ^\L
_{\widehat{\D} ^{(\bullet)}  _{\X} ( T_1) } \E  ^{(\bullet)},
$$
où
$\widehat{\D} ^{(\bullet)} _{\X} ( T_2)  \smash{\widehat{\otimes}} ^\L
_{\widehat{\D} ^{(\bullet)}  _{\X} ( T_1) } \E  ^{(\bullet)} $ désigne le système inductif
$(\widehat{\D} ^{(m)} _{\X} ( T_2)  \smash{\widehat{\otimes}} ^\L
_{\widehat{\D} ^{(m)}  _{\X} ( T_1) } \E  ^{(m)} ) _{m\in \N}$.

 Lorsque $T _1$ est vide, on note simplement $oub _{T_2}$ et $(\hdag T_2)$.
On dispose pour ces deux foncteurs des formules de transitivité en les diviseurs évidentes.
D'après \cite[1.1.8]{caro_courbe-nouveau}, on bénéficie aussi de l'isomorphisme
 $(\hdag T _2,\,T _1) \riso (\hdag T _2) \circ oub _{T_1}$.
 Cela implique que les foncteurs
 de la forme $oub _{T_1}$ sont pleinement fidèles (voir \cite[1.1.8]{caro_courbe-nouveau}).
 Il en est donc de même des ceux de la forme $oub _{T_1, T _2}$.
 On pourra ainsi sans aucune ambiguïté omettre d'indiquer les inclusions
 canoniques $oub _{T_1, T _2} \,:\,\smash[b]{\underset{^{\longrightarrow }}{LD }}  ^\mathrm{b} _{\Q, \mathrm{qc}} (
\smash{\widehat{\D}} _{\X} ^{(\bullet)} (T _2))
\subset
\smash[b]{\underset{^{\longrightarrow }}{LD }}  ^\mathrm{b} _{\Q, \mathrm{qc}} (
\smash{\widehat{\D}} _{\X} ^{(\bullet)}(T _1) )$
et identifier les deux foncteurs
$(\hdag T _2,\,T _1) $ et $(\hdag T _2) $.
\end{vide}

\begin{vide}
  [Produits tensoriels internes]
Soient $\E ^{(\bullet)}  , \, \FF ^{(\bullet)}
\in \smash[b]{\underset{^{\longrightarrow }}{LD }}  ^\mathrm{b} _{\Q, \mathrm{qc}}
( \smash{\widehat{\D}} _{\X} ^{(\bullet)} (T ))$.
On définit le {\og produit tensoriel interne\fg} de
$\E ^{(\bullet)}$ et $\FF ^{(\bullet)}$ en posant dans
$\smash[b]{\underset{^{\longrightarrow }}{LD }}  ^\mathrm{b} _{\Q, \mathrm{qc}}
( \smash{\widehat{\D}} _{\X} ^{(\bullet)} (T ))$ :
\begin{equation}
  \label{prod-tens-int}
\E ^{(\bullet)}
 \smash{\overset{\L}{\otimes}}^{\dag} _{\O _{\X } ( \hdag T ) _{\Q}}
  \FF ^{(\bullet)}
  :=
  \E ^{(\bullet)}  \smash{\widehat{\otimes}} ^\L _{\smash{\widehat{\B}} _{\X} ^{(\bullet)} (T )} \FF ^{(\bullet)} ,
\end{equation}
où
$\E ^{(\bullet)}  \smash{\widehat{\otimes}} ^\L _{\smash{\widehat{\B}} _{\X} ^{(\bullet)} (T )} \FF ^{(\bullet)}$
désigne le système inductif
$(\E ^{(m)}  \smash{\widehat{\otimes}} ^\L _{\smash{\widehat{\B}} _{\X} ^{(m)} (T )} \FF ^{(m)} ) _{m \in  \N }$.

Soient $T _1 \subset T _2$ deux diviseurs et
$\E ^{(\bullet)}
\in \smash[b]{\underset{^{\longrightarrow }}{LD }}  ^\mathrm{b} _{\Q, \mathrm{qc}}
( \smash{\widehat{\D}} _{\X} ^{(\bullet)} (T _1))$. Pour simplifier les notations, on écrira abusivement
$\O _{\X } ( \hdag T _2) _{\Q}
 \smash{\overset{\L}{\otimes}}^{\dag} _{\O _{\X } ( \hdag T_1 ) _{\Q}}
  \E ^{(\bullet)}$
pour
$\smash{\widehat{\B}} _{\X} ^{(\bullet)} (T _2 )  ^{(\bullet)}
 \smash{\overset{\L}{\otimes}}^{\dag} _{\O _{\X } ( \hdag T _1) _{\Q}}
  \E ^{(\bullet)}$.
D'après \cite[1.1.7.1]{caro_courbe-nouveau}, on dispose alors de l'isomorphisme canonique :
\begin{equation}
\label{courbe-1171}
  \O _{\X } ( \hdag T _2) _{\Q}
 \smash{\overset{\L}{\otimes}}^{\dag} _{\O _{\X } ( \hdag T_1 ) _{\Q}}
  \E ^{(\bullet)} \riso (\hdag T _2,\,T _1) (\E ^{(\bullet)}).
\end{equation}

\end{vide}

\begin{vide}
[Produits tensoriels externes]

Soient $\X _1$, $\X_2$ deux $\V$-schémas formels lisses,
$\E _1 ^{(\bullet)}
\in \smash[b]{\underset{^{\longrightarrow }}{LD }}  ^\mathrm{b} _{\Q, \mathrm{qc}}
( \smash{\widehat{\D}} _{\X _1} ^{(\bullet)})$,
$\E _2 ^{(\bullet)} \in \smash[b]{\underset{^{\longrightarrow }}{LD }}  ^\mathrm{b} _{\Q, \mathrm{qc}}
( \smash{\widehat{\D}} _{\X _2} ^{(\bullet)})$,
$p _1 \,:\, \X _1\times _\S \X _2\rightarrow \X_1$
et
$p _2 \,:\, \X _1\times _\S \X _2\rightarrow \X_2$
les projections canoniques.
Le {\og produit tensoriel externe\fg} de $\E _1  ^{(\bullet)} $ et $\E  ^{(\bullet)} _2$ est défini en posant (voir \cite[4.3.5]{Beintro2}) :
\begin{equation}
  \label{def-boxtimes}
  \E  ^{(\bullet)} _1\smash{\overset{\L}{\boxtimes}} _{\O _\S}^{\dag} \E  ^{(\bullet)} _2 :=
p _1 ^! (\E  ^{(\bullet)} _1 ) \smash{\overset{\L}{\otimes}}^{\dag}_{\O _{\X _1\times _\S \X_2,\, \Q}} p _2^!(\E  ^{(\bullet)} _2)
[- d _{\X _1} -d _{\X _2}].
\end{equation}

\end{vide}

\begin{vide}
  [Images directes et images inverses extraordinaires]

Le faisceau $\smash{\widehat{\D}} ^{(m)} _{\X ' \rightarrow \X} (T',T):=
\B ^{(m)} _{\X'} (T') \widehat{\otimes} _{\O _{\X'}} f ^* \smash{\widehat{\D}} _{\X} ^{(m)} (T )$
est muni par fonctorialité d'une structure de
$(\smash{\widehat{\D}} _{\X'} ^{(m)} (T '), f^{-1} (\smash{\widehat{\D}} _{\X} ^{(m)} (T )))$-bimodule.

On définit
le foncteur
$f _{T} ^! \,:\,\smash[b]{\underset{^{\longrightarrow }}{LD }}  ^\mathrm{b} _{\Q, \mathrm{qc}}
( \smash{\widehat{\D}} _{\X} ^{(\bullet)} (T ))
\rightarrow
\smash[b]{\underset{^{\longrightarrow }}{LD }}  ^\mathrm{b} _{\Q, \mathrm{qc}}
( \smash{\widehat{\D}} _{\X'} ^{(\bullet)} (T '))$
{\og image inverse extraordinaire par $f$ à coefficients surconvergent le long de $T$\fg}
en posant,
pour tout $\E ^{(\bullet)} \in \smash[b]{\underset{^{\longrightarrow }}{LD }}  ^\mathrm{b} _{\Q, \mathrm{qc}}
( \smash{\widehat{\D}} _{\X} ^{(\bullet)} (T ))$, :
\begin{equation}
\label{iminvqc}
  f _{T} ^! ( \E ^{(\bullet)}) :=
\smash{\widehat{\D}} ^{(\bullet)} _{\X ' \rightarrow \X} (T',T)
\smash{\widehat{\otimes}} ^\L _{f ^{\text{-}1} \smash{\widehat{\D}} ^{(\bullet)} _\X (T)}
f ^{\text{-}1} \E ^{(\bullet)} [ d_{\X ' /\X}].
\end{equation}

Posons
$\smash{\widehat{\D}} ^{(m)} _{\X \leftarrow \X'} (T,T'):=
\B ^{(m)} _{\X'} (T') \widehat{\otimes} _{\O _{\X'}}
(\omega _{\X'} \otimes _{\O _{\X'}} f _\mathrm{g} ^* (\smash{\widehat{\D}} _{\X} ^{(m)} (T )\otimes _{\O _\X} \omega _{\X} ^{-1})$,
l'indice {\og $\mathrm{g}$\fg} signifiant que l'on a choisit pour calculer l'image inverse la structure gauche de
$\smash{\widehat{\D}} _{\X} ^{(m)} (T )$-module à gauche.
Par fonctorialité, $\smash{\widehat{\D}} ^{(m)} _{\X \leftarrow \X'} (T,T')$ est muni d'une structure de
$(f^{-1} (\smash{\widehat{\D}} _{\X} ^{(m)} (T )), \smash{\widehat{\D}} _{\X'} ^{(m)} (T '))$-bimodule.
On définit
le foncteur
$f _{T,+}  \,:\,\smash[b]{\underset{^{\longrightarrow }}{LD }}  ^\mathrm{b} _{\Q, \mathrm{qc}}
( \smash{\widehat{\D}} _{\X'} ^{(\bullet)} (T '))
\rightarrow
\smash[b]{\underset{^{\longrightarrow }}{LD }}  ^\mathrm{b} _{\Q, \mathrm{qc}}
( \smash{\widehat{\D}} _{\X} ^{(\bullet)} (T ))$
{\og image directe par $f$ à coefficients surconvergent le long de $T$\fg}
en posant,
pour tout
$\E ^{\prime (\bullet)} \in \smash[b]{\underset{^{\longrightarrow }}{LD }}  ^\mathrm{b} _{\Q, \mathrm{qc}}
( \smash{\widehat{\D}} _{\X'} ^{(\bullet)} (T '))$ :
\begin{equation}
\label{imdirqc}
f _{T,+}  ( \E ^{\prime (\bullet)} ):= \R f _* (
\smash{\widehat{\D}} ^{(\bullet)} _{\X \leftarrow \X'} (T,T')
\smash{\widehat{\otimes}} ^\L _{\smash{\widehat{\D}} ^{(m)} _{\X '} (T')}
\E ^{\prime (\bullet)}).
\end{equation}

D'après \cite[1.1.9, 1.1.10]{caro_courbe-nouveau},
on dispose des isomorphismes canoniques de foncteurs :
\begin{equation}
  \label{fT=f}
  oub _{T} \circ f _{T,+} \riso f _+  \circ oub _{T'}
  \hspace{1cm}
  \text{ et }
  \hspace{1cm}
  oub _{T'} \circ f ^! _{T} \riso f ^! \circ oub _T .
\end{equation}
En omettant d'indiquer les inclusions de la forme
$oub _{T} \,:\,\smash[b]{\underset{^{\longrightarrow }}{LD }}  ^\mathrm{b} _{\Q, \mathrm{qc}} (
\smash{\widehat{\D}} _{\X} ^{(\bullet)} (T ))
\subset
\smash[b]{\underset{^{\longrightarrow }}{LD }}  ^\mathrm{b} _{\Q, \mathrm{qc}} (
\smash{\widehat{\D}} _{\X} ^{(\bullet)} )$,
on pourra ainsi par abus de notations écrire
$f ^!$ pour $f _T ^!$ (resp. $f _+$ pour $f _{T+}$).

Lorsque $f _0$ ne se relève pas, les bimodules
$\smash{\widehat{\D}} ^{(m)} _{\X ' \rightarrow \X} (T',T)$ et
$\smash{\widehat{\D}} ^{(m)} _{\X \leftarrow \X'} (T,T')$
peuvent néanmoins être défini par recollement (e.g., voir \cite[2.1.6]{Be2}).
On construit ainsi de manière identique les foncteurs image inverse extraordinaire et image directe par $f _0$,
que l'on note
respectivement par $f _0 ^!$ et $f _{0+}$.
On obtient en particulier la notion d'image inverse par $F _X$, où $F _X$ désigne la puissance $s$-ième de l'endomorphisme de Frobenius absolu de $X$.
Ce foncteur se note simplement $F ^*$. Nous verrons que $F^*$ commute à l'image directe et l'image inverse extraordinaire
(voir \ref{trans-f!pref+} et \ref{commFrobf+}).

\end{vide}

\begin{vide}
  [Commutation à Frobenius et transitivité de l'image inverse extraordinaire]
\label{trans-f!pref+}
  Soient $f$ : $\X '\rightarrow \X$,
  $g$ : $\X'' \rightarrow \X'$ et $h$ : $\X'''\rightarrow \X ''$ trois morphismes
de $\V$-schémas formels lisses.

On dispose, pour tout
$\E ^{(\bullet)} \in \smash[b]{\underset{^{\longrightarrow }}{LD }}  ^\mathrm{b} _{\Q, \mathrm{qc}} (
\smash{\widehat{\D}} _{\X} ^{(\bullet)})$,
 de l'isomorphisme canonique dit de {\og transitivité de l'image inverse extraordinaire\fg},
$g  ^! \circ f  ^! \riso (f\circ g ) ^! $, vérifiant la condition
d'associativité : les deux isomorphismes
$h ^!  g  ^! f  ^! \riso h ^!  (f\circ g ) ^! \riso (f\circ g \circ h) ^! $
et
$h ^!  g  ^!f  ^! \riso (g\circ h ) ^! f ^! \riso (f\circ g \circ h) ^! $
sont identiques (voir \cite[1.2.2]{caro_courbe-nouveau}).

En particulier, en désignant par $F$ les morphismes absolus de Frobenius de $X$ ou $Y$,
on dispose des isomorphismes canoniques de commutation à Frobenius :
$F ^* \circ f ^! (\E ^{(\bullet)}) \riso f ^! \circ F^* (\E ^{(\bullet)})$ (car $F ^* =F ^!$).
Les isomorphismes de transitivité de l'image inverse extraordinaire sont compatibles à Frobenius.

\end{vide}

\begin{vide}
[Passage à la limite sur le niveau et cohérence]
\label{lim-coh}
Berthelot défini par passage à la limite
$\O _{\X } (\hdag T) _{\Q}:=\underset{\underset{m}{\longrightarrow}}{\lim}\,
\widehat{\B} _{\X } ^{(m)} ( T)_{\Q}$ le {\og faisceau des fonctions sur $\X$
à singularités surconvergentes le long de $T$\fg} (\cite[4.2.4]{Be1}).
Il construit de même $\D ^{\dag} _{\X } (\hdag T) _{\Q}:=\underset{\underset{m}{\longrightarrow}}{\lim}\,
\smash{\widehat{\D}} _{\X } ^{(m)} ( T)_{\Q}$ et le nomme {\og faisceau des opérateurs différentiels de niveau fini,
à singularités surconvergentes le long de $T$\fg} (\cite[4.2.5]{Be1}).

En passant à la limite sur le niveau puis en tensorisant par $\Q$ sur $\Z$, on obtient
le foncteur noté abusivement (pour simplifier les notations)
$\underset{\longrightarrow}{\lim}$ :
$\smash{\underset{^{\longrightarrow}}{LD}} ^{\mathrm{b}} _{\Q ,\mathrm{qc}}
(\smash{\widehat{\D}} _{\X} ^{(\bullet)}(T))
\rightarrow
D  ( \D ^\dag _{\X} (\hdag T) _{\Q} )$ (e.g., voir \cite[4.2.2]{Beintro2}).
Celui-ci induit une équivalence de
catégories entre $D ^\mathrm{b} _\mathrm{coh} ( \D ^\dag _{\X} (\hdag T) _{\Q} )$
et une sous-catégorie pleine de $\smash{\underset{^{\longrightarrow}}{LD}} ^{\mathrm{b}} _{\Q ,\mathrm{qc}}
(\smash{\widehat{\D}} _{\X} ^{(\bullet)}(T))$, notée
$\smash{\underset{^{\longrightarrow}}{LD}} ^{\mathrm{b}} _{\Q ,\mathrm{coh}}
(\smash{\widehat{\D}} _{\X} ^{(\bullet)}(T))$ (voir \cite[4.2.4]{Beintro2}).
Par abus de notations, il nous arrivera d'omettre le foncteur $\underset{\longrightarrow}{\lim}$, i.e., d'identifier
$\smash{\underset{^{\longrightarrow}}{LD}} ^{\mathrm{b}} _{\Q ,\mathrm{coh}}
(\smash{\widehat{\D}} _{\X} ^{(\bullet)}(T))$
avec
$D ^\mathrm{b} _\mathrm{coh} ( \D ^\dag _{\X} (\hdag T) _{\Q} )$.
De cette manière
$D ^\mathrm{b} _\mathrm{coh} ( \smash{\D} ^\dag _{\X} (\hdag T) _{\Q} )$
est une sous-catégorie pleine de
$\smash{\underset{^{\longrightarrow}}{LD}} ^{\mathrm{b}} _{\Q ,\mathrm{qc}}
(\smash{\widehat{\D}} _{\X} ^{(\bullet)}(T))$.

\end{vide}

\begin{vide}
  [Cohérence : comment se ramener au cas sans diviseur]
\label{4312Be1}
  Nous utiliserons énormément le résultat fondamental suivant de Berthelot (voir \cite[4.3.12]{Be1}) :
  un $\D ^{\dag} _{\X } (\hdag T) _{\Q}$-module cohérent
   est nul si et seulement si sa restriction sur $\Y$ est nulle.
   Cela implique par exemple qu'un morphisme de $\D ^{\dag} _{\X } (\hdag T) _{\Q}$-modules cohérents
   est injectif (resp. surjectif) si et seulement si sa restriction sur $\Y$ est injectif (resp. surjectif).
   On en déduit aussi qu'un morphisme de
   $D ^\mathrm{b} _\mathrm{coh} ( \smash{\D} ^\dag _{\X} (\hdag T) _{\Q} )$ est
   un isomorphisme si et seulement si sa restriction sur $\Y$ est un isomorphisme.
   La philosophie est qu'un $\D ^{\dag} _{\X } (\hdag T) _{\Q}$-module cohérent {\og vit essentiellement sur $\Y$\fg}.
\end{vide}

\begin{vide}
  [Opérations cohomologiques de complexes cohérents]
Considérons les bimodules suivants $\D ^{\dag} _{\X '\rightarrow \X } (\hdag T' ,T) _{\Q}:=\underset{\underset{m}{\longrightarrow}}{\lim}\,
\smash{\widehat{\D}} _{\X '\rightarrow \X } ^{(m)} ( T' ,T)_{\Q}$,
$\D ^{\dag} _{\X \leftarrow \X'} (\hdag T ,T ')_{\Q}:=
\underset{\underset{m}{\longrightarrow}}{\lim}
\smash{\widehat{\D}} _{\X \leftarrow \X '} ^{(m)} ( T ,T')_{\Q}$.
Soient $T \subset T _2$ un diviseur de $X$,
$\E ^{(\bullet)}
\in \smash[b]{\underset{^{\longrightarrow }}{LD }}  ^\mathrm{b} _{\Q, \mathrm{coh}}
( \smash{\widehat{\D}} _{\X} ^{(\bullet)} (T ))$,
$\FF ^{(\bullet)}
\in \smash[b]{\underset{^{\longrightarrow }}{LD }}  ^\mathrm{b} _{\Q, \mathrm{coh}}
( \smash{\widehat{\D}} _{\X} ^{(\bullet)} (T ))$,
$\E ^{\prime (\bullet)}
\in \smash[b]{\underset{^{\longrightarrow }}{LD }}  ^\mathrm{b} _{\Q, \mathrm{coh}}
( \smash{\widehat{\D}} _{\X'} ^{(\bullet)} (T' ))$,
$\E  := \underset{\longrightarrow}{\lim} \E ^{(\bullet)} $,
$\FF  := \underset{\longrightarrow}{\lim} \FF ^{(\bullet)} $,
$\E '  := \underset{\longrightarrow}{\lim} \E ^{\prime (\bullet)} $.
On définit alors les foncteurs suivants en posant :
\begin{gather}
\label{def-hdagcoh}
  (\hdag T _2 , T ) ( \E ) := \D ^{\dag} _{\X } (\hdag T _2 ) _{\Q}
\otimes _{\D ^{\dag} _{\X } (\hdag T  ) _{\Q}} \E,
\\
f ^! _{T} (\E ):=\D ^{\dag} _{\X '\rightarrow \X }  ( \hdag T' , T ) _{\Q}
\otimes ^{\L} _{ f ^{-1} \D ^{\dag} _{\X } (\hdag T ) _{\Q}} f ^{-1} \E [d _{\X '/\X} ],
\\
f  _{T +}( \E '):=
 \R f_* (\D ^{\dag} _{\X \leftarrow \X '}  ( \hdag T , T ') _{\Q}
\otimes ^{\L} _{ \D ^{\dag} _{\X '} (\hdag T ') _{\Q}} \E '),
\end{gather}
où nous avons omis d'inscrire le symbole {\og $\L$\fg} dans le premier produit tensoriel
car l'homomorphisme $\D ^{\dag} _{\X } (\hdag T  ) _{\Q} \rightarrow \D ^{\dag} _{\X } (\hdag T _2 ) _{\Q}$ est plat
(voir \cite[4.3.10-11]{Be1}).
On dispose des isomorphismes
$(\hdag T _2 , T ) ( \E ) \riso \underset{\longrightarrow}{\lim} (\hdag T _2 , T ) ( \E ^{(\bullet)} )$,
$f^! _T ( \E ) \riso \underset{\longrightarrow}{\lim} f^! _T ( \E ^{(\bullet)} )$,
$f_{T+} ( \E ') \riso \underset{\longrightarrow}{\lim} f_{T+} ( \E ^{\prime (\bullet)} )$.
Les notations relatives à ces foncteurs sont donc compatibles.
De plus, on pose (voir \ref{prod-tens-int}, \ref{def-boxtimes}) :
\begin{gather}
  \label{otimes-coh}
\E
 \smash{\overset{\L}{\otimes}}^{\dag} _{\O _{\X } ( \hdag T ) _{\Q}}
  \FF
  :=
  \underset{\longrightarrow}{\lim}  \E ^{(\bullet)}
 \smash{\overset{\L}{\otimes}}^{\dag} _{\O _{\X } ( \hdag T ) _{\Q}}
  \FF ^{(\bullet)} ,\\
    \E   \smash{\overset{\L}{\boxtimes}} _{\O _\S}^{\dag} \FF :=
\underset{\longrightarrow}{\lim}
\E  ^{(\bullet)} \smash{\overset{\L}{\boxtimes}} _{\O _\S}^{\dag} \FF  ^{(\bullet)}.
\end{gather}
\medskip

$\bullet$ Pour tout $\G \in D _{\mathrm{parf}} ( \D ^{\dag} _{\X ,\Q} (\hdag T ))$,
Virrion définit son {\og foncteur dual $\D ^{\dag} _{\X ,\Q} (\hdag T )$-linéaire\fg} (voir \cite[I.3.2]{virrion})
en posant :
\begin{equation}
\label{def-dualdiv}
\DD _{T} ( \G) := \R \mathcal{H}om _{\D ^{\dag} _{\X } (\hdag T ) _{\Q}}
( \G , \D ^{\dag} _{\X } (\hdag T ) _{\Q} \otimes _{\O _{\X , \Q}} \omega _{\X ,\Q} ^{-1} [d _{\X}]).
\end{equation}
Elle a aussi vérifié (voir \cite[II.3.5]{virrion}) que
l'on dispose de l'{\og isomorphisme de bidualité\fg} compatible à Frobenius :
\begin{equation}
  \label{bidual}
  \DD _T \circ \DD _T (\G) \riso \G.
\end{equation}
\medskip

$\bullet$ Rappelons que Noot-Huyghe (\cite{huyghe_finitude_coho}) a prouvé que le faisceau
$\D ^{\dag} _{\X } (\hdag T ) _{\Q}$
est de dimension cohomologique finie. Comme $\D ^{\dag} _{\X } (\hdag T ) _{\Q}$ est aussi cohérent (\cite[5.4]{Be1}),
on obtient
$D _{\mathrm{parf}} ( \D ^{\dag} _{\X } (\hdag T ) _{\Q})=
D _{\mathrm{coh}} ^{\mathrm{b}}( \D ^{\dag} _{\X } (\hdag T ) _{\Q})$.
Ainsi, le foncteur $\DD _{T} $ préserve la cohérence.
Par contre, le foncteur $f _{T+}$ (resp. $f ^! _T$) ne préserve pas la cohérence
lorsque $f $ est une immersion ouverte (resp. immersion fermée) différente de l'identité.
Berthelot (voir \cite{Beintro2}) a néanmoins vérifié la stabilité de la cohérence par le foncteur
$f _{T+}$ (resp. $f ^! _T$) lorsque $f$ est un morphisme propre (resp. un morphisme lisse).
\medskip

$\bullet$ Lorsque cela aura un sens, on notera alors respectivement
$ f ^+ _{ T} : = \DD_{T'} \circ f ^! _{T} \circ \DD_{T}$ et
$f _{ T , !} : =\DD_{T} \circ f  _{T   +} \circ \DD_{T '}$
{\og l'image inverse\fg} et {\og l'image directe extraordinaire\fg}.
\end{vide}

\begin{vide}
[Isomorphisme de dualité relative]
\label{dual-relasection}
Dans toute cette section concernant l'isomorphisme de dualité relative, $f$ est propre.
Virrion construit dans \cite{Vir04}
un isomorphisme de commutation entre les foncteurs duaux respectifs et
le foncteur image directe par $f$. Un tel isomorphisme est appelé {\og isomorphisme de dualité relative\fg}.
En reprenant ses constructions (notamment, celui très technique du morphisme trace),
nous avons vérifié dans \cite[1.2.7]{caro_courbe-nouveau} que cet isomorphisme de dualité relative
s'étend de la manière suivante :
pour tout $\E '\in D _{\mathrm{coh}} ^{\mathrm{b}} ( \D ^{\dag} _{\X '} (\hdag T')_{\Q} )$,
il existe alors dans $D _{\mathrm{coh}} ^{\mathrm{b}} (\D ^{\dag} _{\X }(\hdag T)_{\Q}  )$
l'isomorphisme canonique :
\begin{equation}
\label{iso-dualrel}
f _{T,+} \circ \DD _{T '} (\E ')\riso \DD _{T}\circ f _{T , +} (\E ' ).
\end{equation}
Lorsque le diviseur est vide, on retrouve un cas particulier du contexte de \cite{Vir04}.
Lorsque le diviseur $T$ n'est pas vide, il n'est pas clair que les catégories utilisées dans \cite{Vir04}
correspondent aux complexes à cohomologie bornée et cohérente. C'est précisément ce qui nous
avait conduit à vérifier \ref{iso-dualrel}
pour $\E '\in D _{\mathrm{coh}} ^{\mathrm{b}} ( \D ^{\dag} _{\X '} (\hdag T')_{\Q} )$.

En outre, lorsque $f $ est une immersion fermée, nous avons vérifié dans \cite[2.4.3]{caro-frobdualrel}
que \ref{iso-dualrel} commute à Frobenius (plus précisément, aux images inverses par Frobenius $F ^*$)

Le théorème de dualité relative implique le fait suivant :
pour tous $\E' \in  D _{\mathrm{coh}} ^{\mathrm{b}} ( \D ^{\dag} _{\X '} (\hdag T')_{\Q} )$,
$\E \in D _{\mathrm{coh}} ^{\mathrm{b}} ( \D ^{\dag} _{\X } (\hdag T)_{\Q} )$,
on dispose alors de l'isomorphisme canonique d'adjonction
fonctoriel en $\E$ et $\E'$ :
\begin{gather}
\label{f+adjf!}
\mathrm{adj} _{f,T} : \mathrm{Hom} _{ \D ^{\dag} _{\X  ,\,\Q} ( \hdag T) }( f _{T, +} (\E'),\E )
\riso
\mathrm{Hom} _{ \D ^{\dag} _{\X ',\,\Q} ( \hdag T') }(\E'  ,f _{T } ^! (\E) ),
\end{gather}
où $\mathrm{Hom} _{ \D ^{\dag} _{\X  ,\,\Q}} (-,-) :=
H ^0 \circ \R \mathrm{Hom} _{ D( \D ^{\dag} _{\X  ,\,\Q})} (-,-) =
\mathrm{Hom} _{D( \D ^{\dag} _{\X  ,\,\Q})} (-,-)$.

D'après \cite[1.2.15]{caro_courbe-nouveau} (voir aussi \cite[1.2.8]{caro-construction} avec
la remarque \cite[1.2.9]{caro-construction}),
cet isomorphisme d'adjonction est transitif pour la composition des morphismes propres.

\end{vide}

\begin{vide}
  [Commutation à Frobenius de l'image directe par un morphisme propre d'un complexe cohérent]
  \label{commFrobf+}
  Berthelot a construit (voir \cite[4.3.9.1]{Beintro2} ou \cite[3.4.4]{Be2}),
  pour tout
  $\E ^{(\bullet)} \in \smash[b]{\underset{^{\longrightarrow }}{LD }}  ^\mathrm{b} _{\Q, \mathrm{qc}} (
  \smash{\widehat{\D}} _{\X} ^{(\bullet)} )$,
  l'isomorphisme canonique de commutation de l'image directe par $f $ à Frobenius :
  $f _+ F^* (\E ^{(\bullet)} )\riso F ^* f _+ (\E ^{(\bullet)} ) $.
\medskip

  Lorsque $f$ est propre, d'après \cite[1.2.12.1]{caro_courbe-nouveau}, pour tout
  $\E \in D _{\mathrm{coh}} ^{\mathrm{b}} ( \D ^{\dag} _{\X } (\hdag T)_{\Q} )$, on dispose d'une seconde construction
  de l'isomorphisme
  $f _{T+} F^* (\E )\riso F ^* f _{T+} (\E  ) $.
  Ce dernier isomorphisme est caractérisé par le fait que,
  pour tous $\E' \in  D _{\mathrm{coh}} ^{\mathrm{b}} ( \D ^{\dag} _{\X '} (\hdag T')_{\Q} )$,
$\E \in D _{\mathrm{coh}} ^{\mathrm{b}} ( \D ^{\dag} _{\X } (\hdag T)_{\Q} )$,
 le morphisme d'adjonction de \ref{f+adjf!}
  \begin{gather}
\notag
\mathrm{adj} _{f,T} : \mathrm{Hom} _{ \D ^{\dag} _{\X  ,\,\Q} ( \hdag T) }( f _{T, +} (\E'),\E )
\riso
\mathrm{Hom} _{ \D ^{\dag} _{\X ',\,\Q} ( \hdag T') }(\E'  ,f _{T } ^! (\E) ),
\end{gather}
est alors compatible aux isomorphismes de commutation à Frobenius (voir \cite[1.2.21]{caro_courbe-nouveau}).
De plus, les morphismes d'adjonction entre l'image directe et l'image inverse
extraordinaire par un morphisme propre (i.e., lorsqu'ils existent, ceux de la forme
$f _{T, +} f _{T } ^! (\E) \rightarrow \E$ ou $\E '\rightarrow f _{T } ^! f _{T, +}  (\E') $)
sont alors, pour cette seconde construction,
aussi compatibles à Frobenius (voir \cite[1.2.13, 1.2.14]{caro_courbe-nouveau}).
En outre, d'après \cite[1.2.11]{caro-construction}, ces morphismes d'adjonction sont transitifs pour la
composition de morphismes propres (voir aussi \cite[1.2.8]{caro-construction} avec
la remarque \cite[1.2.9]{caro-construction}).

  Il n'est pas évident que ces deux constructions coïncident.
  Cependant, lorsque $f $ est une immersion fermée,
  cela a déjà été vérifié (voir \cite[2.5.4]{caro-frobdualrel}).
En outre, toujours lorsque $f$ est une immersion fermée, nous avons établi (voir \cite[2.4.3]{caro-frobdualrel})
que,
pour tout $\E '\in D _{\mathrm{coh}} ^{\mathrm{b}} ( \D ^{\dag} _{\X '} (\hdag T')_{\Q} )$,
l'isomorphisme de dualité relative
\begin{equation}
\notag
f _{T,+} \circ \DD _{T '} (\E ')\riso \DD _{T}\circ f _{T , +} (\E ' )
\end{equation}
de \ref{iso-dualrel} est compatible à Frobenius.

 Nous n'utiliserons dans cet article que l'isomorphisme $f _{T+} F^* (\E )\riso F ^* f _{T+} (\E  ) $
  construit dans \cite[1.2.12.1]{caro_courbe-nouveau}.

  Concernant la transitivité de l'image directe par un morphisme propre et sa compatibilité à Frobenius,
  cela sera traité de manière plus approfondi dans la section \ref{compatsurhol} dans le cas des complexes surholonomes.

\end{vide}

\begin{vide}
  [Structure de Frobenius, holonomie, critère d'holonomie de Virrion]
\label{cri-hol}
Rappelons la convention suivante de Berthelot (\cite[5.1]{Beintro2}) :
un $F$-$\D^{\dag} _{\X} ( \hdag T) _{\Q}$-\textit{module} est la donnée
d'un $\D^{\dag} _{\X} ( \hdag T) _{\Q}$-module $\E$ et d'un isomorphisme $\D^{\dag} _{\X} ( \hdag T) _{\Q}$-linéaire
$\Phi$ : $ \E \riso F^* \E$. Les morphismes de $F$-$\D^{\dag} _{\X} ( \hdag T) _{\Q}$-modules sont les morphismes
$\D^{\dag} _{\X} ( \hdag T) _{\Q}$-linéaires commutant à l'action de Frobenius.
De même, nous appellerons
$F$-$\D^{\dag} _{\X} ( \hdag T) _{\Q}$-\textit{complexe} la donnée d'un complexe
$\E\in D ^{\mathrm{b}} _{\mathrm{coh}}(\D^{\dag} _{\X} ( \hdag T) _{\Q})$
et d'un isomorphisme $\Phi$ : $ \E \riso F^* \E $ dans $D ^{\mathrm{b}} _{\mathrm{coh}}(\D^{\dag} _{\X} ( \hdag T) _{\Q})$.
On notera cette catégorie
$F$-$D ^{\mathrm{b}} _{\mathrm{coh}} (  \D _{\X , \Q } ^{\dag} (\hdag T))$.
Par abus de notations, on pourra omettre d'indiquer la structure de Frobenius $\Phi$
et ainsi écrire $\E$ à la place $(\E,\Phi)$.
\medskip

$\bullet$
Soit $\E$ un $F$-$\D ^{\dag} _{\X,\Q}$-module cohérent. Berthelot définit dans \cite[5.2.7]{Beintro2}
la variété caractéristique de $\E$, que l'on notera $\mathrm{Car} (\E)$. En outre, d'après \cite[5.3.4]{Beintro2},
l'inégalité de Bernstein est vérifiée : si $\E$ est non nul
alors, pour tout point $x$ du support de $\E$, on a
\begin{equation}
\mathrm{dim}  _x \mathrm{Car} (\E) \geq \mathrm{dim} _x X.
\end{equation}

Conformément à Virrion (voir \cite[III.1.2]{virrion}),
on définit
la {\og dimension de $\E$\fg} en posant $\dim \E := \sup _x ( \dim _x \mathrm{Car}(\E))$.
De plus, la {\og codimension de $\E$\fg} est par définition $\mathrm{codim} (\E):= 2 d_X - \dim (\E)$
(c'est donc la codimension de la variété caractéristique $\mathrm{Car}(\E)$ associée à $\E$ dans le fibré cotangent).
Ces deux définitions s'étendent naturellement aux $F\text{-}\D ^\dag _{\PP, \,\Q}$-modules cohérents à droite via
les foncteurs quasi-inverses $- \otimes _{\O _{\X} }\omega _{\X}$ et $- \otimes _{\O _{\X} }\omega _{\X} ^{-1}$.

L'inégalité de Bernstein  assure que, si $\E \neq 0$, alors $\dim \E \geq d _X$.
Autrement dit, un $F$-$\D ^{\dag} _{\X,\Q}$-module cohérent $\E$ est holonome si et seulement si
$\dim (\E)\leq d$ si et seulement si
$\mathrm{codim} (\E) \geq d$.

D'après la terminologie de Berthelot,
un $F$-$\D ^{\dag} _{\X,\Q}$-module cohérent $\E$ est \text{holonome} si
$\E=0$ ou si $\dim \E = d_X$.
De plus, un $F$-complexe $\E \in F$-$D ^\mathrm{b} _{\mathrm{coh}} (\D ^{\dag} _{\X,\Q})$ est
\textit{holonome}
si ses faisceaux de cohomologie le sont. Nous noterons $F$-$D^{\mathrm{b}} _{\mathrm{hol}}(\D ^{\dag} _{\X,\Q})$
la sous-catégorie pleine de $F$-$D^{\mathrm{b}} _{\mathrm{coh}}(\D ^{\dag} _{\X,\Q})$
des $F$-$\D ^{\dag }_{\X,\Q}$-complexes holonomes.
\medskip

$\bullet$ On dispose enfin du {\og critère homologique\fg} de l'holonomie dû à Virrion (voir \cite[III.4.2]{virrion}) qui nous sera très utile :
un $F$-$\D ^{\dag} _{\X,\Q}$-module cohérent $\E$
est holonome si et
seulement si, pour tout $i \neq  d _{\X}$, on ait
$ \mathcal{E} xt ^i _{\D _{\X ,\Q} ^{\dag}} ( \E , \D _{\X ,\Q} ^{\dag} ) =0.$
Ainsi, si $\E$ est un $F$-$\D ^{\dag} _{\X,\Q}$-module holonome
on dispose alors de l'isomorphisme canonique
$\mathcal{H} ^0 \DD (\E) \riso  \DD (\E) $.
De plus, Virrion a vérifié que
$\mathcal{H} ^0 \DD (\E) $ est aussi un $F$-$\D ^{\dag} _{\X,\Q}$-module holonome qui sera noté
$\E ^*$ (voir \cite[III.4.3]{virrion}).
Cela implique que l'holonomie est stable par le foncteur dual $\DD$ (voir \cite[III.4.4]{virrion}).

\end{vide}

\begin{vide}
  [Analogue arithmétique de Berthelot du théorème de Kashiwara]
\label{kashiwara}
On suppose ici que $f _0$ est une immersion fermée qui ne se relève pas forcément en un morphisme
de $\V$-schémas formels lisses.
  \begin{enumerate}
    \item  Pour tout $\D ^\dag _{\X} (\hdag T) _{\Q}$-module cohérent
    (resp. $F$-$\D ^\dag _{\X,\Q}$-module holonome)
    $\E$ à support dans $X '$,
    pour tout $\D ^\dag _{\X ' } (\hdag T ') _{ \Q}$-module cohérent $\E'$
    (resp. $F$-$\D ^\dag _{\X ',\Q}$-module holonome)
    $\mathcal{H} ^k {f _0} _+ (\E ) =0$ et $\mathcal{H} ^k  f _0 ^!(\E) =0$.
    \item Les foncteurs ${f _0} _+$ et $ f _0 ^!$ sont des équivalences quasi-inverses entre la catégorie des
    $\D ^\dag _{\X } (\hdag T ) _{ \Q}$-modules cohérents à support dans $X '$
    (resp. $F$-$\D ^\dag _{\X,\Q}$-modules holonomes à support dans $X'$) et celle des
    $\D ^\dag _{\X ' } (\hdag T') _{ \Q}$-modules cohérents
    (resp. $F$-$\D ^\dag _{\X',\Q}$-modules holonomes).
  \end{enumerate}
Le cas non respectif a été vérifié par Berthelot (voir \ref{kashiwara} ou \cite[3.1.6]{caro_surcoherent}).
Le cas respectif s'en déduit
grâce au critère d'holonomie de Virrion (voir \ref{cri-hol}) et à l'isomorphisme de dualité relative.
\end{vide}

\begin{vide}
  [Foncteur cohomologique local]
  \label{fct-cohloc}
  Nous exposons ici les principaux résultats de \cite[2.2]{caro_surcoherent}.
  Soit $Z$ un sous-schéma fermé de $X$.
  Le {\og foncteur cohomologique local à support strict dans $Z$\fg} construit dans \cite[2.2.6]{caro_surcoherent}
sera noté $\R \underline{\Gamma} ^\dag _Z \, :\,
\smash{\underset{^{\longrightarrow}}{LD}} ^{\mathrm{b}} _{\Q ,\mathrm{qc}}
  ( \smash{\widehat{\D}} _{\X } ^{(\bullet)}(T))
  \rightarrow \smash{\underset{^{\longrightarrow}}{LD}} ^{\mathrm{b}} _{\Q ,\mathrm{qc}}
  ( \smash{\widehat{\D}} _{\X } ^{(\bullet)}(T))$ (on dispose d'une deuxième construction
  par Berthelot donnée dans \cite[4.4.4, 4.4.5]{Beintro2}).
Rappelons en deux mots sa construction.
Lorsque $Z$ est un diviseur de $X$, $\R \underline{\Gamma} ^\dag _Z $ est par définition
le foncteur cohomologique local à support strict dans $Z$ construit par Berthelot (voir \cite[4.4.4, 4.4.5]{Beintro2}).
Si $Z$ est l'intersection des diviseurs $T _1, \dots, T _r$,
alors, pour tout $\E ^{(\bullet)} \in
\smash{\underset{^{\longrightarrow}}{LD}} ^{\mathrm{b}} _{\Q ,\mathrm{qc}}
  ( \smash{\widehat{\D}} _{\X } ^{(\bullet)}(T))$,
  on pose $\R \underline{\Gamma} ^\dag _Z (\E ^{(\bullet)} )
  = \R \underline{\Gamma} ^\dag _{T _1} \circ \cdots \circ
\R \underline{\Gamma} ^\dag _{T _r} (\E ^{(\bullet)} ) $.
Cette définition de $\R \underline{\Gamma} ^\dag _Z (\E ^{(\bullet)} )$
a bien un sens car celle-ci est (à isomorphisme canonique près) indépendante
du choix des diviseurs $T _1 ,\dots, T _r$ tels que $Z = \cap _{l=1,\dots, r} T _r$
(voir \cite[2.2.4--6]{caro_surcoherent}).
Pour tous sous-schémas fermés $Z$, $Z'$ de $X$,
on obtient alors
l'isomorphisme canonique compatible à Frobenius
et fonctoriel en $Z$ et $Z'$ (voir \cite[2.2.8]{caro_surcoherent}) :
\begin{equation}
  \label{GammaZGammaZ'}
  \R \underline{\Gamma} ^\dag _{Z\cap Z'}
\riso
\R \underline{\Gamma} ^\dag _Z \circ \R \underline{\Gamma} ^\dag _{Z'}.
\end{equation}
De plus, pour tous sous-schémas fermés $Z$, $Z'$, $Z''$ de $X$,
les deux morphismes composés canoniques suivants sont égaux :
\begin{gather}
  \notag
  \R \underline{\Gamma} ^\dag _{Z\cap Z'\cap Z''}
  \riso
\R \underline{\Gamma} ^\dag _{Z\cap Z'} \circ \R \underline{\Gamma} ^\dag _{Z''}
\riso
\R \underline{\Gamma} ^\dag _Z \circ \R \underline{\Gamma} ^\dag _{Z'}\circ \R \underline{\Gamma} ^\dag _{Z''},
\\
\label{transZZ'Z''}
\R \underline{\Gamma} ^\dag _{Z\cap Z'\cap Z''}
\riso
\R \underline{\Gamma} ^\dag _{Z} \circ \R \underline{\Gamma} ^\dag _{Z'\cap Z''}
\riso
\R \underline{\Gamma} ^\dag _Z \circ \R \underline{\Gamma} ^\dag _{Z'}\circ \R \underline{\Gamma} ^\dag _{Z''}.
\end{gather}
\medskip

$\bullet$ Lorsque $Z$ est un diviseur, pour tout $\E ^{(\bullet)} \in
\smash{\underset{^{\longrightarrow}}{LD}} ^{\mathrm{b}} _{\Q ,\mathrm{qc}}
  ( \smash{\widehat{\D}} _{\X } ^{(\bullet)}(T))$,
Berthelot a établi que l'on dispose du triangle
distingué de localisation en $Z$ (voir \cite[5.3.6]{Beintro2} ou \cite[2.2]{caro_surcoherent}):
\begin{equation}
  \label{tri-local}
\R \underline{\Gamma} ^\dag _Z (\E ^{(\bullet)})
\rightarrow \E ^{(\bullet)}
\rightarrow (\hdag Z) (\E ^{(\bullet)}) \rightarrow
\R \underline{\Gamma} ^\dag _Z (\E ^{(\bullet)}) [1],
\end{equation}
où le foncteur $(\hdag Z)$ est celui défini dans \ref{qc-vardiv}.
\medskip

$\bullet$
On étend alors naturellement la définition du foncteur
$(\hdag Z)$ au cas où $Z$ est un sous-schéma fermé quelconque
de $X$
en définissant $(\hdag Z) (\E ^{(\bullet)})$ comme égale
au cone du morphisme canonique
$\R \underline{\Gamma} ^\dag _Z (\E ^{(\bullet)})
\rightarrow \E ^{(\bullet)}$. On l'appelle foncteur de localisation en dehors de $Z$
(ou, de façon abusive, foncteur restriction en dehors de $Z$).
Si $Z, Z'$ sont deux sous-schémas fermés de $X$, on bénéficie d'après \cite[2.2.14]{caro_surcoherent}
de l'isomorphisme canonique compatible à Frobenius :
\begin{equation}
  \label{hdagcirchdag}
  (\hdag Z ) \circ (\hdag Z') (\E ^{(\bullet)}) \riso (\hdag Z \cup Z') (\E ^{(\bullet)}).
\end{equation}
\medskip

$\bullet$
D'après \cite[2.2.6.1]{caro_surcoherent}, pour tout $\E ^{(\bullet)}
\in \smash[b]{\underset{^{\longrightarrow }}{LD }}  ^\mathrm{b} _{\Q, \mathrm{qc}}
( \smash{\widehat{\D}} _{\X} ^{(\bullet)} )$,
\begin{equation}\label{gamma=otimegammaO}
\R \underline{\Gamma} ^\dag _{Z } ( \O _{\X, \Q}) \smash{\overset{\L}{\otimes}}
^{\dag}_{\O _{\X , \Q}} \E ^{(\bullet) } \riso \R \underline{\Gamma} ^\dag _{Z }(\E ^{(\bullet)})
,\
 (\hdag Z)( \O _{\X, \Q}) \smash{\overset{\L}{\otimes}}   ^{\dag}_{\O _{\X , \Q}} \E ^{(\bullet)}
 \riso (\hdag Z) (\E ^{(\bullet)}),
\end{equation}
où nous avons identifié le système inductif de $\smash[b]{\underset{^{\longrightarrow }}{LD }}  ^\mathrm{b} _{\Q, \mathrm{qc}}
( \smash{\widehat{\D}} _{\X} ^{(\bullet)} )$ constant égal à $\O _{\X} $ avec $\O _{\X,\Q}$.
Cela implique que les foncteurs de la forme $\R \underline{\Gamma} ^\dag _{Z }$
et $(\hdag Z)$ commutent deux à deux (voir \cite[2.2.6.1]{caro_surcoherent}).
\medskip

$\bullet$
Sans supposer que $f _0$ se relève, soit $Z$ un sous-schéma fermé de $X$ et
$Z': = f^{-1} (Z)$.
Pour tous
$\E ^{(\bullet)} \in
\smash{\underset{^{\longrightarrow}}{LD}} ^{\mathrm{b}} _{\Q ,\mathrm{qc}}
( \smash{\widehat{\D}} _{\X } ^{(\bullet)})$,
$\E ^{\prime \bullet} \in
\smash{\underset{^{\longrightarrow}}{LD}} ^{\mathrm{b}} _{\Q ,\mathrm{qc}}
  ( \smash{\widehat{\D}} _{\X '} ^{(\bullet)})$,
on dispose des isomorphismes fonctoriels en $Z$, compatibles à Frobenius (voir \cite[2.2.18]{caro_surcoherent}) :
\begin{gather}
\label{commutfonctcohlocal2}
  f _0^! \circ\R \underline{\Gamma} ^\dag _{Z }(\E ) \riso \R \underline{\Gamma} ^\dag _{Z  '}\circ f _0^! (\E ), \
f _0^! \circ (\hdag Z)(\E ) \riso (\hdag Z ') \circ f _0^! (\E ),
\\
\label{commutfonctcohlocal2+}
\R \underline{\Gamma} ^\dag _{Z}\circ f _{0+} (\E')
\riso
f _{0+} \circ \R \underline{\Gamma} ^\dag _{Z'}(\E' )
 , \
(\hdag Z )\circ f _{0+} (\E') \riso  f _{0+} \circ (\hdag Z')(\E' ).
\end{gather}
\medskip

$\bullet$
Soient $Z_1$, $Z_2$ deux sous-schémas fermés de $X$ et
$\E \in \smash[b]{\underset{^{\longrightarrow }}{LD }} ^\mathrm{b}
_{\Q, \mathrm{qc}} (\overset{^g}{} \widehat{\D} _{\X}
^{(\bullet)} )$.
On dispose aussi des triangles distingués de localisation de
Mayer-Vietoris (voir \cite[2.2.16]{caro_surcoherent}) :
\begin{gather}\label{eq1mayer-vietoris}
  \R \underline{\Gamma} ^\dag _{Z _1\cap Z_2}(\E ) \rightarrow
  \R \underline{\Gamma} ^\dag _{Z _1}(\E ) \oplus
\R \underline{\Gamma} ^\dag _{Z _2 }(\E )  \rightarrow
\R \underline{\Gamma} ^\dag _{Z _1\cup Z _2}(\E ) \rightarrow
\R \underline{\Gamma} ^\dag _{Z _1 \cap Z _2}(\E )[1],\\
 (\hdag Z _1 \cap Z _2)(\E ) \rightarrow
(\hdag Z _1)(\E ) \oplus  (\hdag Z _2) (\E )
\rightarrow   (\hdag Z _1 \cup Z _2)(\E ) \rightarrow (\hdag Z _1\cap Z _2)(\E )[1].
\end{gather}
\medskip

$\bullet$
Supposons maintenant que $f_0 \,:\, X' \hookrightarrow X$ soit une immersion fermée.
On bénéficie alors de l'isomorphisme canonique :
\begin{equation}
\label{445Beintro}
  \R \underline{\Gamma} ^\dag _{X'} (\E ^{(\bullet)}) \riso f _{0+} f^!_0 (\E ^{(\bullet)})
\end{equation}
En effet,
$f _{0+} f^!_0 (\E ^{(\bullet)})= f _{0+} ( \O _{\X',\Q}
\smash{\overset{\L}{\otimes}} ^{\dag}_{\O _{\X ', \Q}}
f^!_0 (\E ^{(\bullet)})
\riso
f _{0+} ( \O _{\X',\Q} [d_{X'/X}])
\smash{\overset{\L}{\otimes}} ^{\dag}_{\O _{\X , \Q}}
\E ^{(\bullet)}
\riso
f _{0+} f _0 ^! ( \O _{\X,\Q})
\smash{\overset{\L}{\otimes}} ^{\dag}_{\O _{\X , \Q}}
\E ^{(\bullet)},
$
où le premier isomorphisme est l'exemple de \cite[2.1.4]{caro_surcoherent}.
Avec \ref{gamma=otimegammaO}, il suffit donc de prouver
$\R \underline{\Gamma} ^\dag _{X'} ( \O _{\X,\Q})
\riso f _{0+} f _0 ^! ( \O _{\X,\Q}).$
Comme $\R \underline{\Gamma} ^\dag _{X'} ( \O _{\X,\Q})\in
D ^\mathrm{b} _\mathrm{coh} (\D ^\dag _{\X,\,\Q})$ et est à support dans $X'$ (voir \cite{Becohdiff}),
cela découle du théorème de Kashiwara et de \ref{commutfonctcohlocal2} (en remarquant aussi
que $\R \underline{\Gamma} ^\dag _{X'}$ est l'identité sur
$\smash{\underset{^{\longrightarrow}}{LD}} ^{\mathrm{b}} _{\Q ,\mathrm{qc}}
( \smash{\widehat{\D}} _{\X '} ^{(\bullet)})$).
Comme Berthelot a obtenu un isomorphisme similaire (voir \cite[4.4.4, 4.4.5]{Beintro2}),
cela implique que les deux constructions du foncteur cohomologique local
à support strict dans un sous-schéma fermé $Z$
se rejoignent lorsque $Z$ est lisse.

\end{vide}

\begin{vide}
[Surcohérence]
\label{defi-surcoh}
La notion de {\og surcohérence\fg} (voir \cite[3.1.1]{caro_surcoherent}) est définie de la manière suivante :
  soit $\E$ un ($F$-)$\D ^\dag _{\X } (\hdag T) _{ \Q}$-module cohérent (resp. un
objet de ($F$-)$D ^\mathrm{b} _{\mathrm{coh}} (\D ^\dag _{\X } (\hdag T) _{\Q})$). On dit que $\E$ est
{\og $\D ^\dag _{\X } (\hdag T) _{\Q}$-surcohérent\fg} si pour tout morphisme lisse de $\V$-schémas
formels lisses $g$ : $\PP \rightarrow \X$, pour tout
diviseur $H $ de $P$, $(\hdag H ) ( g ^* \E)$ est un
($F$-)$\D ^\dag _{\PP } (\hdag g ^{-1} (T)) _{ \Q}$-module cohérent (resp. un
objet de ($F$-)$D ^\mathrm{b} _{\mathrm{coh}} (\D ^\dag _{\PP} (\hdag g ^{-1} (T)) _{\Q})$).
On note ($F$-)$D ^\mathrm{b} _{\mathrm{surcoh}} (\D ^\dag _{\X} (\hdag T) _{ \Q})$,
la sous-catégorie pleine de
($F$-)$D ^\mathrm{b} _{\mathrm{coh}} (\D ^\dag _{\X} (\hdag T) _{ \Q})$ des complexes surcohérents.

Un complexe est surcohérent si et seulement si ses espaces de cohomologie le sont
(cela vient du fait que les foncteurs de la forme $g ^*$ et $(\hdag H ) $ sont exacts si $g$ est lisse
et $H$ est un diviseur).
La surcohérence est stable par image directe par un morphisme propre,
par image inverse extraordinaire et par foncteur cohomologique local (voir \cite[3.1.7, 3.1.9]{caro_surcoherent}).
\end{vide}

\begin{vide}
  [Isocristaux surconvergents associés aux $\D$-modules arithmétiques :
Cas de la compactification lisse]
\label{Cas-comp-lisse}

Le point de départ est le théorème suivant de Berthelot
(voir \cite{Be4} ou pour une version publiée \cite[4.4.5]{Be1} et \cite[2.2.12]{caro_courbe-nouveau})
qui améliore la caractérisation \cite[4.4.5]{Be1} (et \cite[4.6]{Be2} pour la commutation à Frobenius)
des isocristaux surconvergents :

\begin{theo**}[Berthelot]
\label{be4eqiso}
On note $\X _K$ la fibre générique de $\X$ comme
$K$-espace analytique rigide, et $\mathrm {sp}$ : $\X _K \rightarrow \X$ le morphisme de spécialisation.

Les foncteurs $\sp _*$ et $\sp ^*$ induisent des équivalences quasi-inverses entre la catégorie des isocristaux sur $Y$,
surconvergents le long de $T$, et celle des $\D ^{\dag}  _{\X } ( \hdag T) _{\Q}  $-modules cohérents
$\O  _{\X } ( \hdag T) _{\Q}  $-cohérents. Ceux-ci commutent en outre à Frobenius.

De plus, un $\D ^{\dag}  _{\X } ( \hdag T) _{\Q}  $-module cohérent $\E$ est
$\O  _{\X } ( \hdag T) _{\Q}  $-cohérent si et seulement si la restriction
$\E {|\Y}$ est $\O _{\Y , \Q}$-cohérente.
\end{theo**}

On suppose que $\X$ un $\V$-schéma formel séparé et lisse.
Soient $Z$ un sous-schéma fermé lisse de $X$ tel que $T \cap Z$ soit un diviseur de $Z$ (cette dernière condition
n'est pas vraiment restrictive). On pose $U:=Z \setminus T$.

Nous étendons (voir \cite[2.5.10]{caro-construction}) le théorème \ref{be4eqiso} de la manière suivante :
il existe un foncteur canonique pleinement fidèle noté $\sp _{Z \hookrightarrow \X, T,+}$
de la catégorie des isocristaux sur $U$ surconvergents le long de $T \cap Z$
dans celle des $\D ^{\dag}  _{\X } ( \hdag T) _{\Q}  $-modules surcohérents à support dans $Z$
(pour la surcohérence, il faut aussi consulter \cite[6.1.4]{caro_devissge_surcoh}).
En particulier, lorsque le diviseur $T$ est vide, on obtient un foncteur
pleinement fidèle noté $\sp _{Z \hookrightarrow \X,+}$
de la catégorie des isocristaux convergents sur $Z$
dans celle des $\D ^{\dag}  _{\X ,\Q}  $-modules surcohérents à support dans $Z$.

Lorsque $Z \hookrightarrow \X$ se relève en une immersion fermée
$u\,:\, \ZZ \hookrightarrow \X$ de $\V$-schémas formels lisses,
il suffit pour tout isocristal $E$ sur $U$ surconvergent le long de $T \cap Z$
de poser $\sp _{Z \hookrightarrow \X, T,+} (E):= u _{T+} \circ \sp _*(E)$. La pleine fidélité
se déduit du théorème de Kashiwara.
Les isocristaux surconvergents et les $\D$-modules arithmétiques
ont l'aspect {\it cristallin} suivant :
les images inverses par deux relèvements (induisant la même réduction modulo $\pi$) sont canoniquement isomorphes
(voir respectivement \cite[2.2.17]{Berig} et \cite[2.1.5]{Be2}).
Cela implique que le module $\sp _{Z \hookrightarrow \X, T,+} (E)$ ne dépend pas à isomorphisme canonique près du choix du
relèvement $u$.
On construit alors par recollement $\sp _{Z \hookrightarrow \X, T,+} (E)$
(pour les détails techniques, voir \cite[2]{caro-construction}).

De plus, le foncteur $\sp _{Z \hookrightarrow \X, T,+} $
commute à Frobenius (voir \cite[4.1.9]{caro-construction}), aux foncteurs duaux et images inverses (voir \cite[4]{caro-construction}).
\end{vide}

\begin{vide}
  [$F$-isocristaux surcohérents]
\label{Fisocsurcoh}
Soient $\X$ un $\V$-schéma formel propre et lisse, $T$ un diviseur de $X$, $\Y := \X \setminus T$,
$Z$ un sous-schéma fermé. On suppose que $U:=Z \setminus T$ est lisse (mais $Z$ n'est plus forcément lisse).

$\bullet$ On dispose (voir \cite[3.2.1]{caro_devissge_surcoh} et aussi
\cite[2.3.2]{caro-2006-surcoh-surcv} pour la simplification de la définition dans le cas propre avec structure de Frobenius)
la catégorie $F$-$\mathrm{Isoc} ^{\dag \dag}( \X, T, Z/K)$ de la manière suivante :
  \begin{itemize}
    \item les objets sont les
  $F$-$\smash{\D} ^\dag _{\X} (\hdag T) _\Q$-modules surcohérents $\E$ à support dans $Z$ tels qu'il existe un isocristal convergent $G$ sur $U$
vérifiant
    $\E |{\Y} \riso \sp _{U\hookrightarrow \Y, +} (G)$ ;
   \item les flèches sont les morphismes $F$-$\smash{\D} ^\dag _{\PP} (\hdag T) _\Q$-linéaires.
  \end{itemize}
Cette catégorie $F$-$\mathrm{Isoc} ^{\dag \dag}( \X, T, Z/K)$ ne dépend canoniquement pas du choix
de $(\X, T,Z)$ tels que $Z \setminus T=U$. Elle sera notée
$F$-$\mathrm{Isoc} ^{\dag \dag}( U/K)$ et ses objets sont les {\og $F$-isocristaux surcohérents sur $U$\fg}.
\medskip

$\bullet$ Lorsque $U$ est un $k$-schéma séparé lisse quelconque, on construit plus généralement la catégorie
$F$-$\mathrm{Isoc} ^{\dag \dag}( U/K)$ par recollement (voir \cite[2.2.4]{caro-2006-surcoh-surcv})
\end{vide}

\begin{vide}
  [Isocristaux surconvergents associés aux $\D$-modules arithmétiques :
Cas général]
On dispose du théorème suivant qui étend d'une certaine manière \ref{be4eqiso} (voir \cite[2.3.2]{caro-2006-surcoh-surcv}) :

\begin{theo}
\label{eqcat-gen}
Soit $U $ un $k$-schéma séparé et lisse. On dispose d'une équivalence canonique de catégories
$$\sp _{U,+}\ :\ F\text{-}\mathrm{Isoc} ^{\dag}( U/K )\cong F\text{-}\mathrm{Isoc} ^{\dag \dag}( U/K ).$$
\end{theo}

Donnons une esquisse de la preuve. Tout d'abord, les deux catégories étant construites par recollement,
on se ramène à supposer $U$ affine et lisse.
On construit dans ce cas explicitement le foncteur $\sp _{U,+}$ via le choix d'un plongement $U \hookrightarrow \X$,
avec $\X$ égale au complété $p$-adique de l'espace projectif $\P ^n _\V$ pour un certain entier $n$.
Lorsque que l'adhérence de $U$ dans $X$ est lisse, cette construction est canoniquement isomorphe
au foncteur construit dans \ref{Cas-comp-lisse}.
En procédant ensuite par récurrence sur la dimension de $U$,
quitte à se restreindre à un ouvert dense de $U$,
on peut alors se ramener au cas où $U$ vérifie certaines hypothèses techniques géométriques (voir la notion {\og la désingularisation idéale\fg}
de \cite[7.1.1]{caro_devissge_surcoh}).
Ces conditions géométriques nous permettent par {\it descente finie et étale} de se supposer que la compactification de $U$ est lisse.
Ce dernier cas a déjà été traité dans \ref{Cas-comp-lisse}.

\end{vide}

\section{Un critère d'holonomie}

Soit $\PP$ un $\V$-schéma formel lisse.
\begin{defi}
\label{defi-fef}
  Un $\D ^\dag _{\PP,\,\Q}$-module cohérent $\E$ (resp. un complexe $\E$ de
  $D _\mathrm{coh} ^\mathrm{b} (\D ^\dag _{\PP,\,\Q})$) est dit à {\og fibres extraordinaires finies\fg} si,
  pour tout point fermé $x$ de $P$, pour tout relèvement $i _x$ de l'immersion fermée $\Spec k(x) \hookrightarrow P$ induite par $x$,
  les espaces de cohomologie de $i ^! _x (\E )$ sont des $K$-espaces vectoriels de dimension finie.

Puisque la {\it surcohérence} (\cite[3.1.1]{caro_surcoherent})
est préservée par image inverse extraordinaire (\cite[3.1.7]{caro_surcoherent}),
on remarque qu'un $\D ^\dag _{\PP,\,\Q}$-module surcohérent (resp. un complexe de
  $D _\mathrm{coh} ^\mathrm{b} (\D ^\dag _{\PP,\,\Q})$) est à fibres extraordinaires finies.
\end{defi}

\begin{lemm}\label{lemmsurcoh=>hol}
  Soit $\E$ un $F\text{-}\D ^\dag _{\PP, \,\Q}$-module cohérent à fibres extraordinaires finies.
  En notant $Z$, le sous-schéma fermé réduit de $P$
  définissant son support, il existe un ouvert $\U$ de $\PP$,
  tel que $Y:=Z \cap U$ soit affine, lisse et dense dans $Z$ et
  tel que la restriction de $\E$ sur $\U$
  soit isomorphe à l'image par
  $\sp _{Y\hookrightarrow \U,+}$ (voir \ref{Cas-comp-lisse})
  d'un $F$-isocristal convergent sur $Y$. En particulier,
  $\E |{\U}$ est holonome et non nul.
\end{lemm}
\begin{proof}
Soit $\U'$ un ouvert de $\PP$ tel que $Y':= Z \cap U'$ soit affine,
lisse et dense dans $Z$. D'après Elkik (voir \cite{elkik}), on peut
alors relever $Y'$ en un $\V$-schéma formel affine et lisse $\Y'$.
Comme $\U'$ est lisse et $\Y'$ est affine, il existe un morphisme de
$\V$-schémas formels lisses de la forme $v\,:\, \Y '\hookrightarrow
\U'$ relevant $Y '\hookrightarrow \U'$. Comme $\E |\U'$ est un
$F\text{-}\D ^\dag _{\U',\Q}$-module cohérent à support dans $Y'$ et
à fibres extraordinaires finies, il résulte du théorème de Kashiwara
(voir \ref{kashiwara}) que $v ^! (\E|\U')$ est un $F\text{-}\D ^\dag
_{\Y',\Q}$-module cohérent et à fibres extraordinaires finies.
D'après \cite[2.2.17]{caro_courbe-nouveau}, il existe alors un
ouvert $\U $ de $\U'$ tel que $\Y:= \Y' \cap U$ soit affine et dense
dans $\Y'$ et tel que $v ^! (\E|\U') | \Y$ soit même $\O
_{\Y,\Q}$-cohérent. Notons $u\,:\, \Y \hookrightarrow \U$
l'immersion fermée induite par $v$. Ainsi $u ^! (\E |\U)$ est un
$F\text{-}\D ^\dag _{\Y,\Q}$-module cohérent , $\O _{\Y,
\Q}$-cohérent. D'après \cite[5.3.5.(i)]{Beintro2}, cela implique que
$u ^! (\E |\U)$ est holonome. D'après la version cohérente du
théorème de Kashiwara (voir \ref{kashiwara}), comme $\E |\U$ est à
support dans $Y$, $u _+  u ^! (\E |\U)\riso \E |\U$ et donc, d'après
la version holonome de \ref{kashiwara}, $\E |\U$ est holonome.

Enfin, d'après \ref{Cas-comp-lisse}, en notant $G:= \sp ^* (u ^! (\E|\U))$ le $F$-isocristal convergent sur $Y$ associé
à $u ^! (\E|\U)$, on obtient :
$\sp _{Y\hookrightarrow \U,+} (G) = u _+ ( \sp _* (G)) \riso u _+  u
^! (\E |\U)\riso \E |\U$. On a ainsi prouvé $\sp _{Y\hookrightarrow \U,+}
(G)  \riso \E |\U$. D'où le résultat.
\end{proof}

Nous aurons aussi besoin du résultat suivant qui étend et se déduit de \ref{lemmsurcoh=>hol} :
\begin{lemm}
  \label{lemmsurcoh=>hol2}
  Soient $\E _1,\dots, \E _n$ des $F\text{-}\D ^\dag _{\PP, \,\Q}$-modules cohérents à fibres extraordinaires finies non tous nuls.
  Il existe alors un ouvert $\U$ de $\PP$ tel que
  $\E _1 |\U,\dots, \E _n |\U$  soient holonomes et non tous nuls.
\end{lemm}

\begin{proof}
  On procède par récurrence sur le nombre $N$ de faisceaux $\E _i$ non nul (avec $i$ parcourant $1,\dots, n$).
  Lorsque $N =1$, ce résulte de \ref{lemmsurcoh=>hol} (en effet, le module nul est holonome).
  Supposons à présent $N \geq 2$. Quitte à réindexer, supposons
  $\E _1\neq 0$.
  D'après \ref{lemmsurcoh=>hol},
  il existe un ouvert $\U'$ de $\PP$ tel que $\E _1 |\U'$ soit holonome et non nul.
  Le cas où, pour tout $i\geq2$, $\E _i |\U'=0$ est immédiat.
  Autrement, supposons que les $\E _2 |\U', \dots, \E _n |\U'$ soient non tous nuls.
  Par hypothèse de récurrence, il existe un ouvert $\U$ de $\U'$ tel que
  $\E _2 |\U, \dots, \E _n |\U$ soient holonomes et non tous nuls. Le fait que
  $\E _1 |\U$ soit aussi holonome nous permet de conclure.
\end{proof}

\begin{lemm}
  \label{lemmsurcoh=>holrem}
  Pour tout $F\text{-}\D ^\dag _{\PP, \,\Q}$-module cohérent $\E$, $\mathcal{H} ^0 \DD (\E)$ est holonome.
  De plus, $\mathcal{H} ^0 \DD \circ \mathcal{H} ^0 \DD(\E)$ est le plus grand sous-$F\text{-}\D ^\dag _{\X,\,\Q}$-module
 cohérent de $\E$ qui soit holonome.
\end{lemm}
\begin{proof}
Notons $d$ la dimension de $P$.
D'après \cite[III.2.4.(i)]{virrion}, $\mathrm{codim} (\mathcal{E}xt _{\D ^\dag _{\PP, \Q}} ^{d}(\E, \D ^\dag _{\PP,\Q}) )\geq d$.
Comme $\mathcal{H} ^0 \DD (\E) \riso \mathcal{E}xt _{\D ^\dag _{\PP, \Q}} ^{d}(\E, \D ^\dag _{\PP,\Q})\otimes _{\O _{\X} }\omega _{\X} ^{-1}$,
on obtient ainsi $\mathrm{codim} (\mathcal{H} ^0 \DD (\E) ) \geq d$.
Or, d'après \ref{cri-hol} (voir la définition de l'holonomie et l'inégalité de Bernstein),
$\mathcal{H} ^0 \DD (\E) $ est holonome
si et seulement si
$\mathrm{codim} (\mathcal{H} ^0 \DD (\E) ) \geq d$. D'où l'holonomie de
$\mathcal{H} ^0 \DD (\E) $.

La deuxième assertion étant locale en $\PP$, on peut supposer $\PP$ affine.
Résumons d'abord la construction de la suite spectrale
étudiée par Virrion dans \cite[III.3.4]{virrion}.
Puisque $\E$ est cohérent, il existe
une résolution gauche
$\mathcal{L} ^\bullet \rightarrow \E$ où $\mathcal{L} ^\bullet $ est un complexe de $\D ^\dag _{\PP,\Q}$-modules localement projectifs de type fini
de longueur $d$ (concernant la longueur, voir \cite[III.3.3]{virrion}).
La deuxième suite spectrale d'hypercohomologie du foncteur
$\mathcal{H} om _{\D ^\dag _{\PP, \Q}} (-,\D ^\dag _{\PP, \Q})$ par rapport au complexe
$\mathcal{H} om _{\D ^\dag _{\PP, \Q}} ( \mathcal{L} ^\bullet,\D ^\dag _{\PP, \Q})$  donne la suite spectrale
(voir \cite[III.3.4, p.60]{virrion}) :
\begin{equation}
  \label{suitespectVirrion}
  E _2 ^{i,j} = \mathcal{E}xt _{\D ^\dag _{\PP, \Q}} ^{i}(\mathcal{E}xt _{\D ^\dag _{\PP, \Q}} ^{-j}(\E, \D ^\dag _{\PP,\Q}), \D ^\dag _{\PP,\Q} )
\Rightarrow
E ^n =\mathcal{H} ^n (\E).
\end{equation}
Virrion a établi dans la preuve de \cite[III.3.6]{virrion}
que la filtration décroissante $(F^l (\E) ) _{l=1,\dots, d}$
définie par la suite spectrale
\ref{suitespectVirrion} pour $n=0$ possède la propriété suivante :
$F^l (\E)$ est le plus grand sous-$F$-$\D ^\dag _{\PP,\Q}$-module cohérent de $\E$
de codimension supérieure ou égale à $l$.
Comme être de codimension supérieure ou égale à $d$ équivaut à être holonome,
cela implique que
$F^d (\E)$ est le plus grand sous-$F$-$\D ^\dag _{\PP,\Q}$-module cohérent de $\E$ qui soit holonome.
Or, les termes non nuls de $E ^{i,j} _2 $ sont à l'intérieur du triangle
défini par $i=-j$, $j=0$ et $i=d$ (voir \cite[III.3.4, p.61]{virrion}).
Cela implique que $ E _{2} ^{d,-d} = E _{\infty} ^{d,-d} = F^d (\E)$.
L'isomorphisme
$ E _{2} ^{d,-d} \riso \mathcal{H} ^0 \DD \circ \mathcal{H} ^0 \DD(\E)$
nous permet de conclure.
\end{proof}

\begin{prop}\label{surcoh=>hol}
  Soit $\E$ un $F\text{-}\D ^\dag _{\PP, \, \Q}$-module cohérent. Si les espaces de cohomologie de $\DD (\E)$ sont
  à fibres extraordinaires finies (e.g., si $\DD (\E)$ est $\D ^\dag _{\PP, \, \Q}$-surcohérent),
  alors $\E$ est holonome.
\end{prop}
\begin{proof}
Notons $d $ la dimension de $P$.
Pour tout entier $d \geq i \geq 0$,
  notons $\FF _i:= \mathcal{H} ^{-i} \DD (\E)$. Pour $d \geq i \geq 1$, désignons par
  $Z _i$ le support de $\FF _i$ et $Z:=\cup _{i=1,\dots, d} Z _i$ leur réunion.
D'après le critère d'holonomie de Virrion (voir \ref{cri-hol} ou \cite[III.4.2]{virrion}), $\E$ est holonome si et seulement si
$\FF _1,\dots , \FF _d$ sont tous nuls.

Par l'absurde, supposons que $\FF _1,\dots , \FF _d$ sont non tous nuls.
Il découle de \ref{lemmsurcoh=>hol2}
  qu'il existe un ouvert $\U$ de $\PP$ tel que
  $\FF _1 |\U,\dots, \FF _d |\U$ soient holonomes et non tous nuls.
 Par \ref{lemmsurcoh=>holrem}, $\FF _{0}$ est lui aussi holonome.
 Ainsi, les espaces de cohomologie de $\DD (\E|{ \U})$ sont tous holonomes, i.e.,
le complexe $\DD (\E|{ \U})$ est holonome.
Donc, grâce au théorème de bidualité (voir \ref{bidual} ou \cite[II.3.5]{virrion}) et
  à la préservation de l'holonomie par le foncteur $\DD$ (voir \cite[III.4.4]{virrion}),
  $\E |{ \U}$ est holonome.
 Or, d'après le critère d'holonomie de Virrion (voir \ref{cri-hol} ou \cite[III.4.2]{virrion}),
 cela implique que $\FF _1 |\U,\dots, \FF _d |\U$ sont tous nuls.
 On a ainsi abouti à une contradiction.
\end{proof}

\section{Surholonomie et stabilité}
On désignera par $\X$ un $\V$-schéma formel lisse.

\begin{defi}\label{defi-surhol}
  Soit $\E$ un objet de ($F$-)$D (\D ^\dag _{\X,\,\Q})$. On définit par récurrence sur
  l'entier $r\geq 0$, la notion
  de {\it $r$-surholonomie}, de la façon suivante :
  \begin{enumerate}
    \item $\E$ est {\it $0$-surholonome} si et seulement si
    $\E \in (F$-$)D _\mathrm{surcoh} ^\mathrm{b} (\D ^\dag _{\X,\,\Q})$ ;
    \item Pour tout entier $r \geq 1$, $\E$ est {$r$-surholonome} si et seulement si $\E $ est $r-1$-surholonome
    et, pour tout morphisme lisse $f$ : $\X' \rightarrow \X$ de $\V$-schémas formels, pour tout diviseur $T '$ de $X'$, le complexe
    $\DD (\hdag T ') f ^! (\E)$ est $r-1$-surholonome.
  \end{enumerate}
  On dit que $\E$ est {\it $\infty$-surholonome} ou tout simplement {\it surholonome}
  si $\E$ est $r$-surholonome pour tout entier $r$.
  On note ($F$-)$D ^\mathrm{b} _\mathrm{surhol} (\D ^\dag _{\X,\,\Q})$ la sous-catégorie pleine
  de ($F$-)$D ^\mathrm{b} _\mathrm{coh} (\D ^\dag _{\X,\,\Q})$ des complexes surholonomes.
  Enfin un ($F$-)$\D ^\dag _{\X,\,\Q}$-module est $r$-surholonome (resp. surholonome) s'il l'est en tant
  qu'objet de ($F$-)$D ^\mathrm{b} (\D ^\dag _{\X,\,\Q})$.
\end{defi}

\begin{rema}\label{rema-hol}
\begin{enumerate}
  \item \label{rema-hol(i)} Pour tout entier naturel $r ' \geq r$, si un complexe $\E$ de ($F$-)$D (\D ^\dag _{\X,\,\Q})$ est $r'$-surholonome alors
$\E$ est $r$-surholonome. En particulier $\E \in (F$-$)D _\mathrm{surcoh} ^\mathrm{b} (\D ^\dag _{\X,\,\Q})$.

\item Réciproquement, il est raisonnable de conjecturer que la surcohérence est stable par foncteur dual.
  Si cette conjecture était vraie alors, pour tout entier $r$,
  la notion de $r$-surholonomie serait équivalente à celle de surcohérence (i.e. $0$-surholonomie) et
  donc à celle de surholonomie. Par \ref{surcoh=>hol}, cela impliquerait donc aussi que
  les $F$-complexes surcohérents sont holonomes.

\end{enumerate}
\end{rema}

\begin{prop} \label{stab-surhol-Dual-hdag}
\begin{enumerate}
\item \label{stab-surhol-Dual-hdag3}
  Soient $r \in \N \cup \{\infty \}$ et
  $\E \in D ^\mathrm{b} _\mathrm{coh} (\D ^\dag _{\X,\,\Q})$ un complexe $r+1$-surholonome.
  Alors $\DD (\E)$ est $r$-surholonome

  \item \label{stab-surhol-Dual-hdag2}Un facteur direct d'un complexe $r$-surholonome est
  $r$-surholonome.

  \item \label{stab-surhol-Dual-hdag0} Soient $r \in \N \cup \{\infty \}$ et
  $\E ' \rightarrow \E \rightarrow \E'' \rightarrow \E' [1]$ un triangle distingué
  de $D  (\D ^\dag _{\X,\,\Q})$. Si deux des complexes sont $r$-surholonomes alors
  le troisième l'est. En particulier,
  $D ^\mathrm{b} _\mathrm{surhol} (\D ^\dag _{\X,\,\Q})$ est une sous-catégorie
  triangulée de $D ^\mathrm{b} _\mathrm{surcoh} (\D ^\dag _{\X,\,\Q})$.

  \item \label{stab-surhol-Dual-hdag-1}
  Soient $r \geq 1$ un entier ou $r =\infty$ et
  $\E$ un $F\text{-}\D ^\dag _{\X,\,\Q}$-module $r$-surholonome.
  Alors $\E$ est holonome.

  \item \label{stab-surhol-Dual-hdag1}
  Soient $r \in \N \cup \{\infty \}$ et
  $\E \in D ^\mathrm{b} _\mathrm{coh} (\D ^\dag _{\X,\,\Q})$ un complexe $r$-surholonome.
Alors, pour tout sous-schéma fermé $Z$ de $X$,
  $\R \underline{\Gamma} ^\dag _Z (\E)$ et $(\hdag Z)(\E)$ sont $r$-surholonomes.
\end{enumerate}
\end{prop}

\begin{proof}La propriété \ref{stab-surhol-Dual-hdag3}) est tautologique
(avec les notations de \ref{defi-surhol}, il suffit de prendre $T'$ vide et $f$ égal à l'identité).
Les propriétés \ref{stab-surhol-Dual-hdag2}) et
\ref{stab-surhol-Dual-hdag0})
sont vérifiées pour les complexes cohérents.
On vérifie alors \ref{stab-surhol-Dual-hdag2}) et
\ref{stab-surhol-Dual-hdag0}) lorsque $r=0$.
Le cas général s'en déduit par récurrence sur $r \in \N$.
Traitons à présent \ref{stab-surhol-Dual-hdag-1}).
D'après la remarque \ref{rema-hol}.\ref{rema-hol(i)}, il suffit de l'établir lorsque $r=1$.
Grâce à \ref{stab-surhol-Dual-hdag3}),
si $\E$ est un $F\text{-}\D ^\dag _{\X,\,\Q}$-module $1$-surholonome alors
$\DD (\E)$ est surcohérent. Par \ref{surcoh=>hol}, cela implique que $\E$ est holonome.
D'où \ref{stab-surhol-Dual-hdag-1}).
Maintenant, établissons \ref{stab-surhol-Dual-hdag1}).
Le cas $r =0$ a déjà été traité (voir \cite[3.1.5]{caro_surcoherent}).
Soient $r \geq 1$ un entier et $\E \in D ^\mathrm{b} _\mathrm{coh} (\D ^\dag _{\X,\,\Q})$ un complexe $r$-surholonome.

Traitons à présent le cas où $Z$ est un diviseur.
On obtient par récurrence sur $r$ que $\R \underline{\Gamma} ^\dag _Z (\E)$ et $(\hdag Z)(\E)$ sont $r-1$-surholonomes.
Puis, en utilisant \ref{commutfonctcohlocal2} puis \ref{hdagcirchdag},
pour tout morphisme lisse $f$ : $\X' \rightarrow \X$ de $\V$-schémas formels,
pour tout diviseur $T '$ de $X'$,
on dispose des isomorphismes :
$$\DD (\hdag T ') f ^! (\hdag Z) (\E) \riso
\DD (\hdag T ') (\hdag f ^{-1} (Z)) f ^!  (\E) \riso
\DD (\hdag T ' \cup f ^{-1} (Z)) f ^!  (\E) .$$
Le complexe $\DD (\hdag T ') f ^! (\hdag Z) (\E) $  est donc $r-1$-surholonome.
On a ainsi prouvé que $(\hdag Z) (\E)$ est $r$-surholonome (voir \ref{defi-surhol}).
En utilisant le triangle de localisation en $Z$ (voir \ref{tri-local}) et via \ref{stab-surhol-Dual-hdag0}), il en dérive que
$\R \underline{\Gamma} ^\dag _{Z } (\E)$ est $r$-surholonome.

Passons maintenant au cas où $Z$ est un sous-schéma fermé quelconque de $X$.
D'après \cite[2.2.5]{caro_surcoherent},
il existe des diviseurs $T _1, \dots, \,T _n$ de $X$ tels que
$Z =\cap _{i=1,\dots ,n} T _i$.
En notant $Z':= \cap _{i=2,\dots ,n} T _i$, on dispose du triangle distingué de Mayer-Vietoris
(voir \ref{eq1mayer-vietoris}) :
$$\R \underline{\Gamma} ^\dag _{Z } (\E) \rightarrow \R \underline{\Gamma} ^\dag _{Z '} (\E) \oplus \R \underline{\Gamma} ^\dag _{T_1 } (\E)
\rightarrow \R \underline{\Gamma} ^\dag _{Z'\cup T_1 } (\E)
\rightarrow \R \underline{\Gamma} ^\dag _{Z } (\E)  [1].$$
Comme $Z'$ et $Z '\cup T_1$ sont des intersections de $n-1$ diviseurs,
on termine la preuve grâce à \ref{stab-surhol-Dual-hdag0}) et en procédant par récurrence sur $n$.
\end{proof}

\begin{coro}\label{stab-surhol-dualT}
  Soient $T$ un diviseur de $X$ et $\E \in D ^\mathrm{b} _\mathrm{coh} (\D ^\dag _{\X} (\hdag T) _{\Q})$.
  Alors $\E \in D ^\mathrm{b} _\mathrm{surhol} (\D ^\dag _{\X,\,\Q})$ si et seulement si
  $\DD _T (\E) \in D ^\mathrm{b} _\mathrm{surhol} (\D ^\dag _{\X,\,\Q})$.
\end{coro}
\begin{proof}
Supposons $\E \in D ^\mathrm{b} _\mathrm{surhol} (\D ^\dag _{\X,\,\Q})$.
  On remarque que morphisme canonique $\E \rightarrow (\hdag T) (\E)$ est
  un morphisme de $\D ^\dag _{\X} (\hdag T) _{\Q}$-modules cohérents
  qui est un isomorphisme en dehors de $T$.
Par \cite[4.3.12]{Be1} (voir \ref{4312Be1}), on en déduit que
  $\E \rightarrow (\hdag T) (\E)$ est un isomorphisme.
  Or, d'après \cite[I.4.4]{virrion}, les foncteurs duaux commutent à l'extension des scalaires, i.e.,
$\DD _T \circ (\hdag T) (\E) \riso (\hdag T)\circ \DD (\E)$.
On obtient donc $\DD _T (\E) \riso (\hdag T)\circ \DD (\E)$.
  Il résulte alors de \ref{stab-surhol-Dual-hdag}.i) et \ref{stab-surhol-Dual-hdag}.v) que
  $\DD _T (\E) \in D ^\mathrm{b} _\mathrm{surhol} (\D ^\dag _{\X,\,\Q})$.
  La réciproque découle de l'isomorphisme de bidualité $\DD _T \circ \DD _T (\E) \riso \E$.
\end{proof}

A l'instar de la cohérence, on vérifie que la notion de
surholonomie est locale en $\X$. De plus,
on étend les propriétés standards
des modules cohérents aux modules surholonomes. Par exemple, on a les deux propositions suivantes :

\begin{prop}
\label{propkerimsurcoh}
  Soient $\Phi$ : $\E \rightarrow \FF$ un morphisme de $F\text{-}\D ^\dag _{\X,\,  \Q}$-modules cohérents
  et $r \in \N \cup \{ \infty \}$.
  Si $\E$ et $\FF$ sont $r$-surholonomes, alors
  $\ker \Phi $, $\mathrm{coker} \Phi$ et $\mathrm{Im} \Phi$ sont $r$-surholonomes.
\end{prop}
\begin{proof}
Soient $f$: $\X' \rightarrow \X $ un morphisme lisse de $\V$-schémas formels,
 $T'$ un diviseur de $X'$. Notons $f ^* := \mathcal{H} ^0 f ^! [-d_{X'/X}]$.
Procédons alors à une récurrence sur $r \geq 0$.
Comme $f$ est lisse, le foncteur $f ^*$ est exact et
l'image par $f ^*$ d'un $F\text{-}\D ^\dag _{\X,\,  \Q}$-module cohérent
est un $F\text{-}\D ^\dag _{\X',\,  \Q}$-module cohérent.
Comme l'extension $\D ^\dag _{\X',\,  \Q} \rightarrow \D ^\dag _{\X'} (\hdag T') _{\Q}$ est plate,
le foncteur $(\hdag T')$ de la catégorie des
$F\text{-}\D ^\dag _{\X',\,  \Q}$-modules cohérents
dans celle des $F\text{-}\D ^\dag _{\X'} (\hdag T') _{\Q}$-modules cohérents
est exact.
Il en résulte que
le foncteur $(\hdag T') \circ  f ^*$
qui va de
la catégorie des
$F\text{-}\D ^\dag _{\X,\,  \Q}$-modules cohérents
dans celle des
$F\text{-}\D ^\dag _{\X'} (\hdag T') _{\Q}$-modules cohérents
est exact et commute aux noyaux, conoyaux et images.
Comme les noyaux, conoyaux, images d'un morphisme de $F\text{-}\D ^\dag _{\X',\,  \Q}$-modules cohérents
sont $F\text{-}\D ^\dag _{\X',\,  \Q}$-cohérents, on vérifie alors le cas $r=0$.

Supposons maintenant le théorème validé pour $r-1 \geq 0$ et prouvons-le pour $r$.
Comme le foncteur $(\hdag T') \circ  f ^*$ est exact,
on se ramène par hypothèse de récurrence à prouver que $\DD (\ker \Phi )$, $\DD(\mathrm{coker} \Phi)$ et $\DD (\mathrm{Im} \Phi)$ sont
$r -1$-surholonomes.
Par \ref{stab-surhol-Dual-hdag}.\ref{stab-surhol-Dual-hdag-1}, comme $r \geq 1$,
$\Phi$ est un morphisme de $F\text{-}\D ^\dag _{\X,\,  \Q}$-modules holonomes.
Or, d'après \cite[5.3.5.(ii)]{Beintro2},
une suite exacte de $F\text{-}\D ^\dag _{\X,\,  \Q}$-modules cohérents dont le terme du milieu est holonome
est une suite exacte de $F\text{-}\D ^\dag _{\X,\,  \Q}$-modules holonomes.
Il en résulte que $\ker \Phi $, $\mathrm{coker} \Phi$  et $\mathrm{Im} \Phi$ sont holonomes.
Rappelons que d'après le critère homologique de Virrion de l'holonomie (voir \ref{cri-hol}),
si $\G$ est un $F\text{-}\D ^\dag _{\X,\,  \Q}$-module holonome,
alors en posant $\G ^* := \mathcal{H} ^0 \DD (\G)$ on obtient
$\G ^* \riso \DD (\G)$.
Notons $\Phi ^* := \mathcal{H} ^0 \DD (\Phi)\,:\,\FF ^* \rightarrow \E ^* $ le morphisme dual de $\Phi$.
Comme le foncteur $\mathcal{H} ^0 \DD$ est exact sur la catégorie des
$F\text{-}\D ^\dag _{\X,\,  \Q}$-modules holonomes,
alors $(\mathrm{Im} (\phi )) ^*\riso \mathrm{Im} (\phi ^*) $,  $(\ker \Phi ) ^* \riso \mathrm{coker}(\Phi ^*) $,
$(\mathrm{coker} \Phi ) ^* \riso \ker (\Phi ^*) $.
Or, $\Phi ^*$ est un morphisme de
$F\text{-}\D ^\dag _{\X,\,  \Q}$-modules $r-1$-surholonomes.
Par hypothèse de récurrence, il en découle que
$(\mathrm{Im} (\phi )) ^*$,  $(\ker \Phi ) ^* $,
$(\mathrm{coker} \Phi ) ^*  $
sont $r-1$-surholonomes, ce qu'il fallait démontrer.
\end{proof}

\begin{prop}\label{5surhol}
Soit $\E _1 \rightarrow \E _2 \rightarrow \E _3 \rightarrow \E _4 \rightarrow \E _5$
une suite exacte de $F\text{-}\D ^\dag _{\X ,\, \Q}$-modules cohérents.
Pour tout $r\in \N \cup \{ \infty \}$,
si $\E _1$, $\E _2$, $\E _4$, $\E _5$ sont $r$-surholonomes,
alors $\E _3$ est $r$-surholonome.
\end{prop}
\begin{proof}
La preuve est analogue à celle de \ref{propkerimsurcoh} :
soient $f$: $\X' \rightarrow \X $ un morphisme lisse de $\V$-schémas formels,
 $T'$ un diviseur de $X'$. Traitons le cas $r=0$. Comme le foncteur
$(\hdag T') \circ  f ^*$ est exact,
on obtient la suite exacte
$(\hdag T') \circ  f ^* (\E _1) \rightarrow (\hdag T') \circ  f ^* (\E _2) \rightarrow
(\hdag T') \circ  f ^* (\E _3 ) \rightarrow
(\hdag T') \circ  f ^*( \E _4 )\rightarrow (\hdag T') \circ  f ^* (\E _5)$.
Comme par hypothèse tous les termes sauf celui du milieu sont
des $F\text{-}\D ^\dag _{\X',\,  \Q}$-modules cohérents, celui du milieu l'est aussi. Ainsi, $\E _3$ est
un $F\text{-}\D ^\dag _{\X,\,  \Q}$-module surcohérent.

Supposons maintenant le théorème validé pour $r-1 \geq 0$ et prouvons-le pour $r$.
Comme le foncteur $(\hdag T') \circ  f ^*$ est exact,
on se ramène par hypothèse de récurrence à prouver que
$\DD ( \E _3)$ est $r-1$-surholonome.

Par \ref{stab-surhol-Dual-hdag}.\ref{stab-surhol-Dual-hdag-1}, comme $r \geq 1$,
$\E _1$, $\E _2$, $\E _4$, $\E _5$ sont holonomes.
Notons $\FF _2$ (resp. $\FF _3$) l'image de $\E _2 \rightarrow \E _3$
(resp. $\E _3 \rightarrow \E _4$).
D'après \cite[5.3.5.(ii)]{Beintro2}, comme $\E _2$ (resp. $\E _4$)
est holonome, $\FF _2$ (resp. $\FF _3$) est aussi holonome.
Encore d'après \cite[5.3.5.(ii)]{Beintro2}, comme
on dispose de la suite exacte $0\rightarrow \FF _2 \rightarrow \E _3 \rightarrow  \FF _3 \rightarrow 0$
avec $\FF _2$ et $\FF _3$ holonomes alors $\E _3$ est holonome.

Ainsi, $\DD ( \E _3)\riso \E _3 ^*$. Comme le foncteur $\mathcal{H} ^0 \DD$ est exact sur la catégorie des
$F\text{-}\D ^\dag _{\X,\,  \Q}$-modules holonomes,
on obtient la suite exacte
$\E _5 ^* \rightarrow \E _4 ^* \rightarrow \E _3 ^* \rightarrow \E _2 ^* \rightarrow \E _1^*$.
Par hypothèse de récurrence, $\E _3 ^* $ est donc $r-1$-surholonome.
\end{proof}

\begin{coro}\label{sspectsurhol}
  Pour $N\geq p, q \geq 0$, $r _0\geq 1$,
  soit $\E _{r _0} ^{p,q} \Rightarrow \E ^n$ une suite spectrale de
  $F\text{-}\D ^\dag _{\X,\,  \Q}$-modules cohérents.
Si, pour tous $p,q$, les
 $F\text{-}\D ^\dag _{\X,\,  \Q}$-modules $\E _{r _0} ^{p,q}$ sont $r$-surholonomes, alors
 les $F\text{-}\D ^\dag _{\X,\,  \Q}$-modules $\E ^n$ sont, pour tout $n$, $r$-surholonomes.
\end{coro}

\begin{proof}
Il découle de \ref{propkerimsurcoh} que les $\E _{\infty} ^{p,q} $ sont $r$-surholonomes.
Par \ref{5surhol}, on en déduit que les $\E ^n$ sont $r$-surholonomes.
\end{proof}

\begin{theo}\label{stabinminv}
Soient $\Y$ et $\X$ deux $\V$-schémas formels lisses et
$f _0$ : $Y \rightarrow X$ un morphisme entre leur fibre spéciale.
  Pour tout $\E \in D ^\mathrm{b} _{\mathrm{coh}} (\D ^\dag _{\X ,\, \Q})$,
  pour tout entier $r$,
  si $\E$ est $r$-surholonome alors
  $f _0 ^! (\E) $ est $r$-surholonome.
\end{theo}
    \begin{proof}
    Comme le théorème est local en $\Y$,
    on peut supposer que $f_0$ se relève en un morphisme $f$ : $\Y \rightarrow \X$.
    Procédons alors par récurrence sur l'entier $r \geq 0$.
    Pour $r=0$, c'est \cite[3.1.7]{caro_surcoherent}. Supposons donc le théorème validé
    pour $r-1 \geq 0$ et prouvons-le pour $r$.

Puisque $f$ se décompose en une immersion (le graphe de $f$ : $\Y \hookrightarrow \Y \times _\S \X $)
    suivie de la projection lisse $\Y \times _\S \X \rightarrow \X$,
il suffit de traiter le cas où $f$ est lisse et celui où $f$ est une immersion fermée.
Soient $g$ : ${\Y '}  \rightarrow \Y$ un morphisme lisse et
    $Z$ un diviseur de $Y' $.
Par hypothèse de récurrence sur $r$, il s'agit d'établir que $\DD\circ (\hdag Z) \circ g ^! (f ^! (\E) )$
est $r-1$-surholonome.
Lorsque $f$ est lisse, cela découle de l'isomorphisme
$\DD\circ (\hdag Z) \circ g ^! (f ^! (\E) )\riso \DD\circ (\hdag Z) \circ (f\circ g) ^! (\E) )$
et du fait que $\E$ est $r$-surholonome.
On se ramène donc à traiter le cas où $f$ est une immersion fermée.

     Il suffit de prouver que pour tout morphisme lisse $g$ : ${\Y '}  \rightarrow \Y$ et
    pour tout diviseur $Z$ de $Y' $ le complexe
    $\DD\circ \R \underline{\Gamma} ^\dag _{Z } \circ  g ^! (f ^! (\E) )$ est $r-1$-surholonome.
En effet, en prenant $Z$ égal à l'ensemble vide, cela implique que
$\DD\circ g ^! (f ^! (\E) )$ est $r-1$-surholonome.
En appliquant le foncteur dual $\DD$ au triangle de localisation de
$g ^! (f ^! (\E) )$ en $Z$ (voir \ref{tri-local}),
il découle alors de \ref{stab-surhol-Dual-hdag}.\ref{stab-surhol-Dual-hdag0}
que $\DD\circ (\hdag Z)\circ  g ^! (f ^! (\E) )$ est $r-1$-surholonome.

En fixant $g$ et $Z$, prouvons à présent que $\DD\circ \R \underline{\Gamma} ^\dag _{Z } \circ  g ^! (f ^! (\E) )$ est $r-1$-surholonome.
Comme cela est local en ${\Y '} $,
    on peut donc supposer que $g$ se décompose en une immersion fermée ${\Y '} \hookrightarrow \widehat{\A} _{\Y} ^n$
    suivie de la projection $\widehat{\A} _{\Y} ^n \rightarrow \Y$. En notant $p$ la projection
    $\widehat{\A} _{\X } ^n \rightarrow \X$ et $i $ l'immersion fermée
    ${\Y '}  \hookrightarrow \widehat{\A} _{\X } ^n$, on obtient $f\circ g = p\circ i$.
Ainsi, $\DD\circ  \R \underline{\Gamma} ^\dag _{Z } \circ g ^! (f ^! (\E) )\riso
\DD \circ \R \underline{\Gamma} ^\dag _{Z } \circ i ^! \circ  p ^! (\E)$.
    Comme l'image inverse extraordinaire commute au foncteur cohomologique local (voir \ref{commutfonctcohlocal2}),
$\R \underline{\Gamma} ^\dag _{Z } i ^!  \riso i ^! \R \underline{\Gamma} ^\dag _{Z }$.
    D'où :
    $\DD\circ  \R \underline{\Gamma} ^\dag _{Z }\circ  g ^! (f ^! (\E) )\riso
     \DD \circ i ^! \circ \R \underline{\Gamma} ^\dag _{Z }\circ  p ^! (\E)$.

En utilisant le théorème de Kashiwara (\ref{kashiwara})
      puis le théorème de dualité relative (\cite{Vir04} ou \ref{iso-dualrel}),
on obtient sur $D ^\mathrm{b} _{\mathrm{coh}} (\D ^\dag _{\Y ' ,\, \Q}) $ :
$\DD \riso i^! i_+ \DD \riso i^! \DD i_+ $.
Or, comme $\E$ est surcohérent,
$i ^! \R \underline{\Gamma} ^\dag _{Z } p ^! (\E)
\in D ^\mathrm{b} _{\mathrm{coh}} (\D ^\dag _{\Y ' ,\, \Q}) $.
D'où :
$\DD (i ^! \R \underline{\Gamma} ^\dag _{Z } p ^! (\E))
\riso
i^! \DD i_+ (i ^! \R \underline{\Gamma} ^\dag _{Z } p ^! (\E))$.
Comme $\R \underline{\Gamma} ^\dag _{Z } p ^! (\E) \in D ^\mathrm{b} _{\mathrm{coh}} (\D ^\dag _{\widehat{\A} _{\X } ^n ,\, \Q}) $
et est à support dans $Y'$ (car $Z \subset Y'$),
le théorème de Kashiwara implique
$i^! \DD i_+ (i ^! \R \underline{\Gamma} ^\dag _{Z } p ^! (\E))
\riso
i ^! \DD\circ \R \underline{\Gamma} ^\dag _{Z } p ^! (\E)$. On a ainsi établi
l'isomorphisme
$\DD\circ  \R \underline{\Gamma} ^\dag _{Z } g ^! (f ^! (\E) )
\riso
     i ^! \DD\circ \R \underline{\Gamma} ^\dag _{Z } p ^! (\E)$.

     Or, d'après \ref{stab-surhol-Dual-hdag}, comme $\E$ est $r$-surholonome le complexe
     $\DD\circ  \R \underline{\Gamma} ^\dag _{Z } p ^! (\E)$ est $r-1$-surholonome.
     Par hypothèse de récurrence, $  i ^! \DD\circ \R \underline{\Gamma} ^\dag _{Z } p ^! (\E)$
     est aussi $r-1$-surholonome.
     \end{proof}

\begin{theo}\label{stabinmdir}
Soient $\Y$ et $\X$ deux $\V$-schémas formels lisses et
$f _0$ : $Y \rightarrow X$ un morphisme propre entre leur fibre spéciale.
On suppose que l'une des conditions suivantes est validée :
\begin{enumerate}
  \item \label{stabinmdir1} Le morphisme $f _0$ se relève en un morphisme $f$ : $\Y \rightarrow \X$ de
  $\V$-schémas formels lisses ;
  \item \label{stabinmdir2} Le morphisme $f _0$ est une immersion fermée ;
  \item \label{stabinmdir3} Le $\V$-schéma formel $\Y$ se plonge dans un $\V$-schéma formel propre et lisse.
\end{enumerate}

  Pour tout $\FF \in D ^\mathrm{b} _{\mathrm{coh}} (\D ^\dag _{\Y ,\, \Q})$,
  pour tout entier $r$, si $\FF$ est $r$-surholonome alors
  $f _{0+} (\FF) $ est $r$-surholonome.
\end{theo}
\begin{proof}
Dans le cas \ref{stabinmdir3}), on suppose que $\Y$ se plonge dans un $\V$-schéma formel $\PP$ propre
et lisse. Il en résulte une décomposition de $f _0$ en une immersion ($Y \hookrightarrow P \times X$) suivie
d'un morphisme propre et lisse (la projection $P \times X \rightarrow X$) qui se relève en
$\PP \times \X \rightarrow \X$.
Comme $f _0$ est propre, l'immersion $Y \hookrightarrow P \times X$ est une immersion fermée
(voir \cite[II.4.8.(e)]{HaAG}).
Il suffit donc de traiter les deux cas \ref{stabinmdir2}) et \ref{stabinmdir1}).
Or, le cas \ref{stabinmdir2}) se ramène au cas relevable \ref{stabinmdir1}).
En effet, dans le cas où
$f _0$ est une immersion fermée,
comme le théorème est local en $\X$
on peut supposer $\X$ et donc $\Y$ affines. Il existe alors un relèvement
$f $ : $\Y \rightarrow \X$ de $f _0$.
Il suffit donc de traiter le cas où le morphisme $f _0$ se relève en un morphisme propre $f$ : $\Y \rightarrow \X$.

  Maintenant, traitons le cas \ref{stabinmdir1}) et procédons à une récurrence sur l'entier $r \geq 0$.
  Le cas $r =0$ a été démontré dans \cite[3.1.9]{caro_surcoherent}.
Supposons donc le théorème validé pour $r-1\geq 0$ et prouvons-le pour $r$.
Soient $\FF$ un complexe $r$-surholonome de $\D ^\dag _{\Y ,\, \Q}$-modules,
$g$ : $\X '\rightarrow \X$ un morphisme lisse de $\V$-schémas formels
et $T'$ un diviseur de $X'$. Notons $f '$ : $ \Y \times _{\X} \X ' \rightarrow \X '$
et $g '$ : $ \Y \times _{\X} \X ' \rightarrow \Y $ les projections respectives.
Il s'agit de prouver que
$\DD (\hdag T') g ^! f _+(\FF)$ est $r-1$ surholonome.
Or, en utilisant l'isomorphisme de changement de base de l'image directe d'un morphisme propre
par un morphisme lisse $g ^! f _+ (\FF) \riso f ' _+  g ^{\prime !} (\FF)$
(voir \cite[3.1.8]{caro_surcoherent}), la commutation de l'image directe au foncteur cohomologique
local (voir \ref{commutfonctcohlocal2}) et
au théorème de dualité relative de Virrion
(voir \cite{Vir04} ou \ref{iso-dualrel} pour la version avec diviseur),
on obtient les isomorphismes :
$$\DD (\hdag T') g ^! f _+(\FF) \riso
\DD (\hdag T') f ' _+  g ^{\prime !} (\FF) \riso
\DD f ' _+ (\hdag f ^{\prime -1} T')  g ^{\prime !} (\FF) \riso
f ' _+ \DD (\hdag f ^{\prime -1} T')  g ^{\prime !} (\FF).$$
D'après la proposition \ref{stab-surhol-Dual-hdag} et le théorème \ref{stabinminv}, le complexe
$\DD (\hdag f ^{\prime -1} T')  g ^{\prime !} (\FF)$ est $r-1$-surholonome.
Par hypothèse de récurrence, il en découle que le complexe
$f ' _+ \DD (\hdag f ^{\prime -1} T')  g ^{\prime !} (\FF)$ est $r-1$-surholonome.
\end{proof}

\begin{theo}[Kashiwara]\label{kashiwarasurhol}
Soient $\Y$ et $\X$ deux $\V$-schémas formels lisses,
$f _0$ : $Y \hookrightarrow X$ une immersion fermée entre leur fibre spéciale.
  \begin{enumerate}
    \item \label{kashiwarasurhol1} Pour tout $r \in \N \cup \{\infty\}$,
    pour tout $\D ^\dag _{\X,\Q}$-module $r$-surholonome
    $\E$ à support dans $Y$,
    tout $\D ^\dag _{\Y , \Q}$-module $r$-surholonome $\FF$ et
    tout entier $k \neq 0$, $\mathcal{H} ^k {f _0} _+ (\FF ) =0$ et $\mathcal{H} ^k  f _0 ^!(\E) =0$.
    \item \label{kashiwarasurhol2} Les foncteurs ${f _0} _+$ et $ f _0 ^!$ induisent des équivalences quasi-inverses entre la catégorie des
    $\D ^\dag _{\X , \Q}$-modules $r$-surholonomes à support dans $Y$ et celle des
    $\D ^\dag _{\Y , \Q}$-modules $r$-surholonomes.
  \end{enumerate}
\end{theo}
\begin{proof}
La partie \ref{kashiwarasurhol1}) résulte de la version cohérente du théorème de Kashiwara
dû à Berthelot (voir \ref{kashiwara}) et du fait qu'un $\D ^\dag _{\Y , \Q}$-module $r$-surholonome
(resp. un $\D ^\dag _{\X,\Q}$-module $r$-surholonome à support dans $Y$) est
un $\D ^\dag _{\Y , \Q}$-module cohérent (resp. un $\D ^\dag _{\X,\Q}$-module cohérent à support dans $Y$).
L'assertion \ref{kashiwarasurhol2}) découle à nouveau de \ref{kashiwara} mais aussi de la stabilité
de la $r$-surholonomie par image directe par une immersion fermée et image inverse
extraordinaire (voir \ref{stabinminv} et de \ref{stabinmdir}).
\end{proof}

\begin{prop}
  Pour tout $r \in \N \cup \{ \infty\}$, on a l'équivalence entre les deux assertions :
  \begin{enumerate}
    \item \label{stabotimes} La $r$-surholonomie est stable par produit tensoriel interne ;
    \item \label{stabboxtimes} La $r$-surholonomie est stable par produit tensoriel externe.
  \end{enumerate}

\end{prop}

\begin{proof}
  Comme la $r$-surholonomie est préservée
  par image inverse extraordinaire (voir \ref{stabinminv}), on obtient d'abord
  par définition du produit tensoriel externe (voir \ref{def-boxtimes})
  l'implication $\ref{stabotimes}) \Rightarrow \ref{stabboxtimes})$.

  Réciproquement,
notons $p _1$, $p _2$ : $\X \times _\S \X \rightarrow \X$ les projections respectives à gauche et à droite,
  $\delta $ : $ \X \hookrightarrow \X \times _\S \X$ l'immersion diagonale.
  Pour tous complexes
  $ \E _1$ et $\E _2 \in \smash{\underset{^{\longrightarrow}}{LD}} ^{\mathrm{b}} _{\Q ,\mathrm{qc}}
( \smash{\widehat{\D}} _{\X} ^{(\bullet)})$,
\begin{equation}
\label{boxtimes-otimes}
  \delta ^! ( \E _1 \smash{\overset{\L}{\boxtimes}} _{\O _\S}^{\dag} \E _2)
  =
  \delta ^! ( p _1 ^! (\E _1 ) \smash{\overset{\L}{\otimes}}^{\dag}_{\O _{\X \times _\S \X,\, \Q}} p _2^!(\E _2)
  [-2 d _{\X}])
  \riso
  \E _1 \smash{\overset{\L}{\otimes}}^{\dag}_{\O _{\X ,\, \Q}} \E _2 [- d _{\X}],
\end{equation}
l'isomorphisme résultant de la commutation de l'image inverse extraordinaire (décalée de la dimension
relative) aux produits tensoriels internes
(voir par exemple \cite[1.2.22]{caro_surcoherent}).
Si la $r$-surholonomie est préservée par produits tensoriels externes,
comme elle l'est aussi par image inverse extraordinaire (voir \ref{stabinminv}),
par \ref{boxtimes-otimes}, on obtient alors la stabilité par produits tensoriels internes.
On a ainsi vérifié $\ref{stabboxtimes}) \Rightarrow \ref{stabotimes})$.
\end{proof}

\begin{vide}
[Transitivité de l'image directe par un morphisme propre pour les complexes surholonomes]
\label{compatsurhol}
  Soient $f$ : $\X '\rightarrow \X$,
  $g$ : $\X'' \rightarrow \X'$ et $h$ : $\X'''\rightarrow \X ''$, trois morphismes propres
de $\V$-schémas formels lisses.

$\bullet$
Il découle de la remarque \cite[1.2.16]{caro_courbe-nouveau}
que l'on dispose, pour tout $\G \in  D _{\mathrm{surhol}} ^{\mathrm{b}} ( \D ^{\dag} _{\X '',\,\Q} )$,
de l'isomorphisme canonique
$(f\circ g ) _{+} (\G) \riso f _{+} \circ g _{+}(\G)$,
fonctoriel en $\G$.
On obtient ainsi des {\og isomorphismes de transitivité de l'image directe par un morphisme propre
pour les complexes surholonomes\fg}.
En outre, d'après \cite[1.2.16]{caro_courbe-nouveau}, celui-ci est par construction compatible à Frobenius
(voir \ref{commFrobf+} pour la commutation de l'image directe à Frobenius)
et associatif i.e.,
pour tout $\mathcal{H} \in   D _{\mathrm{surhol}} ^{\mathrm{b}} ( \D ^{\dag} _{\X ''',\Q} )$,
les deux morphismes composés
$(f\circ g \circ h) _{+} (\mathcal{H})  \riso (f\circ g) _{+} \circ h _{+} (\mathcal{H})
\riso $ $f _{+}\circ  g _{+}\circ  h _{+}  (\mathcal{H})  $ et
$(f\circ g \circ h) _{+} (\mathcal{H})  \riso f _{+}\circ  (g\circ h) _{+} (\mathcal{H})
\riso f _{+}\circ  g _{+}\circ   h _{+}  (\mathcal{H})  $ sont égaux.
\medskip

$\bullet$
On dispose d'après Berthelot d'une seconde construction plus naturelle de
l'isomorphisme de transitivité :
$(f\circ g ) _{+} (\G ) \riso f _{+} \circ g _{+}(\G )$
(pour une description de la construction, voir par exemple \cite[1.1.2]{caro-construction})
dont la compatibilité à Frobenius est loin d'être évidente.

\medskip

$\bullet$
En fait, d'après la remarque \cite[1.2.9]{caro-construction},
ces deux isomorphismes de transitivité du foncteur image directe sont identiques.
Pour le vérifier, nous avons établi que les morphismes d'adjonction entre image directe et image inverse extraordinaire
par un morphisme propre
sont transitifs, en prenant pour isomorphisme de transitivité de l'image directe ceux construits par Berthelot.

\end{vide}

\begin{vide}
\label{adjonction-surhol}
Soient $f$ : $\X '\rightarrow \X$ un morphisme propre de $\V$-schémas formels lisses,
$\E \in  D _{\mathrm{surhol}} ^{\mathrm{b}}( \D ^{\dag} _{\X ,\,\Q} )$,
$\E'  \in  D _{\mathrm{surhol}} ^{\mathrm{b}} ( \D ^{\dag} _{\X ',\,\Q} )$.
Comme un complexe surholonome est en particulier cohérent et que la surholonomie est
stable par image inverse extraordinaire
et par image directe par un morphisme propre,
par \ref{f+adjf!},
on dispose des morphismes canoniques d'adjonction
$\E'  \rightarrow f ^! f _+ (\E' )$,
$f _+ f ^! (\E) \rightarrow \E$.

En outre, si $f$ est une immersion fermée,
ces morphismes d'adjonction sont des isomorphismes pour
tout $\E \in  D _{\mathrm{surhol}} ^{\mathrm{b}}( \D ^{\dag} _{\X ,\,\Q} )$
à support dans $X '$
et tout
$\E'  \in  D _{\mathrm{surhol}} ^{\mathrm{b}} ( \D ^{\dag} _{\X ',\,\Q} )$.
Enfin, pour tout $\E \in  D _{\mathrm{surhol}} ^{\mathrm{b}}( \D ^{\dag} _{\X ,\,\Q} )$,
on bénéficie comme dans le cas surcohérent (\cite[3.1.12]{caro_surcoherent})
d'un isomorphisme compatible à Frobenius
$f _+ f ^! (\E ) \riso \R \underline{\Gamma} ^\dag _{X'} (\E)$ s'inscrivant dans le diagramme canonique :
\begin{equation}
  \label{+!=gamma}
  \xymatrix @R=0,3cm {
     {f _{+} \circ f  ^! (\E)} \ar[rr]^-{\mathrm{adj}} \ar[d] _-\sim && {\E} \ar@{=}[d] \\
  {\R \underline{\Gamma} ^\dag _{X'  } (\E)}
  \ar[rr]  && {\E}.}
\end{equation}

\end{vide}

\section{$\D$-modules arithmétiques surholonomes sur une variété. Cas simple}
Pour construire les complexes de $\D$-modules arithmétiques surholonomes sur une variété $U$,
nous étudions dans cette section le cas où la variété $U$ se plonge dans un $\V$-schéma formel propre et lisse
(voir \ref{def-surhol-relebis}).

\begin{nota}\label{notations-catzt}
  Soient $\X$ un $\V$-schéma formel lisse, $T$ et $Z$ deux sous-schémas fermés de $X$ et $U := Z \setminus T$.
On note alors ($F$-)$\mathfrak{M} ^+ _{\X, T, Z}$ la catégorie des
($F$-)$\D ^\dag _{\X, \Q}$-modules surholonomes $\E$ à support
dans $Z$ et vérifiant $\R \underline{\Gamma} ^\dag _T (\E) =0$.

En utilisant le triangle de localisation \ref{tri-local} en $T$,
on remarque que la condition $\R \underline{\Gamma} ^\dag _T (\E) =0$ équivaut à dire
que le morphisme canonique $\E \rightarrow (\hdag T) (\E) $ est un isomorphisme.
Lorsque $T$ est un diviseur de $X$,
cette condition est équivalente à demander que la structure de $\D ^\dag _{\X, \Q}$-module de $\E$
se prolonge en une structure de $\D ^\dag _{\X} (\hdag T) _{ \Q}$-module.

De même, on note ($F$-)$\mathfrak{C} ^+ _{\X, T, Z}$ la sous-catégorie
pleine de ($F$-)$D ^{\mathrm{b}} _{\mathrm{surhol}} (\D ^\dag
_{\X, \Q})$ des complexes $\E$ vérifiant $(\hdag Z)(\E ) =0$ et
$\R \underline{\Gamma} ^\dag _T (\E) =0$.

\end{nota}

\begin{vide}\label{notation-surcohhol}
On garde les notations \ref{notations-catzt}.
Dans \cite[3.2.1]{caro_surcoherent}, on a défini les catégories analogues à celles de \ref{notations-catzt}
en remplaçant
{\og surholonome\fg} par {\og surcohérent\fg}. On les avait désignées par
($F$-)$\mathfrak{M}  _{\X, T, Z}$ et ($F$-)$\mathfrak{C}  _{\X, T, Z}$.
Comme la surholonomie est une notion plus forte que la surcohérence, on obtient les inclusions :
($F$-)$\mathfrak{M} ^+  _{\X, T, Z}\subset $($F$-)$\mathfrak{M}  _{\X, T, Z}$
et ($F$-)$\mathfrak{C} ^+  _{\X, T, Z}\subset$($F$-)$\mathfrak{C}  _{\X, T, Z}$.

 De plus, si $T$ est un {\it diviseur} et si $\X$ est {\it propre},
 la sous-catégorie pleine
  de $(F\text{-})D ^\mathrm{b} _{\mathrm{coh}} (\D ^\dag _{\X  } (\hdag T ) _{ \Q}) $
  des ($F$-)complexes $\E  $ tels que
  $\DD _{T }  (\E ) \in (F\text{-})\mathfrak{C} _{\X , T , Z }$ avait été notée
  $(F\text{-})\mathfrak{C} ^\vee _{\X , T , Z }$ (voir \cite[3.3.1]{caro_surcoherent}).
  Comme le foncteur $\DD _T$ stabilise la catégorie $(F\text{-})\mathfrak{C} ^+ _{\X, T, Z} $,
 on obtient l'inclusion
  $(F\text{-})\mathfrak{C} ^+ _{\X, T, Z}\subset (F\text{-})\mathfrak{C} ^\vee _{\X, T, Z}$.
  \end{vide}

\begin{rema}
[Frobenius]
\label{remafrob}
  Dans toutes les assertions de cette section, les structures de Frobenius sont inutiles.
  Plus précisément, on peut remplacer {\og $F$-\fg} par {\og ($F$-)\fg}.
  Afin de ne pas alourdir les notations, nous avons parfois omis d'indiquer les parenthèses.
\end{rema}

\begin{defi}
  Avec les notations \ref{notations-catzt}, on définit
  {\og le foncteur dual sur $F\text{-}\mathfrak{C} ^+  _{\X, T, Z} $\fg}
  en posant $\DD _T:= (\hdag T) \circ \DD$ :
  $F\text{-}\mathfrak{C} ^+  _{\X, T, Z} \rightarrow F\text{-}\mathfrak{C} ^+  _{\X, T, Z}$,
  où $\DD$ est le foncteur dual $\D ^\dag _{\X ,\Q}$-linéaire de \ref{def-dualdiv} (i.e., le diviseur est vide).

  Lorsque $T$ est un diviseur, on vérifie de manière analogue à la preuve de \ref{stab-surhol-dualT} l'isomorphisme
$\DD _T (\E) \riso (\hdag T)\circ \DD (\E)$, où
$\DD$ (resp. $\DD _T$) est le foncteur dual $\D ^\dag _{\X ,\Q}$-linéaire
(resp. $\D ^\dag _{\X} (\hdag T) _{\Q}$-linéaire)
de \ref{def-dualdiv}.
Cela unifie donc nos notations.
\end{defi}

\begin{prop}
\label{bidua+rel}
Avec les notations de \ref{notations-catzt},
soit $\E \in \mathfrak{C} ^+  _{\X, T, Z} $.
  On bénéficie de l'isomorphisme canonique
$\iota \ :\ \E \riso \DD _T \circ \DD _T (\E)$, dit de {\og bidualité\fg},
compatible à Frobenius.

De plus, soit $f$ : $\Y \rightarrow \X$ un morphisme propre de $\V$-schémas formels lisses.
On dispose alors de l'isomorphisme
$f _+ \circ \DD _{f _0 ^{-1} (T )} (\E) \riso \DD _{T} \circ f _+  (\E)$
dit de {\og dualité relative\fg}.
\end{prop}
\begin{proof}
On dispose du morphisme canonique $\DD (\E) \rightarrow (\hdag T) \DD (\E)$ désigné par $u$.
En notant $\E _1 =(\hdag T) \DD (\hdag T) \DD (\E)$, $\E _2=(\hdag T) \DD \DD (\E)$,
on obtient le morphisme
$v:=(\hdag T)\DD (u)\,:\, \E _2 \rightarrow \E _2$.
Dans un premier temps, vérifions que $v$ est un isomorphisme.

Soit $T' \supset T$ un diviseur de $X$ contenant $T$.
Par \ref{hdagcirchdag},
$(\hdag T') (\hdag T) \DD \riso (\hdag T') \DD $.
Or, d'après \cite[I.4.4]{virrion}, le foncteur dual commute à l'extension des scalaires
et donc $(\hdag T') \DD \riso \DD _{T'} (\hdag T')$.
On obtient donc par composition :
$(\hdag T') (\hdag T) \DD \riso \DD _{T'} (\hdag T')$.
Ainsi,
$(\hdag T') (v)\riso \DD _{T'} (\hdag T') (u) $.
Par \ref{hdagcirchdag}, $(\hdag T')(u)$ est un isomorphisme.
Il en résulte alors que
$(\hdag T') (v)$ est un isomorphisme.

Soient $T _1, \dots , T_r$ des diviseurs de $X$ tels que $T \subset \cap _{i=1} ^r T_i$.
Par récurrence sur $r$, vérifions maintenant que $(\hdag \cap _{i=1} ^r T_i) (v)$ est un isomorphisme.
Le cas où $r =1$ vient d'être traité. Supposons $r \geq 2$.
Les triangles distingués de Mayer-Vietoris (voir \ref{eq1mayer-vietoris}) utilisés
pour $T _1$ et $T'_1:=\cap _{i=2} ^r T_i$ donne par fonctorialité :
\begin{equation}
\label{MVxymat}
\xymatrix {
  {(\hdag T _1 \cap T '_1)(\E _1 ) }
  \ar[r] \ar[d] ^{(\hdag T _1 \cap T '_1)(v)}
&
{(\hdag T _1)(\E _1 ) \oplus  (\hdag T '_1) (\E _1 )}
\ar[r] \ar[d] ^{(\hdag T _1)(v ) \oplus  (\hdag T '_1) (v)}
&
{ (\hdag T _1 \cup T '_1)(\E _1 ) }
\ar[r] \ar[d] ^{(\hdag T _1 \cup T '_1)(v)}
&
{(\hdag T _1\cap T '_1)(\E _1 )[1]}
\ar[d] ^{(\hdag T _1 \cap T '_1)(v)}
\\
  {(\hdag T _1 \cap T '_1)(\E _2 ) }
  \ar[r]
&
{(\hdag T _1)(\E _2 ) \oplus  (\hdag T '_1) (\E _2 )}
\ar[r]
&
{ (\hdag T _1 \cup T '_1)(\E _2 ) }
\ar[r]
&
{(\hdag T _1\cap T '_1)(\E _2 )[1].}
}
\end{equation}
Par hypothèse de récurrence,
$(\hdag T _1)(v ) \oplus  (\hdag T '_1) (v)$ et
$(\hdag T _1 \cup T '_1)(v)$ sont des isomorphismes.
Via \ref{MVxymat}, il en résulte que $(\hdag \cap _{i=1} ^r T_i) (v)$ est un isomorphisme.

Comme $T$ est une intersection de diviseurs, il en résulte que
$(\hdag T) (v)$ est un isomorphisme.
Or, d'après \ref{hdagcirchdag}, $\E _1 \riso (\hdag T) (\E _1)$, $\E _2 \riso (\hdag T) (\E _2)$.
D'où : $v \riso (\hdag T) (v)$.
Ainsi $v$ est un isomorphisme.
On construit alors l'isomorphisme de bidualité $\iota$ en posant :
$$\iota \ :\ (\hdag T) \DD (\hdag T) \DD (\E) \underset{v}{\riso} (\hdag T) \DD \DD (\E) \liso (\hdag T) \E \liso \E,$$
où le deuxième isomorphisme se déduit de l'isomorphisme de bidualité de Virrion (voir
\ref{bidual}) et le dernier de \ref{hdagcirchdag}.

Avec l'isomorphisme de dualité relative (voir \ref{iso-dualrel})
et la commutation de l'image directe aux foncteurs localisations (voir \ref{commutfonctcohlocal2+}),
on construit l'isomorphisme de dualité relative via le composé :
\begin{equation}
  \notag
  f _+ \circ \DD _{f _0 ^{-1} (T )} (\E) =
f _+
\circ (\hdag f _0 ^{-1} (T ))
\circ
\DD (\E)
\underset{\ref{commutfonctcohlocal2+}}{\riso}
(\hdag T )
\circ
f _+
\circ
\DD (\E)
\underset{\ref{iso-dualrel}}{\riso}
(\hdag T )
\circ
\DD
\circ
f _+ (\E)
=
\DD _{T} \circ f _+  (\E).
\end{equation}

\end{proof}

\begin{prop}\label{propindX}
  Soient $\X$ un $\V$-schéma formel lisse, $Z$, $Z'$, $T $ et $T'$ des sous-schémas fermés de $X$
  tels que $Z \setminus T =Z'\setminus T'$.
  On bénéficie alors des égalités $F\text{-}\mathfrak{M} ^+ _{\X, T, Z}= F\text{-}\mathfrak{M} ^+ _{\X, T', Z'}$
  et $F\text{-}\mathfrak{C} ^+ _{\X, T, Z}= F\text{-}\mathfrak{C} ^+ _{\X, T', Z'}$.
  En particulier, en notant $\overline{U}$ l'adhérence schématique de $U :=Z \setminus T$ dans $X$,
  $F\text{-}\mathfrak{M} ^+ _{\X, T, Z}=F\text{-}\mathfrak{M} ^+ _{\X, T, \overline{U}}$
  et
  $F\text{-}\mathfrak{C} ^+ _{\X, T, Z}=F\text{-}\mathfrak{C} ^+ _{\X, T, \overline{U}}$.
\end{prop}

\begin{proof}
  D'après \cite[3.2.4]{caro_surcoherent},
  on dispose des égalités $F\text{-}\mathfrak{M}  _{\X, T, Z}= F\text{-}\mathfrak{M}  _{\X, T', Z'}$
  et $F\text{-}\mathfrak{C} _{\X, T, Z}= F\text{-}\mathfrak{C}  _{\X, T', Z'}$.
Les inclusions
  $F\text{-}\mathfrak{M} ^+  _{\X, T, Z}\subset F\text{-}\mathfrak{M}  _{\X, T, Z}$,
  $F\text{-}\mathfrak{C} ^+  _{\X, T, Z}\subset F\text{-}\mathfrak{C}  _{\X, T, Z}$
  (voir \ref{notation-surcohhol})
  nous permettent alors de conclure.
\end{proof}

\begin{lemm}
  \label{transf12!Gamma}
Soient $f _1\,:\,\X _2 \rightarrow \X _1$, $f _2\,:\,\X _3 \rightarrow \X _2$
deux morphismes de $\V$-schémas formels lisses, $Z _1$ un sous-schéma fermé de $X _1$.
On pose $Z _2:= f _1 ^{-1} (Z_1)$, $Z _3:= f _2 ^{-1} (Z_2)$.

Les diagrammes canoniques de foncteurs
$F\text{-}\mathfrak{C} ^+  _{\X_1, T_1, Z_1} \rightarrow
F\text{-}\mathfrak{C} ^+  _{\X_3, T_3, Z_3} $ :
\begin{equation}
  \label{transf12!Gamma-diag}
  \xymatrix {
  {(f _1 \circ f _2) ^! \R \underline{\Gamma} ^\dag _{Z _1}}
  \ar[d] _-\sim
  \ar[rr] _-\sim
  &&
  {\R \underline{\Gamma} ^\dag _{Z _3}(f _1 \circ f _2) ^! }
  \ar[d] _-\sim
  \\
  {f _2 ^! f _1 ^! \R \underline{\Gamma} ^\dag _{Z _1}}
  \ar[r] _-\sim
  &
  { f _2 ^! \R \underline{\Gamma} ^\dag _{Z _2}f _1 ^! }
  \ar[r] _-\sim
  &
  {\R \underline{\Gamma} ^\dag _{Z _3} f _2 ^! f _1 ^! ,}
}
\xymatrix {
  {(f _1 \circ f _2) ^! (\hdag Z _1)}
  \ar[d] _-\sim
  \ar[rr] _-\sim
  &&
  {(\hdag Z _3) (f _1 \circ f _2) ^! }
  \ar[d] _-\sim
  \\
  {f _2 ^! f _1 ^! (\hdag Z _3) }
  \ar[r] _-\sim
  &
  { f _2 ^! (\hdag Z _2) f _1 ^! }
  \ar[r] _-\sim
  &
  {(\hdag Z _3) f _2 ^! f _1 ^! ,}
}
\end{equation}
où les isomorphismes horizontaux sont les isomorphismes de commutation à l'image inverse extraordinaire
du foncteur cohomologique local ou du foncteur de localisation
 (voir \ref{commutfonctcohlocal2}),
est commutatif.
\end{lemm}

\begin{proof}
Le cas du foncteur de localisation se traitant de manière identique, contentons-nous de prouver
le cas du foncteur cohomologique local.
Rappelons d'abord que le morphisme du haut à gauche de \ref{transf12!Gamma-diag} correspond par construction
(plus précisément, voir la construction dans la première partie de la preuve de \cite[2.2.18.1]{caro_surcoherent})
au composé des morphismes canoniques :
$$(f _1 \circ f _2) ^! \R \underline{\Gamma} ^\dag _{Z _1}
\liso
\R \underline{\Gamma} ^\dag _{Z _3}(f _1 \circ f _2) ^! \R \underline{\Gamma} ^\dag _{Z _1}
\riso
\R \underline{\Gamma} ^\dag _{Z _3}(f _1 \circ f _2) ^! .$$
Il s'agit d'établir que le rectangle du fond (et vertical) du diagramme ci-dessous
\begin{equation}
\label{diag1defg+g!hautbas}
  \xymatrix@C=0,1cm @R=0,3cm{
  { f _2 ^! f _1 ^! \R \underline{\Gamma} ^\dag _{Z _1}}
  \ar[rr] _-\sim
  \ar[d] _-\sim
  &
  &
  {(f _1 \circ f _2) ^! \R \underline{\Gamma} ^\dag _{Z _1}}
  \ar[d] _-\sim
  \\
  { \R \underline{\Gamma} ^\dag _{Z _3} f _2 ^! f _1 ^! }
  \ar[rr] _-\sim
  &
  &
  {\R \underline{\Gamma} ^\dag _{Z _3}(f _1 \circ f _2) ^! }
  \\
  &
  { \R \underline{\Gamma} ^\dag _{Z _3} f _2 ^! f _1 ^! \R \underline{\Gamma} ^\dag _{Z _1}}
  \ar[lu] _-\sim
  \ar[rr] _-\sim
  \ar[luu] _-\sim
  &
  &
  {\R \underline{\Gamma} ^\dag _{Z _3}(f _1 \circ f _2) ^! \R \underline{\Gamma} ^\dag _{Z _1}}
  \ar[lu] _-\sim
  \ar[luu] _-\sim
 }
\end{equation}
est commutatif. Le rectangle horizontal et celui oblique sont commutatifs par fonctorialité,
tandis que le triangle de droite l'est par définition. Toutes les flèches étant des isomorphismes,
il est alors équivalent de prouver que le triangle de gauche est commutatif.
Considérons pour cela le diagramme canonique suivant :
\begin{equation}
\label{diag2defg+g!hautbas}
\xymatrix@C=0,5cm {
{ \R \underline{\Gamma} ^\dag _{Z  _3} f _2 ^! f _1 ^! \R \underline{\Gamma} ^\dag _{Z _1}}
\ar[d] _-\sim
&
&
{ \R \underline{\Gamma} ^\dag _{Z  _3} f _2 ^!
\R \underline{\Gamma} ^\dag _{Z  _2}f _1 ^! \R \underline{\Gamma} ^\dag _{Z _1}}
\ar[ll] _-\sim
\ar[rr] _-\sim
\ar[rd] _-\sim
\ar[ld] _-\sim
&
&
{ \R \underline{\Gamma} ^\dag _{Z  _3} f _2 ^! f _1 ^! \R \underline{\Gamma} ^\dag _{Z _1}}
\ar[d] _-\sim
\\
{ f _2 ^! f _1 ^! \R \underline{\Gamma} ^\dag _{Z _1}}
&
{f _2 ^! \R \underline{\Gamma} ^\dag _{Z  _2}f _1 ^! \R \underline{\Gamma} ^\dag _{Z _1}}
\ar[r] _-\sim
\ar[l] _-\sim
&
{f _2 ^! \R \underline{\Gamma} ^\dag _{Z  _2}f _1 ^!}
&
{\R\underline{\Gamma} ^\dag _{Z  _3} f _2 ^!\R \underline{\Gamma} ^\dag _{Z  _2}f _1 ^!}
\ar[r] _-\sim
\ar[l] _-\sim
&
{\R\underline{\Gamma} ^\dag _{Z  _3} f _2 ^! f _1 ^!.}
}
\end{equation}
Les trois petits {\og carrés\fg} de ce diagramme sont commutatifs par fonctorialité. Par exemple,
celui du milieu l'est par fonctorialité du morphisme canonique
$\R\underline{\Gamma} ^\dag _{Z  _3} f _2 ^!\R \underline{\Gamma} ^\dag _{Z  _2}f _1 ^!
\riso f _2 ^! \R \underline{\Gamma} ^\dag _{Z  _2}f _1 ^!$
via $\R \underline{\Gamma} ^\dag _{Z _1} \rightarrow Id$.
Le composé du haut de \ref{diag2defg+g!hautbas} est l'identité de
$\R \underline{\Gamma} ^\dag _{Z  _3} f _2 ^! f _1 ^! \R \underline{\Gamma} ^\dag _{Z _1}$.
Le composé du bas de \ref{diag2defg+g!hautbas} est
l'isomorphisme canonique $ f _2 ^! f _1 ^! \R \underline{\Gamma} ^\dag _{Z _1}
\riso
\R\underline{\Gamma} ^\dag _{Z  _3} f _2 ^! f _1 ^!$
du triangle de gauche de \ref{diag1defg+g!hautbas}.
La commutativité de \ref{diag2defg+g!hautbas}
se traduit donc par celle du triangle de gauche de
\ref{diag1defg+g!hautbas}.
\end{proof}

\begin{theo}
[Opérations cohomologiques sur les schémas]
\label{defg+g!hautbas}
  Soient $f _1$ : $\X _2 \rightarrow \X _1$ un morphisme propre de $\V$-schémas formels lisses,
  pour $i=1,\,2$,
  $T _i$
 et $Z _i$ des sous-schémas fermés de $X _i$ et $U _i := Z _i \setminus T _i$.
  On suppose que le morphisme $f _1$ induit un morphisme
  $g _1$ : $U _2 \rightarrow U _1$.

  Les foncteurs $f _{1+} \text{ et }\DD _{T _1}\circ  f _{1+} \circ \DD _{T _2}$
  induisent alors des foncteurs
  $F\text{-}\mathfrak{C} ^+ _{\X _2 , T _2 , Z _2 } \rightarrow
  F\text{-}\mathfrak{C} ^+ _{\X _1 , T _1 , Z _1 }.$
  De plus, les foncteurs
  $\R \underline{\Gamma} ^\dag _{Z _2}  \circ  (\hdag T _2)\circ   f _1^!
  \text{ et }
  \DD _{T _2}\circ \R \underline{\Gamma} ^\dag _{Z _2} \circ  (\hdag T _2)\circ f _1^! \circ \DD _{T _1}$
  induisent des foncteurs
  $F\text{-}\mathfrak{C} ^+ _{\X _1 , T _1 , Z _1 } \rightarrow
  F\text{-}\mathfrak{C} ^+ _{\X _2 , T _2 , Z _2 }.$
  Ces foncteurs seront notés abusivement
  $g _{1+}$, $g _{1!}$, $g _1 ^!$ et $g _1 ^+$
  (en fait, ces notations ne sont pas toujours abusives : voir la remarque \ref{nota-abus}).
  On les appellera respectivement {\og image directe\fg},
  {\og image directe extraordinaire\fg}, {\og image inverse extraordinaire\fg} et {\og image inverse\fg} par $g _1$.

De plus, soient
$f _2$ : $\X _3 \rightarrow \X _2$ (resp. $f _3$ : $\X _4 \rightarrow \X _3$)
un morphisme propre de $\V$-schémas formels lisses et,
  pour $i=3,\,4$,
  $T _i$ et $Z _i$ deux sous-schémas fermés de $X _i$ et $U _i := Z _i \setminus T _i$.
  On suppose que le morphisme $f _2$ (resp. $f _3$) induit un morphisme
  $g _2$ : $U _3 \rightarrow U _2$ (resp. $g _3$ : $U _4 \rightarrow U _3$).

  On dispose alors de l'isomorphisme de transitivité canonique $g _2 ^! \circ g _1 ^! \riso (g _1 \circ g _2 ) ^!$.
  Celui-ci vérifie la condition d'associativité : les deux isomorphismes composés
  $g _3 ^! \circ g _2 ^! \circ g _1 ^! \riso g _3 ^! \circ (g _1 \circ g _2 ) ^! \riso (g _1 \circ g _2 \circ g _3) ^!$
  et
  $g _3 ^! \circ g _2 ^! \circ g _1 ^! \riso  (g _2 \circ g _3 ) \circ  g _1 ^! \riso (g _1 \circ g _2 \circ g _3) ^!$
  sont identiques.

  De même, pour les foncteurs images inverses, images directes ou images directes extraordinaires,
  on bénéficie d'isomorphismes de transitivité canoniques analogues,
  ces isomorphismes vérifiant la condition d'associativité correspondante.
\end{theo}

\begin{proof}Commençons par prouver que les factorisations sont valables.
Soient $\E _1 \in F\text{-}\mathfrak{C} ^+ _{\X _1 , T _1 , Z _1 }$,
$\E _2 \in F\text{-}\mathfrak{C} ^+ _{\X _2 , T _2 , Z _2 }$.
Par \ref{stabinmdir}, $f _{1+} (\E _2) $ est un complexe surholonome.
  D'après \ref{propindX}, on peut supposer que $Z _2$ est l'adhérence schématique
  de $U _2$ dans $X _2$. Il en résulte l'inclusion $Z _2 \subset f _1 ^{-1} (Z _1)$.
  Le complexe $f _{1+} (\E _2)$ est donc à support dans $Z _1$.
  En outre, puisque $U _2 \cap f_1 ^{-1} (T _1)= \emptyset$,
  grâce à \ref{propindX}, quitte à remplacer $T _2 $ par
  $T _2 \cup f_1 ^{-1} (T _1)$, on se ramène au cas où $T _2 $ contient $f_1 ^{-1} (T _1)$.
Ainsi
$\R \underline{\Gamma} ^\dag _{T _1} \circ f _{1+} (\E_2)
\riso
f _{1+} \circ \R \underline{\Gamma} ^\dag _{f ^{-1}(T _1)}  (\E_2)
\riso 0$
(voir \ref{commutfonctcohlocal2+} pour le premier isomorphisme
et \ref{GammaZGammaZ'} pour le second).
Nous avons ainsi vérifié que
$f _{1+} (\E _2) \in F\text{-}\mathfrak{C} ^+ _{\X _1 , T _1 , Z _1 }$.
  Puisque le foncteur $\DD _{T _1}$ (resp. $\DD _{T _2}$)
  préserve $F\text{-}\mathfrak{C} ^+ _{\X _1 , T _1 , Z _1 } $
  (resp. $F\text{-}\mathfrak{C} ^+ _{\X _2 , T _2 , Z _2 } $),
  on en déduit que $\DD _{T _1} \circ f _{1+} \circ \DD _{T _2}(\E _2)
  \in F\text{-}\mathfrak{C} ^+ _{\X _1 , T _1 , Z _1 }$.
  De manière analogue, via la stabilité de la surholonomie par image inverse extraordinaire et
  foncteur cohomologique local (voir \ref{stab-surhol-Dual-hdag}.\ref{stab-surhol-Dual-hdag1}, \ref{stabinminv}),
  on établit que
  $\R \underline{\Gamma} ^\dag _{Z _2}  \circ  (\hdag T _2)\circ   f _1^! (\E _1),
  \DD _{T _2}\circ \R \underline{\Gamma} ^\dag _{Z _2} \circ  (\hdag T _2)\circ f _1^! \circ \DD _{T _1} (\E _1)
  \in F\text{-}\mathfrak{C} ^+ _{\X _2 , T _2 , Z _2 } $.
\\

Concernant la transitivité du foncteur image directe par $g_1$
et la condition d'associativité correspondante, cela découle
par définition de \ref{compatsurhol}.
Avec l'isomorphisme de bidualité (voir \ref{bidua+rel}),
on en déduit que les propriétés analogues sont aussi satisfaites pour le foncteur image directe
extraordinaire par $g _1$.

Construisons maintenant l'isomorphisme de transitivité du foncteur images inverses extraordinaires.
Par commutation aux images inverses extraordinaires
des foncteurs cohomologiques locaux à support strict et des foncteurs de localisation
(voir \ref{commutfonctcohlocal2}) et des propriétés usuelles de composition
(voir \ref{GammaZGammaZ'}, \ref{hdagcirchdag}), on obtient :
$$(\R \underline{\Gamma} ^\dag _{Z _3}  \circ  (\hdag T _3)\circ   f _2^!)
\circ
(\R \underline{\Gamma} ^\dag _{Z _2}  \circ  (\hdag T _2)\circ   f _1^!)
\riso
\R \underline{\Gamma} ^\dag _{Z _3\cap f _2 ^{-1} (Z _2)}
  \circ  (\hdag T _3 \cup f _2 ^{-1} (T _2)) \circ
  f _2 ^! \circ f _1 ^!.$$
D'après \ref{trans-f!pref+}, on bénéficie de l'isomorphisme de transitivité
$f _2 ^! \circ f _1 ^! \riso(f _1 \circ f _2) ^!$.
De plus, comme $Z _3 \setminus T _3 =Z _3\cap f _2 ^{-1} (Z _2) \setminus (T _3 \cup f _2 ^{-1} (T _2))$,
par \ref{propindX}, \ref{GammaZGammaZ'} et \ref{hdagcirchdag}, on en déduit :
$$(\R \underline{\Gamma} ^\dag _{Z _3}  \circ  (\hdag T _3)\circ   f _2^! )
\circ ( \R \underline{\Gamma} ^\dag _{Z _2}  \circ  (\hdag T _2)\circ   f _1^! )
\riso
\R \underline{\Gamma} ^\dag _{Z _3}
  \circ  (\hdag T _3 ) \circ
  (f _1 \circ f _2) ^!,$$
ce qui donne par définition l'isomorphisme canonique :
$g _2 ^! \circ g _1 ^! \riso (g _1 \circ g _2 ) ^!$.
La condition d'associativité de cet isomorphisme résulte des différents cas particuliers
\ref{trans-f!pref+}, \ref{transZZ'Z''}, \ref{transf12!Gamma} (il s'agit de vérifier
la commutativité d'un diagramme que l'on décompose en carrés élémentaires commutatifs grâce aux différents cas
particuliers).

Enfin, on déduit du cas de l'image inverse extraordinaire, via l'isomorphisme de bidualité (voir \ref{bidua+rel}),
que les propriétés analogues sont aussi satisfaites pour le foncteur image inverse par $g _1$.

\end{proof}

\begin{rema}
  \label{nota-abus}
  Avec les notations \ref{defg+g!hautbas},
  d'après la proposition \ref{def-surhol-rele2}, lorsque $\X _1$ et $\X _2$ sont propres,
  les foncteurs $g _{1+}$, $g _{1!}$, $g _1 ^!$ et $g _1 ^+$ ne dépendent que de $g _1$.
  Avec cette hypothèse de propreté, les notations $g _{1+}$, $g _{1!}$, $g _1 ^!$ et $g _1 ^+$
  de \ref{defg+g!hautbas} ne sont donc plus abusives.
\end{rema}

\begin{theo}
\label{theoindpdt1}
  Soient $f$ : $\X _2 \rightarrow \X _1$ un morphisme propre et lisse de $\V$-schémas formels lisses et,
  pour $i=1,\,2$,
  $T _i$ un diviseur de $X _i$ et $Z _i$ un sous-schéma fermé de $X _i$.
  On suppose de plus que $f$ induise l'isomorphisme
  $Z _2\setminus T _2 \riso Z _1 \setminus T _1$.
  On se donne $\E _1 \in (F\text{-})\mathfrak{M} ^+ _{\X _1 , T _1 , Z _1}$ et
$\E _2 \in (F\text{-})\mathfrak{M} ^+ _{\X _2, T _2, Z _2}$.
  Alors,
  \begin{enumerate}
    \item \label{1theoindpdt1}Pour tout $l \neq 0$, $\mathcal{H} ^l (f _+) (\E _2) =0$ et
   $\mathcal{H} ^l  ((\hdag T _2) \circ \R \underline{\Gamma} ^\dag _{Z _2} \circ  f ^!(\E _1))=0$ ;
    \item On dispose d'isomorphismes canoniques
    $f _+\circ (\hdag T _2) \circ \R \underline{\Gamma} ^\dag _{Z _2} \circ  f ^!(\E _1) \riso \E _1$
    et
    $\E_2 \riso (\hdag T _2) \circ  \R \underline{\Gamma} ^\dag _{Z _2} \circ  f ^!\circ f _+ (\E _2)$.
  \end{enumerate}
  Les foncteurs $f_+$ (i.e., avec les notations de \ref{defg+g!hautbas}, $Id _+$)
  et
  $\R \underline{\Gamma} ^\dag _{Z _2} \circ f ^!$ (i.e. $Id ^!$)
  induisent donc
  des équivalences quasi-inverses entre les catégories
  $(F\text{-})\mathfrak{M} ^+ _{\X _2, T _2, Z _2}$ et $(F\text{-})\mathfrak{M} ^+ _{\X _1, T _1, Z _1}$.
  Lorsque $Z _2$ est l'adhérence schématique de $Z _2 \setminus T_2$ dans $X _2$,
  le foncteur $(\hdag T _2)$ est inutile.

On dispose des mêmes résultats en remplaçant {\og $\mathfrak{M} ^+$\fg} par {\og $\mathfrak{M} $\fg}.
\end{theo}

\begin{proof}
En notant $Z ' _2$ l'adhérence schématique de
$Z _2\setminus T _2$ dans $X_2$,
comme $Z ' _2 \setminus f ^{-1} _0 (T _1) = Z _2 \setminus T _2$,
grâce à \ref{propindX} (et \ref{GammaZGammaZ'}, \ref{hdagcirchdag}), on vérifie que pour
tout $\E _1 \in (F\text{-})\mathfrak{M} _{\X _1 , T _1 , Z _1}$
(et donc pour tout $\E _1 \in (F\text{-})\mathfrak{M} ^+ _{\X _1 , T _1 , Z _1}$)
le morphisme canonique
\begin{equation}
  \label{nonhdag.3.2.6}
  \R \underline{\Gamma} ^\dag _{Z '_2} \circ  f ^!(\E _1)
\rightarrow
(\hdag T _2) \circ \R \underline{\Gamma} ^\dag _{Z _2} \circ  f ^!(\E _1)
\end{equation}
est un isomorphisme.
On se ramène ainsi au cas où $Z _2$ est l'adhérence schématique de $Z _2 \setminus T_2$ dans $X _2$.
Dans ce cas, le foncteur $(\hdag T _2)$  est bien inutile.
Le cas où {\og $\mathfrak{M} ^+$\fg} est remplacé par {\og $\mathfrak{M} $\fg} a été traité
dans \cite[3.2.6]{caro_surcoherent} (l'énoncé de \cite[3.2.6]{caro_surcoherent} est correcte
lorsque $Z _2$ est l'adhérence schématique de $Z _2 \setminus T_2$ dans $X _2$, autrement il faut rajouter
le foncteur $(\hdag T _2) $). Rappelons que d'après \cite[3.2.2.1]{caro_surcoherent},
le théorème \cite[3.2.6]{caro_surcoherent} reste valable sans structure de Frobenius.
Le cas non respectif, i.e. celui concernant {\og $\mathfrak{M} ^+$\fg}, s'en déduit aussitôt grâce à la stabilité
de la surholonomie par image inverse extraordinaire et image directe par un morphisme propre, par foncteur cohomologique
local à support strict et foncteur de localisation
(voir \ref{stab-surhol-Dual-hdag}, \ref{stabinminv}, \ref{stabinmdir}).
\end{proof}

La proposition \ref{j!} ci-dessous nous sera nécessaire dans la procédure de recollement de \ref{defi-rem-sch-sep}
(qui permet d'après \ref{defi2-rem-sch-sep} de construire, pour toute variété $U$,
la catégorie $(F\text{-})\mathfrak{M} ^+ _{U }$).
Cette proposition \ref{j!} est une conséquence du théorème \ref{theoindpdt1} ainsi que
du lemme ci-après.
\begin{lemm}\label{UiU}
  Soient $\X$ un $\V$-schéma formel lisse, $T$ un diviseur de $X$, $Z$ un sous-schéma fermé
  de $X$ et $U := Z \setminus T$. Soient $n$ diviseurs $T _1,\dots, T_n $ de $X$
  tels que $\cup _{i=1,\dots,n} U _i = U$, avec $U _i := Z \setminus T _i$.
  Soit  $\E \in D ^\mathrm{b} _\mathrm{coh} (\D ^\dag _{\X,\Q} )$ à support dans $Z$.
  Si, pour tout $i=1,\dots, n$, $(\hdag T _i) (\E) = 0$ alors $(\hdag T) (\E) =0$.
\end{lemm}

\begin{proof}
  Comme
  $(\hdag T) (\E)  \in D ^\mathrm{b} _\mathrm{coh} (\D ^\dag _{\X} (\hdag T) _{\Q} )$,
  il résulte de \ref{4312Be1} que
$(\hdag T) (\E) =0$ si et seulement si $\E$ est nul en dehors de $T$.
On se ramène ainsi au cas où $T $ est vide.
D'après \cite[2.2.9]{caro_surcoherent},
comme $\E$ est à support dans $Z$,
$\R \underline{\Gamma} ^\dag _{Z }(\E) \riso \E $.
De plus, pour tout $i=1, \dots, n$, le triangle de localisation en $T_i$ (voir \ref{tri-local}) donne
$\R \underline{\Gamma} ^\dag _{T _i }(\E) \riso  \E $.
  Par \ref{GammaZGammaZ'}, on en déduit :
$\R \underline{\Gamma} ^\dag _{Z \cap \cap _{i=1,\dots,n} T _i} (\E)\riso \E$.
  Comme $Z \cap \cap _{i=1,\dots,n} T _i= \emptyset$,
  $\R \underline{\Gamma} ^\dag _{Z \cap \cap _{i=1,\dots,n} T _i} (\E) =0$.
  D'où $\E =0$.
\end{proof}

\begin{prop}\label{j!}
  Soient $f$ : $\X' \rightarrow \X$ un morphisme propre et lisse de $\V$-schémas formels lisses,
  $T$ (resp. $T '$) un diviseur de $X$ (resp. $X'$), $Z$ (resp. $Z'$) un sous-schéma fermé de $X$ (resp. $X'$)
  et $U := Z \setminus T$ (resp. $U' := Z' \setminus T'$).
  On suppose que $f$ induise une immersion ouverte $j$ : $U' \hookrightarrow U$.

Pour tout $\E  \in (F\text{-})\mathfrak{M}  _{\X  , T  , Z }$, pour tout $l \neq 0$,
$\mathcal{H} ^l ( (\hdag T') \R \underline{\Gamma} ^\dag _{Z '} \circ  f ^!(\E ))=0$.
Le foncteur $j ^! := (\hdag T') \R \underline{\Gamma} ^\dag _{Z '} \circ  f ^!$ induit donc le foncteur
$j ^!\,:\,(F\text{-})\mathfrak{M}  _{\X  , T  , Z }\rightarrow (F\text{-})\mathfrak{M}  _{\X ' , T ' , Z '}$.
De même en remplaçant la surcohérence par la surholonomie, i.e., {\og $\mathfrak{M}$\fg}
par {\og $\mathfrak{M}^+$\fg}.
\end{prop}

\begin{proof}
Comme un module surholonome est surcohérent,
le cas surholonome découle immédiatement du cas surcohérent. Traitons donc ce dernier.
Soit $\E  \in (F\text{-})\mathfrak{M}  _{\X  , T  , Z }$.
Comme le morphisme canonique $\E \rightarrow \E (\hdag T)$ est un isomorphisme,
avec \ref{hdagcirchdag} et \ref{commutfonctcohlocal2}, on obtient
l'isomorphisme
$(\hdag T') \R \underline{\Gamma} ^\dag _{Z '} \circ  f ^! (\E)
\riso
(\hdag T' \cup f _0 ^{-1} (T)) \R \underline{\Gamma} ^\dag _{Z '} \circ  f ^! (\E)$.
Comme $Z ' \setminus T'  =Z ' \setminus (T' \cup f _0 ^{-1} (T))$,
par \ref{theoindpdt1},
$(F\text{-})\mathfrak{M}  _{\X ' , T ' , Z '}= (F\text{-})\mathfrak{M}  _{\X ' , T '\cup f _0 ^{-1} (T), Z '}$.
Quitte à remplacer $T'$ par $T' \cup f _0 ^{-1} (T)$,
on peut donc supposer $f _0 ^{-1} (T) \subset T'$.
Par stabilité de la surcohérence par image inverse extraordinaire et foncteur cohomologique local,
$(\hdag T') \R \underline{\Gamma} ^\dag _{Z '} \circ  f ^!(\E ))
\in D ^\mathrm{b} _\mathrm{coh} (\D ^\dag _{\X'} (\hdag T') _{\Q} )$.
On en déduit, d'après \ref{4312Be1}, que
$\mathcal{H} ^l ( (\hdag T') \R \underline{\Gamma} ^\dag _{Z '} \circ  f ^!(\E ))=0$
si et seulement si sa restriction en dehors de $T'$ est nulle
et donc si et seulement si sa restriction en dehors de $f _0 ^{-1} (T)$ est nulle.
  Quitte à remplacer $f$ par le morphisme propre et lisse induit
$\X ' \setminus f _0 ^{-1} (T) \rightarrow \X \setminus T$,
  on se ramène ainsi au cas où $T$ est vide, i.e., $U =Z$.

Comme $j\,:\,U' \hookrightarrow Z$ est une immersion ouverte, $Z \setminus j(U')$ est un sous-schéma fermé de $X$.
  Il existe donc des diviseurs $T _1 ,\dots, T _n$ de $X$
  tels que $\cap _{i=1,\dots ,n} T _i = Z \setminus j(U')$.
  Pour tout $i =1,\dots ,n$, notons $T ' _i := T' \cup f _0 ^{-1} (T _i)$ les diviseurs
  de $X'$ correspondants et $U ' _i := Z '\setminus T '_i$.
L'immersion ouverte $j$ induit donc, pour tout $i =1,\dots ,n$,
l'isomorphisme
$U '_i \riso j( U ') \setminus T _i$.
Or, comme $Z \setminus T _i \subset Z \setminus \cap _{i=1,\dots ,n} T _i = j(U')$,
on obtient $Z \setminus T _i \subset j(U') \setminus T _i $. D'où
$Z \setminus T _i = j(U') \setminus T _i $.
Le morphisme $f$ induit donc l'isomorphisme
$U '_i \riso Z \setminus T _i$.
Ainsi, d'après \ref{theoindpdt1},
$(F\text{-})\mathfrak{M}  _{\X ' , T ' _i, Z '}\cong (F\text{-})\mathfrak{M}  _{\X , T _i, Z }$
et,
pour tout $l \neq 0$,
comme $\E (\hdag T _i) \in (F\text{-})\mathfrak{M}  _{\X , T _i, Z }$ alors
$\mathcal{H} ^l ( (\hdag T '_i) \R \underline{\Gamma} ^\dag _{Z '} \circ  f ^!( \E (\hdag T _i) ))=0$.
  Comme $T '_i$ est un diviseur, le foncteur $(\hdag T ' _i)$ est exact sur la catégorie des modules cohérents
(voir \ref{def-hdagcoh}).
On obtient donc, via aussi
\ref{hdagcirchdag} et \ref{commutfonctcohlocal2}, l'isomorphisme canonique
  $$(\hdag T '_i) (\mathcal{H} ^l (\R \underline{\Gamma} ^\dag _{Z '} \circ  f ^!(\E )))
   \riso
   \mathcal{H} ^l ( (\hdag T '_i) \R \underline{\Gamma} ^\dag _{Z '} \circ  f ^!( \E (\hdag T _i) ))=0.$$
D'après \ref{UiU}, comme $\cup _{i =1,\dots ,n} U '_i =U'$, on obtient alors
$(\hdag T ') (\mathcal{H} ^l (\R \underline{\Gamma} ^\dag _{Z '} \circ  f ^!(\E )))=0$.
Comme $T'$ est un diviseur, il en résulte
$ \mathcal{H} ^l ((\hdag T ') \R \underline{\Gamma} ^\dag _{Z '} \circ  f ^!(\E ))
\riso
(\hdag T ') (\mathcal{H} ^l (\R \underline{\Gamma} ^\dag _{Z '} \circ  f ^!(\E )))=0$.
\end{proof}

\begin{rema}
  \label{rema-j!}
  Dans la proposition \ref{j!}, le fait que $T$ et $T'$ soient des diviseurs est indispensable.
  En effet, voici un contre-exemple : soient $\X$ un $\V$-schéma formel lisse tel que $\dim X \geq 2$,
  $T$ un sous-schéma fermé lisse de $X$ de codimension pure égale à $2$.
  On obtient ainsi une immersion ouverte $j\,:\, X \setminus T \subset X$.
  D'après \ref{proptheo}, $\O _{\X ,\Q} \in F\text{-}\mathfrak{M} ^+ _{\X  , \emptyset  , X }$.
  Avec les notations analogues à \ref{j!}, on pose $j^! (\O _{\X ,\Q} ):=\O _{\X ,\Q} (\hdag T)  $.
  On vérifie grâce au triangle distingué de localisation (voir \ref{tri-local})
  $$\R \underline{\Gamma} ^\dag _T (\O _{\X ,\Q} )
\rightarrow \O _{\X ,\Q}
\rightarrow \O _{\X ,\Q} (\hdag T)  \rightarrow
\R \underline{\Gamma} ^\dag _T (\O _{\X ,\Q} ) [1],$$
que $\mathcal{H} ^1 (\O _{\X ,\Q} (\hdag T) ) \riso
\mathcal{H} ^2 _T (\O _{\X ,\Q}) \not =0$.
En effet, si $u\,:\,\T \hookrightarrow \X$ est un relèvement de $T \hookrightarrow X$
(ce qui est localement sur $\X$ possible, par exemple lorsque $\X$ est affine),
$\mathcal{H} ^2 _T (\O _{\X ,\Q}) \riso u_+ (\O _{\T,\Q})\not =0$.
En d'autres termes, $\mathcal{H} ^1 (j^! (\O _{\X ,\Q} ))\not =0$.
\end{rema}

\begin{vide}
  \label{pre-def-surhol-rele}
  Pour $i =1,2$, soient $\X _i$ un $\V$-schéma formel propre et lisse, $T_i$ un diviseur de $X_i$,
  $Z_i$ un sous-schéma fermé de $X_i$.
  On suppose en outre que $Z _2\setminus T _2=Z _1\setminus T _1$.
Les catégories $(F\text{-})\mathfrak{M} ^+ _{\X _1 , T _1 , Z _1}$
et $(F\text{-})\mathfrak{M} ^+ _{\X _2 , T _2 , Z _2}$ sont alors canoniquement isomorphes.

En effet, de manière analogue à \cite[3.2.9]{caro_surcoherent},
quitte à remplacer $\X _2$ par $\X_1 \times \X _2$ (et utiliser les projections de $\X_1 \times \X _2$ sur $\X _1$ et $\X _2$),
 on se ramène grâce à \ref{theoindpdt1} au cas où il existe un morphisme propre et lisse
$f\,:\, \X _2 \rightarrow \X_1$ induisant l'identité sur $Z _2\setminus T _2$.
Le théorème \ref{theoindpdt1} nous permet dans ces conditions de conclure.

\end{vide}

\begin{defi}
\label{def-surhol-rele}
D'après \ref{pre-def-surhol-rele},
si $U$ est une variété telle qu'il existe
un $\V$-schéma formel $\X$ propre et lisse, un diviseur $T$ de $X$, un sous-schéma fermé $Z$ de $X$
tels que $U= Z \setminus T$,
la catégorie
$(F\text{-})\mathfrak{M} ^+ _{\X  , T , Z }$
ne dépend canoniquement que de $U$.
On la notera donc simplement
$(F\text{-})\mathfrak{M} ^+ _{U }$.

Ses objets seront nommés {\og (F\text{-})$\D _{U}$-modules arithmétiques surholonomes\fg}.
Ils constituent un analogue $p$-adique de la catégorie des
$\D _{U _\C }$-modules holonomes, où $U _\C$ est une variété sur le corps des complexes $\C$.
Nous étendrons dans \ref{rem-sch-sep} la construction de $(F\text{-})\mathfrak{M} ^+ _{U }$
pour une variété $U$ quelconque.
\end{defi}

\begin{theo}
\label{theoindpdt2}
  Soient $f$ : $\X _2 \rightarrow \X _1$ un morphisme propre de $\V$-schémas formels lisses,
  $T _1$ et $Z _1$ deux sous-schémas fermés de $X _1$, $T _2$ et $Z_2$ deux sous-schémas fermés de $X _2$.
De plus, on suppose que $f$ induise l'isomorphisme $Z _2\setminus T _2 \riso Z _1 \setminus T _1$.

Alors,
 les foncteurs $f_+$ et
 $(\hdag T _2)\circ   \R \underline{\Gamma} ^\dag _{Z _2} \circ f ^!$ induisent
   des équivalences quasi-inverses entre les catégories
 $(F\text{-})\mathfrak{C} ^+ _{\X _2, T _2, Z _2}$ et $(F\text{-})\mathfrak{C} ^+ _{\X _1, T _1, Z _1}$.
 Lorsque $Z _2$ est l'adhérence schématique de $Z _2 \setminus T_2$ dans $X _2$,
  le foncteur $(\hdag T _2)$  est inutile.
On dispose des mêmes résultats en remplaçant {\og $\mathfrak{C} ^+$\fg} par {\og $\mathfrak{C} $\fg}.
\end{theo}

\begin{proof}
L'isomorphisme \ref{nonhdag.3.2.6} (de la preuve de \ref{theoindpdt1})
reste valable pour $\E _1 \in (F\text{-})\mathfrak{C} _{\X _1 , T _1 , Z _1}$.
On se ramène alors (de manière analogue à \ref{theoindpdt1}) au cas où
$Z _2$ est l'adhérence schématique de $Z _2 \setminus T_2$ dans $X _2$.
Dans ce cas, le foncteur $(\hdag T _2)$ est inutile et
le théorème a alors déjà été traité pour les complexes surcohérents
dans \cite[3.2.8]{caro_surcoherent}.
Le cas concernant {\og $\mathfrak{C} ^+$\fg} (à la place de {\og $\mathfrak{C}$\fg})
en résulte grâce à la stabilité
de la surholonomie.
\end{proof}

\begin{rema}
  \label{rema-theoind}
La différence entre les deux théorèmes \ref{theoindpdt1} et \ref{theoindpdt2} est que le premier
concerne les {\og modules\fg} tandis que le second traite des {\og complexes\fg}.
Pour que les foncteurs que nous utilisons dans ces deux théorèmes
(i.e., ceux de la forme $f _+$ et $(\hdag T _2) \circ \R \underline{\Gamma} ^\dag _{Z _2} \circ  f ^!$)
préservent bien le fait d'être un module (i.e., la condition \ref{theoindpdt1}.\ref{1theoindpdt1}
doit être vérifiée),
les hypothèses de \ref{theoindpdt1} sont alors plus fortes que celles de \ref{theoindpdt2} :
le morphisme $f $ est lisse en plus d'être propre
et les sous-schémas fermés $T _1$ et $T _2$ sont en outre des diviseurs.
Le fait que $T _1$ et $T _2$ soient des diviseurs permet d'utiliser \ref{4312Be1}, ce qui
est techniquement très appréciable
(dans la preuve de \cite[3.2.6]{caro_surcoherent}, cela permet de traiter assez facilement
le cas où $Z _2 \setminus T _2$ est lisse).
De plus, l'hypothèse de lissité de $f $ permet de conserver des diviseurs, i.e., le sous-schéma fermé
$f  ^{-1} (T _1)$ est alors un diviseur de $X _2$.
Enfin, le contre-exemple de \ref{rema-j!} indique que l'hypothèse que $T _1$ et $T _2$ soient des diviseurs
permet (dans le cas des modules) de travailler dans un cadre plus adéquat.
\end{rema}

\begin{defi}
  \label{def-surhol-relebis}
Soit $U$ une variété telle qu'il existe un $\V$-schéma formel $\X$ propre et lisse,
deux sous-schémas fermés $Z,\,T$ de $X$ tels que $U= Z \setminus T$.
De façon similaire à \ref{def-surhol-rele} ou à \cite[3.2.9]{caro_surcoherent},
il découle de \ref{theoindpdt2} que la catégorie
$(F\text{-})\mathfrak{C} ^+ _{\X  ,\, T  ,\, Z }$ est indépendante de $\X $, $T $ et $Z $,
ne dépend que de $U $ et sera désignée par $(F\text{-})D ^\mathrm{b} _\mathrm{surhol} (\D _{U })$.
Ses objets sont appelés {\og ($F$-)complexes de $\D$-modules arithmétiques surholonomes sur $U$\fg}
ou simplement {\og ($F$-)complexes surholonomes sur $U$\fg}.
Elles correspondent à un analogue $p$-adique de
la catégorie des complexes à cohomologie bornée et holonome
sur une variété complexe $ U_\C$ notée
$D ^\mathrm{b} _\mathrm{hol} (\D _{U _\C})$.
\end{defi}

\begin{prop}
  \label{def-surhol-rele2}
  Soient $g$ : $U _2 \rightarrow U _1$ un morphisme de variétés.
  On suppose que, pour $i=1,2$, $U _i$ se plonge dans un $\V$-schéma formel propre et lisse $\X _i$.
  On note alors $Z _i$ l'adhérence schématique de $U _i$ dans $X_i$ et
  $T _i := Z _i \setminus U_i$.

On dispose alors des foncteurs canoniques
$g _{+},\,g _{!}$ :
$F\text{-}D ^\mathrm{b} _\mathrm{surhol} (\D _{U _2})\rightarrow F\text{-}D ^\mathrm{b} _\mathrm{surhol} (\D _{U _1})$
et
$g  ^!,\,g  ^+$ :
$F\text{-}D ^\mathrm{b} _\mathrm{surhol} (\D _{U _1 })\rightarrow F\text{-}D ^\mathrm{b} _\mathrm{surhol} (\D _{U _2})$.
Ces foncteurs ne dépendent canoniquement que de $g$ et
sont respectivement {\og l'image directe\fg},
{\og l'image directe extraordinaire\fg}, {\og l'image inverse extraordinaire\fg} et {\og l'image inverse\fg} par $g $.

Lorsqu'il existe un morphisme $f\,:\,\X _2 \rightarrow \X _1$ propre et lisse
prolongeant $g$ alors
ceux-ci correspondent à ceux définis dans \ref{defg+g!hautbas}, i.e.,
$g_+ = f _{+}$, $g _!= \DD _{T _1}\circ  f _{+} \circ \DD _{T _2}$,
$g^!= \R \underline{\Gamma} ^\dag _{Z _2}  \circ  (\hdag T _2)\circ   f ^!$,
$g^+=  \DD _{T _2}\circ \R \underline{\Gamma} ^\dag _{Z _2} \circ  (\hdag T _2)\circ f ^! \circ \DD _{T _1}$.

\end{prop}

\begin{proof}
Quitte à remplacer $\X _2$ par $\X _2 \times \X _1$, $T _2$ par
la réunion des images inverses de $T _1$ et $T _2$ sur $X _2 \times X _1$,
$Z _2$ par l'adhérence schématique de $U _2$ dans $X _2 \times X _1$,
grâce à \ref{theoindpdt2},
on peut supposer qu'il existe un morphisme propre et lisse
$f$ : $\X _2\rightarrow \X _1$ prolongeant $g$.

On vérifie à présent que les foncteurs
définis dans \ref{defg+g!hautbas} ne dépendent canoniquement pas du prolongement $f$.
En effet, soit $f'\,:\, \X'_2 \rightarrow \X '_1$ un second morphisme de $\V$-schémas formels propres
et lisses prolongent $g$. Pour $i=1,2$, on note de même
$Z '_i$ l'adhérence schématique de $U _i$ dans $X'_i$ et $T' _i := Z '_i \setminus U_i$.

Alors, quitte à remplacer
$\X ' _i $ par $\X ' _i \times \X  _i$,
$T '_i$ par la réunion des images inverses de $T _i$ et $T '_i$ sur $X _i ' \times X _i$,
$Z ' _i$ par l'adhérence schématique de $U _i$ dans $X ' _i \times X _i$,
on peut supposer qu'il existe
un morphisme $u_1\,:\, \X' _1 \rightarrow \X _1$ (resp. $u_2\,:\, \X' _2 \rightarrow \X _2$) propre et lisse
induisant l'identité sur $U _1$ (resp. $U_2$) tels que
$f \circ u_2 = u_1 \circ f'$ et
tels que les
foncteurs
$u_{i+}$ et
$(\hdag T '_i)\circ   \R \underline{\Gamma} ^\dag _{Z '_i} \circ u_i ^!$ induisent
   des équivalences quasi-inverses
$F\text{-}\mathfrak{C} ^+ _{\X _i , T _i , Z _i}\cong F\text{-}\mathfrak{C} ^+ _{\X '_i , T' _i , Z' _i}$
(voir \ref{theoindpdt2}).

Soit $ \E ' _2 \in F\text{-}\mathfrak{C} ^+ _{\X '_2 , T '_2 , Z' _2}$.
Il lui correspond canoniquement $\E _2 = u _{2+} (\E '_2) \in F\text{-}\mathfrak{C} ^+ _{\X _2 , T _2 , Z _2}$.
On obtient alors l'isomorphisme canonique
$f '_+ (\E ' _2 ) \riso u_{1+} f _+ (\E _2)$.
En d'autres termes, le foncteur $g_+$ est canoniquement indépendant du choix du prolongement $f$.
On procède de même pour les trois autres foncteurs.

\end{proof}

\begin{rema}
 Avec les notations de \ref{def-surhol-rele2} et de sa preuve,
 lorsque $g$ est propre,
on peut même supposer que le prolongement $f$ vérifie $T _2 = f _0 ^{-1} (T _1)$.

En effet,
d'après \cite[II.4.8.(e)]{HaAG}, le fait que $g$ soit propre implique que l'immersion ouverte
$U _2 \hookrightarrow Z _2 \setminus f _0 ^{-1} (T _1)$ est un morphisme propre.
Comme $U _2$ est dense dans $Z _2$, on obtient alors
$U _2 =Z _2 \setminus f _0 ^{-1} (T _1)$.
D'où :
$F\text{-}\mathfrak{C} ^+ _{\X _2 , T _2 , Z _2}
=F\text{-}\mathfrak{C} ^+ _{\X _2 , f _0 ^{-1} (T _1) , Z _2}$ (voir \ref{propindX}).

\end{rema}

\begin{prop}
Soit $U$ une variété
qui se plonge dans un $\V$-schéma formel propre et lisse $\X$.
Soient $u$ : $Y  \hookrightarrow U $ une immersion fermée
et $j$ : $U  \setminus Y  \subset U $ l'immersion ouverte induite.
On dispose alors, pour tout $\E \in F\text{-}D ^\mathrm{b} _\mathrm{surhol} (\D _{U  })$,
du triangle distingué de localisation :
\begin{equation}
  \label{defitrgllocsch}
  u _+ u ^! (\E) \rightarrow \E \rightarrow j_+ j ^! (\E)
\rightarrow
u _+ u ^! (\E) [1].
\end{equation}
\end{prop}

\begin{proof}
Notons $Z$ (resp. $Z ' $) l'adhérence schématique
de $U $ (resp. $Y$) dans $X $ et $T:= Z \setminus U$.
Considérons $\E $ comme un élément de
$F\text{-}\mathfrak{C} ^+ _{\X  , T  , Z }$.
On dispose du triangle de localisation (voir \ref{tri-local}) :
$$\R \underline{\Gamma} ^\dag _{Z '} (\E) \rightarrow \E \rightarrow
(\hdag Z ' ) (\E)
\rightarrow
\R \underline{\Gamma} ^\dag _{Z '} (\E) [1].$$

Comme $Y=Z ' \setminus T$, on remarque que
$F\text{-}\mathfrak{C} ^+ _{Y}=F\text{-}\mathfrak{C} ^+ _{\X  , T  , Z '}$.
Par définition (voir \ref{def-surhol-rele2}),
$u ^! (\E) = \R \underline{\Gamma} ^\dag _{Z '} (\E) $
et $u_+ (\R \underline{\Gamma} ^\dag _{Z '} (\E) ) = \R \underline{\Gamma} ^\dag _{Z '} (\E) $.
De même, comme $U  \setminus Y = Z \setminus (T \cup Z')$ alors
$F\text{-}\mathfrak{C} ^+ _{U\setminus Y}=F\text{-}\mathfrak{C} ^+ _{\X  , T \cup Z'  , Z }$,
$j ^! (\E)=(\hdag Z ' ) (\E)$ et
$j_+ ((\hdag Z ' ) (\E)) =(\hdag Z ' ) (\E)$.

\end{proof}

\begin{vide}
[Isomorphisme de dualité relative]
\label{def-surhol-rele3}
On garde les notations de \ref{def-surhol-rele2}.
Comme il est encore conjectural que l'isomorphisme de dualité relative d'un morphisme propre (voir \ref{iso-dualrel} puis \ref{bidua+rel})
soit compatible à Frobenius, il n'est pas clair que
le foncteur $\DD _{\X_1,T _1}$ :
$F\text{-}\mathfrak{C} ^+ _{\X _1 , T _1 , Z _1} \rightarrow
F\text{-}\mathfrak{C} ^+ _{\X _1 , T _1 , Z _1}$ soit indépendant, à isomorphisme canonique près,
du choix de $\X _1$, $T _1$ et $Z _1$ tels que $U _1=Z _1 \setminus T_1$.
La construction d'un foncteur dual
sur $F\text{-}D ^{\mathrm{b} } _\mathrm{surhol} (\D _{U _1 })$ devient donc problématique.

Pour y remédier, on construit alors
la catégorie $F\text{-}D ^{\mathrm{b} \vee } _\mathrm{surhol} (\D _{U _1})$
ainsi que le foncteur dual $\DD _{U _1}$ induisant l'équivalence de catégories
$\DD _{U _1}$ : $F\text{-}D ^{\mathrm{b} } _\mathrm{surhol} (\D _{U _1 })
\cong F\text{-}D ^{\mathrm{b} \vee } _\mathrm{surhol} (\D _{U _1})$
de la manière suivante.
\medskip

De façon analogue au théorème \ref{theoindpdt2}, lorsque $g$ est un isomorphisme, on vérifie que
les foncteurs $f_!$ et
 $\DD \circ (\hdag T _2)\circ   \R \underline{\Gamma} ^\dag _{Z _2} \circ f ^! \circ \DD$ induisent
   des équivalences quasi-inverses entre les catégories
 $(F\text{-})\mathfrak{C} ^+ _{\X _2, T _2, Z _2}$ et $(F\text{-})\mathfrak{C} ^+ _{\X _1, T _1, Z _1}$.
Il résulte de cette {\og version dualisée\fg} du théorème \ref{theoindpdt2} que
la catégorie
$(F\text{-})\mathfrak{C} ^+ _{\X _1, T _1, Z _1}$ est alors indépendante du choix de $\X _1$, $T _1$ et $Z _1 $
tels que $U _1=Z _1 \setminus T_1$ à isomorphisme canonique près
(ces isomorphismes canoniques sont dans cette deuxième construction de la forme $f_!$ ou
 $\DD \circ (\hdag T _2)\circ   \R \underline{\Gamma} ^\dag _{Z _2} \circ f ^! \circ \DD$).
On la note alors $(F\text{-})D ^{\mathrm{b}\vee} _\mathrm{surhol} (\D _{U })$.

Le foncteur $\DD _{\X_1,T _1}$ :
$F\text{-}\mathfrak{C} ^+ _{\X _1 , T _1 , Z _1} \rightarrow
F\text{-}\mathfrak{C} ^+ _{\X _1 , T _1 , Z _1}$ induit alors
les équivalences de catégories
$$\DD _{U _1}\, :\, F\text{-}D ^{\mathrm{b} } _\mathrm{surhol} (\D _{U _1 })
\cong F\text{-}D ^{\mathrm{b} \vee } _\mathrm{surhol} (\D _{U _1})
\text{ et }
\DD _{U _1}\, :\, F\text{-}D ^{\mathrm{b} \vee} _\mathrm{surhol} (\D _{U _1})
\cong F\text{-} D ^{\mathrm{b} } _\mathrm{surhol} (\D _{U _1}).$$
En effet, via le théorème de bidualité
qui commute à Frobenius (voir \ref{bidual}), on vérifie que le foncteur
$\DD _{U _1}$ ne dépend pas à isomorphisme canonique près du choix du triplet
$(\X _1 , T _1 , Z _1)$ tel que $U _1 =Z _1\setminus T_1$.
\medskip

$\bullet$
En oubliant Frobenius (voir \ref{remafrob}), on a
$D ^{\mathrm{b} \vee} _\mathrm{surhol} (\D _{U _1}) =D ^{\mathrm{b} } _\mathrm{surhol} (\D _{U _1})$
et
$\DD _{U _1}$ : $D ^{\mathrm{b} } _\mathrm{surhol} (\D _{U _1 })
\cong D ^{\mathrm{b} } _\mathrm{surhol} (\D _{U _1})$.
Lorsque $g$ est propre, il découle de \ref{bidua+rel} que
l'on dispose de l'isomorphisme de foncteurs
$D ^\mathrm{b} _\mathrm{surhol} (\D _{U _2 }) \rightarrow D ^{\mathrm{b}\vee} _\mathrm{surhol} (\D _{U _1 })$ :
\begin{equation}
  \label{iso-dualrel-gen}
  \DD _{U _1} \circ g _+ \riso g_+ \circ \DD _{U _2}.
\end{equation}
Cet isomorphisme correspond au {\og théorème de dualité relative\fg}.
\medskip

$\bullet$
Enfin, si le théorème de dualité relative de \ref{iso-dualrel} est compatible à Frobenius, alors
celui de \ref{iso-dualrel-gen} l'est aussi,
(l'isomorphisme $f_! \riso f_+$ commute à Frobenius et donc)
$F\text{-}D ^{\mathrm{b} \vee} _\mathrm{surhol} (\D _{U _1}) =F\text{-}D ^{\mathrm{b} } _\mathrm{surhol} (\D _{U_1 })$.

\end{vide}

\begin{vide}[Coefficient constant]
\label{coeff-const}
    Soient $\X$ un $\V$-schéma formel propre et lisse, $T$ et $Z$ deux sous-schémas fermés de $X$,
    $d _{Y/X}$ la dimension relative de $Z$ sur $X$,
et $U := Z \setminus T$.
Le {\og coefficient constant associé à $U$\fg} est défini en posant
$\O _U := \R \underline{\Gamma} ^\dag _{Z} (\hdag T) (\O _{\X,\,\Q})[-d _{U/X}]$.
Celui-ci ne dépend pas du choix du
triplet $(\X  ,\, T  ,\, Z )$ tel que $Z \setminus T = U$.
En effet, d'après \ref{proptheo}, ce complexe est un objet de
$F\text{-}D ^\mathrm{b} _\mathrm{surhol} (\D _{U })$.
Cela résulte alors de \ref{def-surhol-rele2} en remarquant que
$\O _U = j^! (\O _{\X,\Q}) [-d _{U/X}]$.
Notons que le foncteur $j^! [-d _{U/X}]$ est un analogue de $\L j^*$.

Lorsque $U$ est lisse et $T$ est un diviseur,
le complexe $\O _U= \R \underline{\Gamma} ^\dag _{Z} (\hdag T) (\O _{\X,\,\Q})[-d _{U/X}]$ est
réduit à un terme.
En effet, en notant $\Y:= \X \setminus T$, il suffit d'après \ref{4312Be1} de le voir en dehors de $T$.
Or,
$\R \underline{\Gamma} ^\dag _{Z} (\hdag T) (\O _{\X,\,\Q})[-d _{U/X}] |\Y
\riso
\R \underline{\Gamma} ^\dag _{U} (\O _{\Y,\,\Q})[-d _{U/X}]$.
Ce dernier est réduit à un terme car $U$ est lisse.
Dans ces conditions, $\O _U$ est un $F\text{-}\D _U$-module arithmétique surholonome.
Par exemple, si $Z =X$, $\O _U = \O _{\X} (\hdag T) _\Q$.
\end{vide}

\section{$\D$-modules arithmétiques surholonomes sur une variété}
\label{rem-sch-sep}
  Soit $U$ une variété. Lorsqu'il n'existe pas de $\V$-schéma formel propre et lisse $\X$,
  de sous-schéma fermé $Z$ de $X$ et de diviseur $T$ de $X$ tels que $U = Z \setminus T$, on peut
  néanmoins étendre la construction de
$(F\text{-})\mathfrak{M}  _{U }$, $(F\text{-})\mathfrak{M} ^+ _{U }$
par recollement.
Nous détaillons la procédure dans cette section.

\begin{nota}
\label{nota-rem-sch-sep}
On choisit un recouvrement fini ouvert $(U _\alpha) _{\alpha \in \Lambda}$ de $U$ tel que,
pour tout $\alpha\in \Lambda$, il existe
un $\V$-schéma formel propre et lisse $\X _\alpha$,
  un sous-schéma fermé $Z _\alpha$ de
  $X _\alpha$ et un diviseur $T_\alpha$ de $X_\alpha$ tel que $U _\alpha = Z _\alpha\setminus T_\alpha$.
  De tels recouvrements existent bien. En effet,
 si $U _\alpha$ est un ouvert affine de $U$, il existe une immersion fermée de la forme
 $U _\alpha \hookrightarrow \A ^n _k$. Il suffit alors de prendre $\X _\alpha$ égal au
  complété $p$-adique de l'espace projectif $\P ^n _\V$, $T _\alpha$ le diviseur de $\P ^n _k$
  complémentaire de $\A ^n _k$ et $Z _\alpha$ l'adhérence schématique de $U _\alpha $ dans $\P ^n _k$.

  Pour tous $\alpha, \beta$ et $\gamma \in \Lambda$,
  on note $p _1 ^{{\alpha \beta}}$ : $\X _\alpha \times _\S \X _\beta \rightarrow \X _\alpha$ et
  $p _2 ^{{\alpha \beta}}$ : $\X _\alpha \times _\S \X _\beta \rightarrow \X _\beta$
  les projections canoniques, $Z _{{\alpha \beta}}$
  l'adhérence schématique de
  $U _\alpha \cap U _\beta$ dans $X _\alpha \times X _\beta$ (via l'immersion
  $U _\alpha \cap U _\beta \hookrightarrow U _\alpha \times U _\beta
  \hookrightarrow X _\alpha \times X _\beta$),
  $T _{{\alpha \beta}}=
  \smash{p _1 ^{\alpha \beta }} ^{-1} (T _\alpha) \cup \smash{p _2 ^{\alpha \beta }} ^{-1} (T _\beta)$.
  On remarque que $(F\text{-})\mathfrak{M} ^+ _{U _\alpha \cap U _\beta }=
  (F\text{-})\mathfrak{M} ^+ _{\X _\alpha \times \X _\beta , T _{\alpha \beta}, Z _{\alpha \beta}}$.
  De même, en notant $Z _{\alpha \beta \gamma}$ l'adhérence schématique de
  $U _\alpha \cap U _\beta \cap U _\gamma$ dans $X _\alpha \times X _\beta \times X _\gamma$, et
  $T _{\alpha \beta \gamma}$ la réunion des images inverses
  de $T _\alpha$, $T _\beta$ et $T_\gamma$ par les projections de
  $\X _\alpha \times \X _\beta \times \X _\gamma$ sur $\X _\alpha$, $\X _\beta$
  et $\X _\gamma$, on remarque que
  $(F\text{-})\mathfrak{M} ^+ _{U _\alpha \cap U _\beta \cap U _\gamma}=
  (F\text{-})\mathfrak{M} ^+ _{\X _\alpha \times \X _\beta \times \X _\gamma ,
  T _{\alpha \beta \gamma}, Z _{\alpha \beta \gamma}}$.

  Soient
  $j _1 ^{\alpha \beta} $ : $U _{\alpha}  \cap U _{\beta} \subset U _\alpha$,
  $j _2 ^{\alpha \beta} $ : $U _{\alpha}  \cap U _{\beta} \subset U _\beta$,
  $j _{12} ^{\alpha \beta \gamma} $ : $U _{\alpha}  \cap U _{\beta}\cap U _{\gamma} \subset
  U _{\alpha}  \cap U _{\beta}$,
  $j _{23} ^{\alpha \beta \gamma} $ :
  $U _{\alpha}  \cap U _{\beta}\cap U _{\gamma}\subset U _\beta \cap U _\gamma$
  et
  $j _{13} ^{\alpha \beta \gamma} $ : $U _{\alpha } \cap U _{\beta}  \cap U _{\gamma}
  \subset
  U _\alpha \cap U _\gamma$
  les immersions ouvertes.
\end{nota}

\begin{defi}
\label{defi-rem-sch-sep}
Avec les notations \ref{nota-rem-sch-sep},
  on définit la catégorie $(F\text{-})\mathfrak{M} ^+ (U,\, (U _\alpha, \X _\alpha, T_\alpha,  Z_\alpha) _{\alpha \in \Lambda})$
  de la façon suivante :
\begin{itemize}
  \item Un objet est constitué par la donnée, pour tout $\alpha \in \Lambda$,
  d'un objet $\E _\alpha$ de
  $(F\text{-})\mathfrak{M} ^+ _{\X _\alpha, T _\alpha, Z _\alpha }$ et,
  pour tous $\alpha,\beta \in \Lambda$,
  d'un isomorphisme
  $\theta _{\alpha \beta}$ :
  $j _2 ^{\alpha \beta !} (\E _\beta)
  \riso j _1 ^{\alpha \beta !} (\E _\alpha)$
  de $(F\text{-})\mathfrak{M} ^+ _{U _\alpha \cap U _\beta }$
  (i.e., d'un isomorphisme
$\D ^\dag _{\X _\alpha \times _\S \X _\beta, \, \Q}$-linéaire
$\R \underline{\Gamma} ^\dag _{Z _{{\alpha \beta}}}
\circ (\hdag T _{{\alpha \beta}})\circ p _2 ^{\alpha \beta !}  (\E _\beta)\riso
\R \underline{\Gamma} ^\dag _{Z _{{\alpha \beta}}}
\circ (\hdag T _{{\alpha \beta}})\circ p _1 ^{\alpha \beta !}  (\E _\alpha)$),
ces isomorphismes vérifiant dans $(F\text{-})\mathfrak{M} ^+ _{U _\alpha \cap U _\beta \cap U _\gamma}$
la condition de cocycle
$j _{13}
^{\alpha \beta \gamma !} (\theta _{  \alpha \gamma} )=
j _{12} ^{\alpha \beta \gamma !} (\theta _{  \alpha \beta} )
\circ j _{23} ^{\alpha \beta \gamma !} ( \theta _{ \beta \gamma })$.
La famille d'isomorphismes $(\theta _{\alpha \beta}) _{\alpha ,\beta \in \Lambda}$
est appelée {\og donnée de recollement\fg} de $ (\E _\alpha) _{\alpha \in \Lambda} $.

\item Une flèche
$ ((\E _\alpha) _{\alpha \in \Lambda}, (\theta _{  \alpha \beta}) _{\alpha,\beta  \in \Lambda}) \rightarrow
  ((\E '_\alpha) _{\alpha \in \Lambda}, (\theta '_{  \alpha \beta}) _{\alpha,\beta  \in \Lambda})$
de la catégorie
$(F\text{-})\mathfrak{M} ^+
  (U,\, (U _\alpha, \X _\alpha, T_\alpha,  Z_\alpha) _{\alpha \in \Lambda})$
  est une famille de morphismes
  $ \E _\alpha \rightarrow   \E '_\alpha$
  commutant aux données de recollement respectives.

\end{itemize}

On remarque grâce à \ref{j!} que les foncteurs $j _1 ^{\alpha \beta !}$, $j _{12} ^{\alpha \beta \gamma !}$ etc.
sont exacts, ce qui donne un sens à notre définition.
\end{defi}

\begin{prop}
\label{prop-rem-sch-sep}
La catégorie
$(F\text{-})\mathfrak{M} ^+
  (U,\, (U _\alpha, \X _\alpha, T_\alpha,  Z_\alpha) _{\alpha \in \Lambda})$
  ne dépend pas du choix de
$(U _\alpha, \X _\alpha, T_\alpha,  Z_\alpha) _{\alpha \in \Lambda}$
vérifiant les conditions de \ref{nota-rem-sch-sep}.
\end{prop}

\begin{proof}
La preuve se fait en trois étapes suivantes :

i) S'il existe un $\V$-schéma formel propre et lisse $\X$,
un diviseur $T$ de $X$ tel que $U = Z \setminus T$ où $Z$ est l'adhérence de $U$ dans $X$
alors,
pour tout recouvrement fini ouvert $(\X _\alpha) _{\alpha \in \Lambda}$ de $\X$ tel que,
quitte à prendre un recouvrement plus fin,
$T _\alpha := X \setminus [(X\setminus T) \cap X _\alpha]$ soit un diviseur de
$X$ (par exemple, puisque $X$ est propre et lisse,
si $(X\setminus T) \cap X _\alpha$ est affine : \cite[1.3.1]{caro-2006-surcoh-surcv}),
on dispose, en notant
$U _{\alpha }:=X _\alpha \cap U$,
 d'une équivalence entre
les catégories $(F\text{-})\mathfrak{M} ^+   _{\X, \, T , \, Z}$
et $(F\text{-})\mathfrak{M} ^+   (U,\, (U _\alpha, \X , T_\alpha,  Z) _{\alpha \in \Lambda})$.

  En effet, soient
  $p_1$ et $p_2$ : $\X \times _\S\X \rightarrow \X $
  les projections respectives à gauche et droite canoniques,
$T _{{\alpha \beta}}=
  p _1  ^{-1} (T _\alpha) \cup p _2  ^{-1} (T _\beta)$,
  $Z'$ l'adhérence de $U$ dans $X \times _k X $,
  $(\E _\alpha) _{\alpha \in \Lambda}$ une famille
  d'objets de $(F\text{-})\mathfrak{M} ^+ _{\X , T _\alpha, Z  }$,
  munie d'isomorphismes
  $\D ^\dag _{\X \times _\S \X , \, \Q}$-linéaires,
$ \theta _{\alpha \beta}$ :
$\R \underline{\Gamma} ^\dag _{Z' }
\circ (\hdag T _{{\alpha \beta}})\circ p _2 ^{ !}  (\E _\beta)\riso
\R \underline{\Gamma} ^\dag _{Z'}
\circ (\hdag T _{{\alpha \beta}})\circ p _1 ^{ !}  (\E _\alpha)$,
ceux-ci vérifiant
la condition de cocycle. En notant $\delta $ : $\X \hookrightarrow \X \times _\S \X$ l'immersion
fermée diagonale,
comme $p _1 \circ \delta = p _2 \circ \delta$, il résulte de
\ref{fct-cohloc}, \ref{commutfonctcohlocal2}, \ref{commutfonctcohlocal2+}
que $\delta ^! (\theta _{\alpha \beta})$
est canoniquement isomorphe à un isomorphisme de la forme
$\eta _{\alpha \beta}\,:\, (\hdag T _\alpha) (\E _\beta) \riso (\hdag T _\beta) (\E _\alpha)$.
La famille $\eta _{\alpha \beta}$ vérifie en outre la condition de cocycle
$(\hdag T _\beta ) (\eta _{\alpha \gamma}) =
(\hdag T _\gamma) (\eta _{\alpha \beta}) \circ (\hdag T _\alpha) (\eta _{\beta \gamma})$.
On conclut alors l'étape i) via le lemme suivant.
\begin{lemm}
Avec les notations ci-dessus,
  soit $(\E _\alpha) _{\alpha \in \Lambda}$ une famille de $\D ^\dag _{\X} (\hdag T _\alpha) _{\Q}$-modules
  surholonomes
  munie d'isomorphismes
$\eta _{\alpha \beta}$ : $(\hdag T _\alpha) (\E _\beta) \riso (\hdag T _\beta) (\E _\alpha)$
avec $\alpha ,\beta \in \Lambda$ et
vérifiant, pour tous $\alpha ,\beta ,\gamma \in \Lambda$, l'égalité
$(\hdag T _\beta ) (\eta _{\alpha \gamma}) =
(\hdag T _\gamma) (\eta _{\alpha \beta}) \circ (\hdag T _\alpha) (\eta _{\beta \gamma})$.
Il existe alors un unique $\D ^\dag _{\X} (\hdag T )_{\Q}$-module $\E$ (à isomorphisme canonique près)
qui soit surholonome et
muni d'une famille d'isomorphismes $\eta _\alpha$ : $ (\hdag T _\alpha) (\E) \riso \E _\alpha$ telle que
$(\hdag T _\beta) (\eta  _\alpha) \circ (\hdag T _\alpha) (\eta ^{-1}  _\beta)=\eta _{\alpha \beta}$.
\end{lemm}
\begin{proof}
  On pose $\E =\ker ( \xymatrix { { \oplus _{\alpha \in \Lambda} \E _\alpha} \ar@<1ex>[r]\ar@<-1ex>[r]
  &
  {\oplus _{\alpha, \beta \in \Lambda} (\hdag T _\beta) (\E _\alpha)}})$, où la première flèche est
  induite par les morphismes canoniques $\E _\alpha \rightarrow (\hdag T _\beta) (\E _\alpha)$
  tandis que la deuxième se construit terme à terme via les flèches de la forme
  $\E _\beta \rightarrow (\hdag T _\alpha) (\E _\beta) \overset{\eta _{\alpha \beta}}{\longrightarrow}
  (\hdag T _\beta) (\E _\alpha)$. Comme, pour tout $\alpha$,
  $T \subset T _\alpha$, $\E$ est muni d'une structure canonique de
$\D ^\dag _{\X} (\hdag T )_{\Q}$-module. Par \ref{propkerimsurcoh}, $\E$ est en outre
surholonome (en tant que $\D ^\dag _{\X,\Q}$-module).
Comme le morphisme canonique $\E _\alpha \rightarrow (\hdag T _\alpha) (\E _\alpha)$
est un isomorphisme (car, par \ref{4312Be1}, c'est le cas en dehors de $T _\alpha$),
la projection $\E\rightarrow \E _\alpha $ induit le morphisme canonique
$ (\hdag T _\alpha) (\E) \rightarrow \E _\alpha$ que l'on note $\eta _\alpha$.
Vérifions à présent que $\eta _\alpha$ est un isomorphisme
et que
la formule $(\hdag T _\beta) (\eta  _\alpha) \circ (\hdag T _\alpha) (\eta ^{-1}  _\beta)=\eta _{\alpha \beta}$
est correcte.
Comme cela est local
  en $\X$, supposons donc $\X$ affine. Pour tous $\alpha, \beta$,
  pour toute section $e _\alpha$ de $\E _\alpha$,
  notons $(\hdag T _\beta) (e _\alpha)$ l'image de $e _\alpha$ dans $(\hdag T _\beta) (\E _\alpha)$
  via le morphisme canonique $\E _\alpha \rightarrow (\hdag T _\beta) (e _\alpha)$.
  L'ensemble des sections globales de $\E$
  est $\{ (e _\alpha) _{\alpha}, \, \text{avec } e _\alpha\in \E _\alpha\,
  \text{tels que }\forall \beta,\alpha \in \Lambda,\,
  \eta _{\alpha \beta} ( (\hdag T _\alpha) (e_\beta)) = (\hdag T _\beta) (e _\alpha)\}$.
  On calcule alors que la projection $(\hdag T_\alpha) (\E  )\rightarrow (\hdag T_\alpha) (\E_\alpha)$
  est un isomorphisme (en effet, $e _\alpha$ détermine alors $(\hdag T _\alpha) (e _\beta )$).
Ainsi, $\eta _\alpha$ est un isomorphisme. On calcule de même que
la formule $(\hdag T _\beta) (\eta  _\alpha) \circ (\hdag T _\alpha) (\eta ^{-1}  _\beta)=\eta _{\alpha \beta}$
est correcte.

Traitons à présent l'unicité. Soit
$\E '$ un
$\D ^\dag _{\X,\,\Q} (\hdag T )$-module surholonome
muni d'une famille d'isomorphismes $\eta ' _\alpha$ : $ (\hdag T _\alpha) (\E') \riso \E _\alpha$ telle que
$(\hdag T _\beta) (\eta ' _\alpha) \circ (\hdag T _\alpha) (\eta ^{\prime -1}  _\beta)=\eta _{\alpha \beta}$.
Il en dérive les morphismes
$\E ' \rightarrow (\hdag T _\alpha) (\E') \underset{\eta ' _\alpha}{\riso} \E _\alpha$ qui induisent le morphisme
$\E' \rightarrow
\ker ( \xymatrix { { \oplus _{\alpha \in \Lambda} \E _\alpha} \ar@<1ex>[r]\ar@<-1ex>[r]
  &
  {\oplus _{\alpha, \beta \in \Lambda} (\hdag T _\beta) (\E _\alpha)}})=\E$.
Comme, pour tout $\alpha \in \Lambda$, $(\hdag T _\alpha)( \E') \riso (\hdag T _\alpha)( \E)$,
on obtient $\E' \riso \E$.

\end{proof}
ii) Soit
un recouvrement fini ouvert $(U _\alpha) _{\alpha \in \Lambda}$ (resp. $(U '_{\alpha '}) _{\alpha '\in \Lambda '}$)
de $U$ tel qu'il existe
un $\V$-schéma formel propre et lisse $\X _\alpha$ (resp. $\X '_{\alpha'}$),
  un sous-schéma fermé $Z _\alpha$ de $X _\alpha$ (resp $Z '_{\alpha'}$ de $X '_{\alpha'}$)
  et un diviseur $T_\alpha$ de $X_\alpha$ (resp. $T '_{\alpha'}$ de $X '_{\alpha'}$)
  tel que
  $U _\alpha = Z _\alpha\setminus T_\alpha$
  (resp. $U '_{\alpha'}= Z '_{\alpha'} \setminus T '_{\alpha'}$).
  Grâce à i), quitte à prendre un recouvrement plus fin que
  $(U _\alpha) _{\alpha \in \Lambda}$ et $(U '_{\alpha '}) _{\alpha '\in \Lambda '}$,
  on se ramène à traiter le cas où $\Lambda =\Lambda '$ et $U _\alpha = U '_{\alpha}$.

iii) Quitte à remplacer $\X ' _\alpha$ par $\X' _\alpha \times _\S \X _\alpha$,
on peut supposer qu'il existe un morphisme propre et lisse
$\X ' _\alpha \rightarrow \X _\alpha$ dont la restriction à $U _\alpha$ induise l'identité de $U _\alpha$.
On obtient ainsi des foncteurs canoniques (\ref{defg+g!hautbas})
$Id _\alpha ^!$ : $(F\text{-})\mathfrak{M} ^+   _{\X _\alpha, \,  Z _\alpha, \, T_\alpha}
\rightarrow (F\text{-})\mathfrak{M} ^+   _{\X '_\alpha ,\,  Z' _\alpha ,\, T '_\alpha}$
et
$Id _{\alpha +}$ :
$(F\text{-})\mathfrak{M} ^+   _{\X '_\alpha ,\,  Z '_\alpha, \, T'_\alpha}
\rightarrow (F\text{-})\mathfrak{M} ^+   _{\X _\alpha ,\,  Z _\alpha ,\, T _\alpha}$.
Dans ce cas, si
$ (\E _\alpha) _{\alpha \in \Lambda}$ est une famille
d'objets $\E _\alpha$ de $(F\text{-})\mathfrak{M} ^+ _{\X _\alpha, T _\alpha, Z _\alpha }$,
munie d'une donnée de recollement,
on lui associe la famille $ (Id ^! _\alpha (\E _\alpha) ) _{\alpha \in \Lambda}$, qui est munie d'une structure
canonique induite de donnée de recollement (grâce aux isomorphismes de transitivité de l'image inverse
extraordinaire de \ref{defg+g!hautbas} et à son associativité).
De façon analogue, on construit un foncteur quasi-inverse en remplaçant $Id _\alpha ^!$ par $Id _{\alpha +}$.

\end{proof}

\begin{defi}
  \label{defi2-rem-sch-sep}
Avec les notations \ref{nota-rem-sch-sep},
on désigne par $(F\text{-})\mathfrak{M} ^+ _{U }$, la catégorie
$(F\text{-})\mathfrak{M} ^+  (U,\, (U _\alpha, \X _\alpha, T_\alpha,  Z_\alpha) _{\alpha \in \Lambda})$
qui ne dépend canoniquement que de $U$ (voir \ref{prop-rem-sch-sep}).
Un objet de $(F\text{-})\mathfrak{M} ^+ _{U }$ sera nommé
{\og $(F\text{-})\D_U$-module arithmétique surholonome\fg}.
En remplaçant {\og $\mathfrak{M}^+$\fg} par {\og $\mathfrak{M}$\fg},
on construit de même la catégorie
$(F\text{-})\mathfrak{M}  _{U }$
des {\og $(F\text{-})\D_U$-modules arithmétiques surcohérents\fg}.
\end{defi}

\begin{defi}
\label{defi-const-rem-sch-sep}
Avec les notations \ref{nota-rem-sch-sep},
  lorsque $U$ est lisse, la famille $(\R \underline{\Gamma} ^\dag _{Z _\alpha}
(\hdag T _\alpha) (\O _{\X _\alpha ,\,\Q})[-d _{Z _\alpha /X _\alpha }]) _{\alpha \in \Lambda}$
est canoniquement munie d'une donnée de recollement.
Cet objet de
$F\text{-}\mathfrak{M} ^+ _{U }$ sera noté de $\O _U$ (on vérifie que cela ne dépend pas du choix
de la famille de quadruplet $(U_\alpha ,\, \X _\alpha, \,  Z _\alpha, \, T_\alpha)_{\alpha \in \Lambda}$)
et sera appelé {\it coefficient constant associé à $U$}.

\end{defi}

On termine cette section par la description
des catégories $(F\text{-})D ^\mathrm{b} (\mathfrak{M}  _{U })$ et
$(F\text{-})D ^\mathrm{b}  (\mathfrak{M} ^+  _{U })$.

\begin{vide}
\label{Dsurcoh}
On conserve les notations de \ref{nota-rem-sch-sep}.
On désigne par $C ^\mathrm{b} (\mathfrak{M}  _{\X _\alpha , T _{\alpha }, Z _{\alpha }})$
la sous-catégorie pleine de $C ^\mathrm{b} (\D ^\dag  _{\X_\alpha  }(\hdag T _\alpha) _\Q)$
des complexes
$\E  ^\bullet $ tel que $\E ^n  \in
\mathfrak{M}  _{\X _\alpha , T _{\alpha }, Z _{\alpha }}$.
D'après \ref{j!}, on dispose du foncteur exact $j _1 ^{\alpha \beta !}
\,:\,\mathfrak{M}  _{\X _\alpha , T _{\alpha }, Z _{\alpha }} \rightarrow
\mathfrak{M}  _{\X _\alpha \times \X _\beta , T _{\alpha \beta}, Z _{\alpha \beta}}$.
Il induit donc un foncteur
$j _1 ^{\alpha \beta !}
\,:\,
C ^\mathrm{b} ( \mathfrak{M}  _{\X _\alpha , T _{\alpha }, Z _{\alpha }} )
\rightarrow
C ^\mathrm{b} (\mathfrak{M}  _{\X _\alpha \times \X _\beta , T _{\alpha \beta}, Z _{\alpha \beta}})$.
On obtient aussi un foncteur
$j _2 ^{\alpha \beta !}
\,:\,
C ^\mathrm{b} ( \mathfrak{M}  _{\X _\alpha , T _{\alpha }, Z _{\alpha }} )
\rightarrow
C ^\mathrm{b} (\mathfrak{M}  _{\X _\alpha \times \X _\beta , T _{\alpha \beta}, Z _{\alpha \beta}})$
et
de même pour $j _{12}
^{\alpha \beta \gamma !} (\theta _{  \alpha \gamma} )$,
$j _{13}
^{\alpha \beta \gamma !} (\theta _{  \alpha \gamma} )$,
$j _{23}
^{\alpha \beta \gamma !} (\theta _{  \alpha \gamma} )$.
On dispose de foncteurs analogues en remplaçant {\og $\mathfrak{M}$\fg} par {\og $\mathfrak{M} ^+$\fg}.
\end{vide}

\begin{defi}
\label{defi-D(M)}
$\bullet$ Avec les notations de \ref{Dsurcoh},
on définit la catégorie $C ^\mathrm{b} (\mathfrak{M}  (U,\, (U _\alpha, \X _\alpha, T_\alpha,  Z_\alpha) _{\alpha \in \Lambda}))$
de la manière suivante :

\begin{itemize}
  \item Un objet est une famille
$\E = ( \E _\alpha ^\bullet , \theta _{\alpha \beta}) _{\alpha,\beta  \in \Lambda }$
où $\E _\alpha ^\bullet \in
C ^\mathrm{b} ( \mathfrak{M}  _{\X _\alpha , T _{\alpha }, Z _{\alpha }} )$
$\theta _{\alpha \beta}  $ :
  $j _2 ^{\alpha \beta !} (\E ^\bullet   _\beta)
  \riso j _1 ^{\alpha \beta !} (\E ^\bullet  _\alpha)$
  est un isomorphisme de
  $C ^\mathrm{b} (\mathfrak{M}  _{\X _\alpha \times \X _\beta , T _{\alpha \beta}, Z _{\alpha \beta}})$
satisfaisant la condition de cocycle
$j _{13}
^{\alpha \beta \gamma !} (\theta _{  \alpha \gamma} )=j _{12}
^{\alpha \beta \gamma !} (\theta _{  \alpha \beta} ) \circ j _{23}
^{\alpha \beta \gamma !} ( \theta _{ \beta \gamma })$.
La famille d'isomorphismes $\theta _{\alpha \beta}  $ est dite {\it données de recollement}.

\item Un morphisme $f$ : $( \E ^\bullet _\alpha , \theta _{\alpha \beta}) _{\alpha,\beta  \in \Lambda}
\rightarrow ( \FF ^\bullet _\alpha , \theta ' _{\alpha \beta}) _{\alpha,\beta  \in \Lambda }$
est une famille de morphismes $f  _\alpha$ : $\E ^\bullet _\alpha  \rightarrow \FF ^\bullet _\alpha $
de $C ^\mathrm{b} ( \mathfrak{M}  _{\X _\alpha , T _{\alpha }, Z _{\alpha }} )$
commutant aux données de recollement respectives.

De manière analogue à \ref{prop-rem-sch-sep}, on vérifie que
$C ^\mathrm{b} (\mathfrak{M}  (U,\, (U _\alpha, \X _\alpha, T_\alpha,  Z_\alpha) _{\alpha \in \Lambda}))$
ne dépend pas, à équivalence canonique de catégories près, du choix de
$(U _\alpha, \X _\alpha, T_\alpha,  Z_\alpha) _{\alpha \in \Lambda}$
vérifiant les conditions de \ref{nota-rem-sch-sep}
et est d'ailleurs canoniquement équivalente à
$C ^\mathrm{b} (\mathfrak{M}  _{U} )$, où $\mathfrak{M}  _{U} $ est la catégorie abélienne de \ref{defi2-rem-sch-sep}..
On pourra donc les identifier.
\end{itemize}

$\bullet$  Un morphisme $f$ : $( \E ^\bullet _\alpha , \theta _{\alpha \beta}) _{\alpha,\beta  \in \Lambda}
\rightarrow ( \FF ^\bullet _\alpha , \theta '_{\alpha \beta}) _{\alpha,\beta  \in \Lambda }$
est un {\it quasi-isomorphisme} si, pour tout $\alpha \in \Lambda$,
le morphisme $f  _\alpha $ : $ \E ^\bullet _\alpha \rightarrow \FF ^\bullet _\alpha $
est un quasi-isomorphisme.

$\bullet$ Si $f $ et $g$ sont deux morphismes :
$( \E ^\bullet _\alpha , \theta _{\alpha \beta}) _{\alpha,\beta  \in \Lambda}
\rightarrow ( \FF ^\bullet _\alpha , \theta '_{\alpha \beta}) _{\alpha,\beta  \in \Lambda }$,
on dit que $f$ et $g$ sont {\it homotopes} s'il existe une famille de morphismes
$k  _\alpha $ : $\E ^\bullet _\alpha  \rightarrow \FF ^\bullet _\alpha [-1]$
de $C ^\mathrm{b} ( \mathfrak{M}  _{\X _\alpha , T _{\alpha }, Z _{\alpha }} )$
commutant aux données de recollement respectives et
telle que $f  _\alpha -g  _\alpha =d  k  _\alpha +k  _\alpha d $,
où $d$ sont les différentielles des complexes respectifs.
L'homotopie est une relation d'équivalence et se préserve par composition.
Lorsque deux morphismes sont homotopes et si l'un est un quasi-isomorphisme alors l'autre l'est également.

On définit la catégorie $K ^\mathrm{b} (\mathfrak{M}  (U,\, (U _\alpha, \X _\alpha, T_\alpha,  Z_\alpha) _{\alpha \in \Lambda}))$
dont les objets sont ceux de $C ^\mathrm{b} (\mathfrak{M}  (U,\, (U _\alpha, \X _\alpha, T_\alpha,  Z_\alpha) _{\alpha \in \Lambda}))$
et dont les morphismes sont les classes d'homotopie de ceux de
$C ^\mathrm{b} (\mathfrak{M}  (U,\, (U _\alpha, \X _\alpha, T_\alpha,  Z_\alpha) _{\alpha \in \Lambda}))$.
On vérifie aussi que
$K ^\mathrm{b} (\mathfrak{M}  (U,\, (U _\alpha, \X _\alpha, T_\alpha,  Z_\alpha) _{\alpha \in \Lambda}))$
est canoniquement équivalente à
$K ^\mathrm{b} (\mathfrak{M}  _{U})$.

$\bullet$  Si $f=(f_\alpha) _{\alpha \in \Lambda}$ : $( \E ^\bullet _\alpha , \theta _{\alpha \beta}) _{\alpha,\beta  \in \Lambda}
\rightarrow ( \FF ^\bullet _\alpha , \theta '_{\alpha \beta}) _{\alpha,\beta  \in \Lambda }$
est
une flèche
de $K ^\mathrm{b} (\mathfrak{M}  (U,\, (U _\alpha, \X _\alpha, T_\alpha,  Z_\alpha) _{\alpha \in \Lambda}))$,
on définit son cône $\mathcal{C} _f$ de la manière suivante :
en notant
$\mathcal{C} _{f_\alpha}$ le cône de $f _\alpha$
la famille
$\mathcal{C} _f:=(\mathcal{C} _{f_\alpha}) _{\alpha \in \Lambda}$
est canoniquement un objet de $K ^\mathrm{b} (\mathfrak{M}  (U,\, (U _\alpha, \X _\alpha, T_\alpha,  Z_\alpha) _{\alpha \in \Lambda}))$.
Un triangle distingué de $K ^\mathrm{b} (\mathfrak{M}  (U,\, (U _\alpha, \X _\alpha, T_\alpha,  Z_\alpha) _{\alpha \in \Lambda}))$
est par définition un triangle
isomorphe à un triangle de la forme $\E \overset{f}{\longrightarrow} \FF \rightarrow \mathcal{C} _f \rightarrow \E [1]$.
La catégorie $K ^\mathrm{b} (\mathfrak{M}  (U,\, (U _\alpha, \X _\alpha, T_\alpha,  Z_\alpha) _{\alpha \in \Lambda}))$ devient triangulée.
En localisant cette catégorie par rapport aux quasi-isomorphismes, on obtient une catégorie triangulée
que l'on notera $D ^\mathrm{b} (\mathfrak{M}  (U,\, (U _\alpha, \X _\alpha, T_\alpha,  Z_\alpha) _{\alpha \in \Lambda}))$.
De même,
$D ^\mathrm{b} (\mathfrak{M}  (U,\, (U _\alpha, \X _\alpha, T_\alpha,  Z_\alpha) _{\alpha \in \Lambda}))$
est canoniquement équivalente à
$D ^\mathrm{b} (\mathfrak{M}  _{U})$.

$\bullet$ On définit un foncteur
$F ^*$ : $D ^\mathrm{b} (\mathfrak{M}  (U,\, (U _\alpha, \X _\alpha, T_\alpha,  Z_\alpha) _{\alpha \in \Lambda}))
\rightarrow
D ^\mathrm{b} (\mathfrak{M}  (U,\, (U _\alpha, \X _\alpha, T_\alpha,  Z_\alpha) _{\alpha \in \Lambda}))$
en posant, pour tout
$( \E ^\bullet _\alpha , \theta _{\alpha \beta}) _{\alpha,\beta  \in \Lambda}\in
D ^\mathrm{b} (\mathfrak{M}  (U,\, (U _\alpha, \X _\alpha, T_\alpha,  Z_\alpha) _{\alpha \in \Lambda}))$,
$F ^* (( \E ^\bullet _\alpha , \theta _{\alpha \beta}) _{\alpha,\beta  \in \Lambda})
:=
(F ^*  \E ^\bullet _\alpha , \theta ' _{\alpha \beta}) _{\alpha,\beta  \in \Lambda}$, où
$\theta ' _{\alpha \beta}$ est construit par commutativité du diagramme :
$$\xymatrix {
{F ^* j _2 ^{\alpha \beta !} (\E ^\bullet   _\beta)}
\ar[r] _-\sim ^-{F ^* \theta _{\alpha \beta}}
\ar[d] _-\sim
&
{F ^*  j _1 ^{\alpha \beta !} (\E ^\bullet  _\alpha)}
\ar[d] _-\sim
\\
{j _2 ^{\alpha \beta !} (F ^*  \E ^\bullet   _\beta)}
\ar@{.>}[r]^-{\theta '_{\alpha \beta}}
&
{j _1 ^{\alpha \beta !} (F ^*  \E ^\bullet  _\alpha).}
}$$
On définit ensuite la catégorie $F\text{-}D ^\mathrm{b} (\mathfrak{M}  (U,\, (U _\alpha, \X _\alpha, T_\alpha,  Z_\alpha) _{\alpha \in \Lambda}))$
dont les objets sont les couples $( \E ,\phi)$ avec
$\E \in D ^\mathrm{b} (\mathfrak{M}  (U,\, (U _\alpha, \X _\alpha, T_\alpha,  Z_\alpha) _{\alpha \in \Lambda}))$
et $\phi$ est un isomorphisme $F ^* (\E) \riso \E$ dans
$D ^\mathrm{b} (\mathfrak{M}  (U,\, (U _\alpha, \X _\alpha, T_\alpha,  Z_\alpha) _{\alpha \in \Lambda})) $.

La catégorie $F\text{-}D ^\mathrm{b} (\mathfrak{M}  (U,\, (U _\alpha, \X _\alpha, T_\alpha,  Z_\alpha) _{\alpha \in \Lambda}))$
est canoniquement équivalente à
$F\text{-} D ^\mathrm{b} (\mathfrak{M}  _{U})$.

$\bullet$ En remplaçant {\og $\mathfrak{M}$\fg} par {\og $\mathfrak{M} ^+$\fg},
on définit de même la catégorie
$(F\text{-})D ^\mathrm{b} (\mathfrak{M} ^+  (U,\, (U _\alpha, \X _\alpha, T_\alpha,  Z_\alpha) _{\alpha \in \Lambda}))$.
Cette catégorie est canoniquement équivalente à
$(F\text{-}) D ^\mathrm{b} (\mathfrak{M} ^+ _{U})$.
\end{defi}

\begin{rema}\label{M-Csurcoh}
Lorsqu'il existe un $\V$-schéma formel propre et lisse $\X$,
  un sous-schéma fermé $Z$ de $X$ et un diviseur $T$ de $X$ tels que $U = Z \setminus T$,
  de manière analogue au cas des $\D$-modules sur les variétés complexes (voir par exemple \cite[VI.1.13, 2.10, 2.11]{borel}),
  il est raisonnable de conjecturer que
les foncteurs canoniques
$D ^\mathrm{b} (\mathfrak{M}  _{U}) \rightarrow D ^\mathrm{b} _\mathrm{surcoh} (\D  _{U})$,
$D ^\mathrm{b} (\mathfrak{M} ^+ _{U}) \rightarrow D ^\mathrm{b} _\mathrm{surhol} (\D  _{U})$
sont des équivalences de catégories.
On s'intéressera à cette question dans un prochain travail.
\end{rema}

\section{Application aux fonctions $L$}

On suppose dans cette section que $k$ est un corps fini à $p ^s$ éléments et $F$ désigne toujours
la puissance $s$-ième de l'endomorphisme de Frobenius.
Soit $U$ une variété qui se plonge dans un $\V$-schéma formel $\X$ propre et lisse.
Pour tout point fermé $y$ de $U$, on note $i _y$ : $\Spec k(y) \hookrightarrow U$ l'immersion fermée induite,
$p$ : $U \rightarrow \Spec k$ et $p _y$ : $\Spec k(y) \rightarrow \Spec k$ les morphismes structuraux.

\begin{defi}\label{defifonctLP}
Soient $\E \in F\text{-}D ^{\mathrm{b} \vee} _\mathrm{surhol} (\D _{U })$ (voir \ref{def-surhol-rele3})
et $S$ est un sous-ensemble de points de $U$.
On note $S ^\mathrm{0}$ le sous-ensemble des points fermés de $S$.
On définit respectivement la
fonction $L$ associée à $\E$ au dessus de $S$ et
la fonction cohomologique $P$ associée à $\E$ en posant :
\begin{align}
\notag
  L( S, \E, t )&:=\prod _{ y \in S^{\mathrm{0}} } \prod_{j \in \Z }
 \det _K \left( 1 - t F  {| \mathcal{H}^j ( p _{y+} i_{y}^+ (\E))} \right)^{(-1)^{j+1}},\\
 \notag
P(U , \E, t ) &:= \prod _{j\in\Z} \det_K \left( 1- t F {|{\mathcal{H}^j   (  p _{!} \E ) }}\right)^{(-1)^{j+1}},
\end{align}
où $F  {| \mathcal{H}^j ( p _{y+} i_{y}^+ (\E))}$ désigne l'action de Frobenius sur
$\mathcal{H}^j ( p _{y+} i_{y}^+ (\E))$
et $F{|{\mathcal{H}^j   (  p _{!} \E ) }}$ celle sur
$\mathcal{H}^j   (  p _{!} \E ) $
(celles-ci résultant des théorèmes de commutations de Frobenius
aux opérations cohomologiques).
\end{defi}

\begin{rema}
\label{remaL=L}
  S'il existe
  un diviseur $T$ de $X$ et $Z$ un sous-schéma fermé de $X$ tels que
  $U = Z \setminus T$, nous avions défini la catégorie
  $F\text{-}\mathfrak{C} ^\vee  _{U}$ (voir \cite[3.3.1]{caro_surcoherent}).
On retrouve alors, modulo l'inclusion
  $F\text{-}D ^{\mathrm{b} \vee} _\mathrm{surhol} (\D _{U })\subset F\text{-}\mathfrak{C} ^\vee  _{U}$,
  les fonctions $L$ et $P$ de \cite[3.3.3 et 3.3.4]{caro_surcoherent}.
Nous avions établi l'égalité $L=P$ sur $F\text{-}\mathfrak{C} ^\vee  _{U}$ (voir \cite[3.3.8]{caro_surcoherent}).

En général, il existe toujours deux sous-schémas fermés $Z$ et $T$ de $X$ tels que
$U = Z \setminus T$ (on prend $Z =\overline{Y}$ et $T = Z \setminus U$).
Mais on ne peut pas toujours choisir $T$ égal à un diviseur.
Par exemple, si $U$ est le complémentaire d'un sous-schéma fermé de $X$ de codimension égale à $2$.
\end{rema}

\begin{lemm}
  \label{LPtridis}
  Soit
  $\E ' \rightarrow \E \rightarrow \E ''\rightarrow \E' [1]$ un triangle distingué de
$F\text{-}D ^{\mathrm{b} \vee} _\mathrm{surhol} (\D _{U })$.
On a alors les
deux égalités
\begin{align}
L (U ,\E , t) & =L (U ,\E ', t)\times L (U,\E '', t), \\
P (U ,\E , t) & =P (U ,\E ', t)\times P (U ,\E '', t).
\end{align}
\end{lemm}
\begin{proof}
De manière analogue à \cite[3.2.1]{caro_courbe-nouveau}, cela résulte du caractère multiplicatif en les suites exactes du déterminant.
\end{proof}

\begin{lemm}
\label{LPYY'}
Soient $j $ : $U ' \hookrightarrow U$ une immersion
de variétés,
$\E'\in F\text{-}D ^{\mathrm{b} \vee} _\mathrm{surhol} (\D _{U' })$.
On dispose alors des égalités :
\begin{align}
L(U ',\E' ,t)= & L(U , j _! (\E'),t),\\
P(U ',\E' ,t)= & P(U , j _! (\E'),t)\label{LPYY'2}.
\end{align}
\end{lemm}

\begin{proof}
L'égalité $P(U ',\E' ,t)= P(U , j _! (\E'),t)$ résulte de la transitivité
du foncteur image directe extraordinaire.
Traitons à présent l'égalité $L(U ',\E' ,t)= L(U , j _! (\E'),t)$.
Soient $Z$ (resp. $Z'$) l'adhérence de $U$ (resp. $U'$) dans $X$,
$T$ le complémentaire de $U$ dans $Z$. (resp. $T'$ le complémentaire de $U'$ dans $Z'$).
Alors, $j _! (\E') = \DD _T \DD _{T'} (\E ')$ et
$j^+ j_! (\E')= \DD _{T'} \R \underline{\Gamma} ^\dag _{Z'} (\hdag T') \DD _{T} (j _! (\E'))$.
On obtient ainsi l'isomorphisme $j^+ j_! (\E') \riso
\DD _{T'} \R \underline{\Gamma} ^\dag _{Z'} (\hdag T')\DD _{T'}  (\E')
\riso
\DD _{T'} \DD _{T'} (\E ') \riso \E'$ (on utilise l'isomorphisme de bidualité de \ref{bidua+rel}).
Pour tout point $y'$ de $U'$, notons $i' _{y'}$ l'immersion fermé
$\Spec k (y') \hookrightarrow U'$. Pour tout $y' \in U'$,
on en déduit $i^+ _{y'} (j _! (\E')) \riso
i^{\prime +} _{y'} (j^+ j_! (\E')) \riso i^{\prime +} _{y'} (\E')$.
Pour finir la preuve, il reste à établir que si $y \in U \setminus U'$ alors
$i^+ _{y} (j _! (\E')) =0$.
Or, dans ce cas, soit $y \not \in Z'$, soit $y \in T'$.
En notant $\iota_y\,:\, \Spf \V (y) \hookrightarrow \X$ un relèvement de l'immersion fermée canonique
$\Spec k(y) \hookrightarrow X$, on obtient,
via l'isomorphisme de bidualité, l'isomorphisme $i^+ _{y} (j _! (\E')) \riso \DD \iota_y ^! \DD _{T'} (\E ')$.
Le cas $y \not \in Z'$ résulte du fait que $\DD _{T'} (\E ')$ soit à support dans $Z'$.
Lorsque $y \in T'$, on obtient
$\DD \iota_y ^! \DD _{T'} (\E ')=\DD \iota_y ^! (\hdag T') \DD (\E ')=0$
car
$\iota_y ^! \circ  (\hdag T') \riso (\hdag \Spec k(y) ) \circ \iota_y ^! =0$ (voir \ref{commutfonctcohlocal2}).
\end{proof}

\begin{theo}\label{theoL=P}
Pour tout $ \E\in F\text{-}D ^{\mathrm{b} \vee} _\mathrm{surhol} (\D _{U })$,
l'égalité $L( U, \E, t )=P(U , \E, t ) $ est validée.
\end{theo}
\begin{proof}
Afin de se ramener au cas déjà étudié (voir \ref {remaL=L} ou \cite[3.3.8]{caro_surcoherent}),
on procède par récurrence sur la dimension de $U$. Lorsque $\dim U = 0$, cela découle de
\cite[3.3.8]{caro_surcoherent} et de l'égalité
$F\text{-}D ^{\mathrm{b} \vee} _\mathrm{surhol} (\D _{U})=F\text{-}\mathfrak{C} ^\vee  _{U}$.

Soit $U'$ un ouvert affine, lisse et dense dans $U$ et $Y := U \setminus U'$ le sous-schéma fermé réduit induit. Notons
$i$ : $Y \hookrightarrow U$ et $j$ : $U ' \subset U$ les morphismes canoniques.
En dualisant \ref{defitrgllocsch} (i.e., en appliquant le foncteur dual
au triangle de localisation de \ref{defitrgllocsch} appliqué au dual de $\E$), on obtient le triangle de localisation
\begin{equation}
\label{dual-trilocUUiY}
j_! j ^+ (\E)  \rightarrow \E \rightarrow u _! u ^+ (\E)
\rightarrow j_! j ^+ (\E) [1].
\end{equation}
Comme $U'$ est affine et lisse, il existe un $\V$-schéma formel propre et lisse $\X'$,
un diviseur $T'$ de $X'$ et $Z'$ un sous-schéma fermé de $X'$ tel que
$U' = Z '\setminus T'$ (voir le début de \ref{nota-rem-sch-sep}).
Grâce à \cite[3.3.8]{caro_surcoherent} et comme
$F\text{-}D ^{\mathrm{b} \vee} _\mathrm{surhol} (\D _{U '})\subset F\text{-}\mathfrak{C} ^\vee  _{U'}$,
on en déduit que $L( U', j ^+ (\E), t )=P(U ', j ^+ (\E), t ) $.
De plus, comme $\dim Y < \dim U$, par hypothèse de récurrence,
$L( Y, u ^+ (\E), t ) =P(Y, u ^+ (\E), t ) $.
On conclut grâce à \ref{LPtridis}, \ref{dual-trilocUUiY}.

\end{proof}

\section{Surholonomie des $F$-isocristaux surconvergents unités}

Soient $\PP$ un $\V$-schéma formel séparé et lisse,
$X$ un sous-schéma fermé lisse de $P$, $T$ un diviseur de $P$ tel que $T _X :=T
\cap X$ soit un diviseur de $X$. On note $Y := X \setminus T _X$,
$u_0\,:\, X \hookrightarrow P$ l'immersion fermée canonique.

\begin{lemm}\label{prop-desc}
  Soient $f_0$ : $X' \rightarrow X$ un morphisme projectif, surjectif génériquement fini et étale de $k$-schémas lisses
  tel que $f_0 ^{-1} (T _X)$ soit un diviseur de $X'$.
  Soient $E$ un isocristal sur $Y$ surconvergent le long de $T _X$
  et $E ' := f_0 ^* (E)$ l'isocristal sur $X' \setminus f_0 ^{-1} (T _X)$ surconvergent le long de $f_0 ^{-1} (T _X)$
  déduit.

\begin{enumerate}
  \item \label{prop-desc(1)} Il existe alors un morphisme propre et lisse $g$ : $\PP' \rightarrow \PP$ et
  une immersion fermée $u' _0$ : $X ' \hookrightarrow P'$ telle que $g _0 \circ u' _0 = u _0 \circ f_0$.

  \item \label{prop-desc(2)} Pour tout $n \in \N \cup \{ \infty\}$,
  en notant $T' := g _0 ^{-1} (T)$,
  si
  $\sp _{X '\hookrightarrow \PP', T ',+ }(E')$ (défini dans \ref{Cas-comp-lisse})
  est $n$-surholonome
  (resp. $\sp _{X '\hookrightarrow \PP', T ',+ }(E')$ et $\DD \sp _{X '\hookrightarrow \PP', T ',+ } (E')$ sont $n$-surholonomes)
  alors
  $\sp _{X \hookrightarrow \PP, T ,+ }(E)$
  (resp. $\sp _{X \hookrightarrow \PP, T ,+ }(E)$ et $\DD \sp _{X \hookrightarrow \PP, T ,+ } (E)$)
  est
  $n$-surholonome (resp. sont $n$-surholonomes).
\end{enumerate}
\end{lemm}

\begin{proof}
L'assertion \ref{prop-desc(1)} est immédiate :
on choisit une immersion fermée $X ' \hookrightarrow \P ^r _X$,
dont le composée avec la projection $\P ^r _X \rightarrow X$ donne $f_0$ puis on prend
$u' _0$ égal au composé $X ' \hookrightarrow \P ^r _X \hookrightarrow \P ^r _P$,
$\PP '$ égal à l'espace projectif formel de dimension $r$ au-dessus de $\PP$ et $g$ : $\PP '
\rightarrow \PP$ égal à la projection.
Notons $\E :=\sp _{X \hookrightarrow \PP, T ,+ }(E)$
et $\E ':=\sp _{X' \hookrightarrow \PP', T ',+ }(E')$.

Traitons maintenant \ref{prop-desc(2)}.
Comme la $n$-surholonomie en tant que $\D ^\dag _{\PP,\Q}$-module est
locale en $\PP$, on peut supposer $\PP$ affine.
Dans ce cas, $u_0$ se relève en une immersion fermée
$u\,:\, \X \hookrightarrow \PP$ de $\V$-schémas formels lisses.
Comme $u^! _T(\E)\riso \sp _{X \hookrightarrow \X,T,+ } (E)$
(voir \cite[2.8]{caro_unite} ou \cite[4.1.8]{caro-construction}),
via le théorème de Kashiwara
(voir \ref{kashiwara}),
$\E \riso u_{T,+} (\sp _{X \hookrightarrow \X,T,+ } (E))$.
En utilisant l'isomorphisme canonique $u_{T,+} \riso u_+$ (voir \ref{fT=f}),
le théorème de dualité relative (voir \ref{iso-dualrel}),
la préservation de la $n$-surholonomie par image directe par un morphisme propre,
on se ramène ainsi au cas où $X =P$.

Prouvons à présent que $\E$ est un facteur direct de $g _{ T,+} (\E')$.
Par \ref{f+adjf!}, on dispose du morphisme d'adjonction
$g _{ T +} \circ  g _T ^!  ( \E )\rightarrow \E$.
En composant avec
$g _{ T +} \circ  \R \underline{\Gamma} ^\dag _{X '} \circ g _T ^!  ( \E ) \rightarrow
g _{ T +} \circ  g _T ^!  ( \E )$, on obtient ainsi le morphisme canonique
\begin{equation}
\label{prop-desc-iso1}
  g _{ T +} \circ  \R \underline{\Gamma} ^\dag _{X '} \circ g _T ^!  ( \E )\rightarrow \E.
\end{equation}
En remplaçant, dans \ref{prop-desc-iso1}, $\E$ par $\DD _T (\E)$, il s'en suit le morphisme canonique
$g _{ T +} \circ  \R \underline{\Gamma} ^\dag _{X '} \circ g _T ^! \circ \DD _T (\E) \rightarrow \DD _T (\E)$.
En appliquant le foncteur dual $\DD _T $ à ce morphisme, via l'isomorphisme de bidualité
et l'isomorphisme de dualité relative (voir \ref{bidual} et \ref{iso-dualrel}), cela donne
\begin{equation}
  \label{prop-desc-iso2}
\E \rightarrow  g _{ T,+} \circ \DD _{T'} \circ  \R \underline{\Gamma} ^\dag _{X '}\circ  g  ^! _T \circ  \DD _T (\E ).
\end{equation}

Or, avec \cite[4.1.8]{caro-construction} et \cite[4.3.5]{caro-construction}, on dispose des isomorphismes :
$$\DD _{T'} \circ  \R \underline{\Gamma} ^\dag _{X '}\circ   g_T   ^! \circ \DD _T (\E )
\riso \sp _{X '\hookrightarrow \PP', T ',+ }((f_0 ^* (E ^\vee) )^\vee)
\riso
\sp _{X '\hookrightarrow \PP', T ',+ }(E')
\riso \R \underline{\Gamma} ^\dag _{X '}\circ   g _T  ^!  ( \E ).$$
Il en découle que le morphisme
\ref{prop-desc-iso2} est de la forme
$\E \rightarrow g _{ T +}  \circ  \R \underline{\Gamma} ^\dag _{X '} \circ  g _T ^!  ( \E )$.
En le composant avec \ref{prop-desc-iso1},
on obtient un morphisme canonique :
$\E \rightarrow \E,$ que l'on note $\theta$.
Pour vérifier que $\E$ est un facteur direct de $g _{ T,+} (\E')$, il suffit d'établir que
$\theta$ est un isomorphisme.

Par \ref{4312Be1},
il suffit de vérifier que $\theta$ est un isomorphisme en dehors de $T$. On se ramène ainsi à supposer
$T$ vide.
De plus, comme $\E = \sp _* (E)$ (où $\sp $ : $ \PP _K \rightarrow \PP$ est le morphisme de spécialisation),
d'après \ref{be4eqiso},
ce morphisme $\theta$ est l'image par $\sp _* $
d'un morphisme d'isocristaux surconvergents de la forme $E \rightarrow E$.
Comme le foncteur restriction à un ouvert dense est fidèle sur la catégorie des
isocristaux convergents (cela résulte par exemple \cite[4.1.3]{tsumono}),
pour vérifier que $\theta$ est un isomorphisme,
on peut alors en outre supposer que le morphisme $f _0$ est fini et étale.
Comme $X'$ est alors affine, l'immersion fermée $u' _0$ : $X ' \hookrightarrow P'$ se relève en une immersion
fermée de $\V$-schémas formels lisses $u' $ : $\X ' \hookrightarrow \PP'$.

Comme $u ^{\prime !} \circ g  ^!  ( \E ) \riso \sp _* (f _0^* (E))$
est un isocristal convergent sur $X'$ (i.e., via \ref{be4eqiso}, est un $\D ^\dag _{\X',\Q}$-module cohérent
$\O _{\X',\Q}$-cohérent), c'est en particulier un $\D ^\dag _{\X',\Q}$-module cohérent.
On vérifie alors (de manière analogue à \cite[3.1.12]{caro_surcoherent})
que l'on dispose d'un isomorphisme canonique (d'ailleurs compatible à Frobenius)
$u ' _+ u ^{\prime !} ( g  ^!  ( \E ))
\riso \R \underline{\Gamma} ^\dag _{X'} (g  ^!  ( \E ))$ s'inscrivant dans le diagramme canonique :
\begin{equation}
  \label{+!=gamma2}
  \xymatrix @R=0,3cm {
     {u ' _+ u ^{\prime !} ( g  ^!  ( \E ))}
     \ar[rr]^-{\mathrm{adj}} \ar[d] _-\sim && {\E} \ar@{=}[d] \\
  {\R \underline{\Gamma} ^\dag _{X'  } (g  ^!  ( \E ))}
  \ar[rr]  && {\E},}
\end{equation}
où le morphisme du haut est le morphisme d'adjonction de $u'$ (via \ref{f+adjf!}),
celui du bas est le morphisme canonique (par fonctorialité en $X'$ de $\R \underline{\Gamma} ^\dag _{X'  } $).

Le morphisme
$g _{ +} \circ  \R \underline{\Gamma} ^\dag _{X '} \circ g  ^!  ( \E )\rightarrow \E$
est ainsi canoniquement isomorphe au morphisme induit par adjonction entre les foncteurs images directes et images
inverses extraordinaires par un morphisme propre (à savoir $u'$ et $g$) :
$g _{ +} \circ  u '_+ \circ u^{\prime !} \circ g  ^!  ( \E )\rightarrow \E$.
Par transitivité de ces morphismes d'adjonction (voir \ref{dual-relasection} ou plus précisément
\cite[1.2.11]{caro-construction}),
celui-ci est canoniquement isomorphe au morphisme d'adjonction $f _+ f ^! ( \E ) \rightarrow \E$, où $f = g \circ u '$.
Il en résulte, par construction de \ref{prop-desc-iso2} à partir de
\ref{prop-desc-iso1},
que $\E \rightarrow g _{  +}  \circ  \R \underline{\Gamma} ^\dag _{X '} \circ  g  ^!  ( \E )$
est canoniquement isomorphisme au morphisme canonique d'adjonction
$\E \rightarrow f _! f ^+ (\E) \riso f _+ \circ f ^! (\E)$.

Comme $f_0$ est fini et étale, on dispose d'un foncteur image directe $f _{0*}$ de la catégorie des
isocristaux convergents sur $X'$ vers les isocristaux convergents sur $X$.
Via l'équivalence de catégories \ref{be4eqiso} (ici, les diviseurs sont vides),
$f _{0*}$ (resp. $f ^* _0$) correspond à $f _+$ ou $f _!$
(resp. $f ^!$ ou $f ^+$).
De plus, le morphisme d'adjonction $f _{0*} f _{0} ^* (E) \rightarrow E$
(resp. $E\rightarrow f _{0*} f _{0} ^* (E) $) donne via cette correspondance le morphisme d'adjonction
$f _+ f ^! ( \E ) \rightarrow \E$ (resp. $\E \rightarrow f _+ f ^! ( \E ) $).
Comme le morphisme composé $E \rightarrow f _{0*} f _{0} ^* (E) \rightarrow E$ est
égal à la multiplication par le degré de $X'$ sur $X$, il en résulte que $\theta$ est
un isomorphisme.

Ainsi, $\E$ est un facteur direct de $g _{ T,+} (\E')$.
Si $\E'$ est $n$-surholonome, comme la $n$-surholonomie est préservée par l'image directe d'un morphisme propre,
cela implique que
$g _{ T,+}(\E') \riso g _{ +}   (\E')$ (voir \ref{fT=f})
est $n$-surholonome.
Par \ref{stab-surhol-Dual-hdag}.\ref{stab-surhol-Dual-hdag2},
$\E$ est donc $n$-surholonome.

Si $\E'$ et $\DD (\E')$ sont $n$-surholonomes, $\E$ est donc $n$-surholonome.
De plus, $\DD (\E)$ est un facteur direct de $\DD \circ g _{ T,+} (\E')$.
Comme $\E'$ est $n$-surholonome, en particulier $\E'$ est $\D ^\dag _{\PP,\Q}$-cohérent
et on bénéficie donc de l'isomorphisme de dualité relative (voir \ref{iso-dualrel}) :
$\DD \circ g _{ +}  (\E') \riso g _{ +} \circ \DD (\E')$.
Or, d'après \ref{fT=f},
$g _{ T,+}(\E') \riso g _{ +}   (\E')$.
Ainsi, $\DD (\E)$ est un facteur direct de $g _{ +} \circ \DD (\E')$.
La proposition \ref{stab-surhol-Dual-hdag}.\ref{stab-surhol-Dual-hdag2} nous permet alors de conclure.
\end{proof}

\begin{rema}
  \label{rema-prop-desc}
$\bullet$ Avec les notations de \ref{prop-desc},
comme
$\sp _{X' \hookrightarrow \PP', T ',+ }(E')\riso \R \underline{\Gamma} ^\dag _{X '} g _T  ^!  ( \sp _{X \hookrightarrow \PP, T ,+ }(E) )$
(voir \cite[2.8]{caro_unite} ou \cite[4.1.8]{caro-construction}),
par \ref{stab-surhol-Dual-hdag} et \ref{stabinminv},
si $\sp _{X \hookrightarrow \PP, T ,+ }(E)$ est $n$-surholonome
alors $\sp _{X' \hookrightarrow \PP', T ',+ }(E')$ est $n$-surholonome.

$\bullet$ Lorsque $E$ est en outre muni d'une structure de Frobenius, on dispose alors d'un résultat analogue à
  \ref{prop-desc} dans un cadre plus général (voir \cite[6.3.1]{caro_devissge_surcoh}).
  En effet, la structure de Frobenius permet d'utiliser le théorème de pleine fidélité de Kedlaya
  (voir \cite{kedlaya_full_faithfull}).
\end{rema}

\begin{prop}
\label{proptheo}
Pour tout isocristal convergent $G$ sur $X$,
le $\D ^\dag _{\PP ,\Q}$-module $\sp _{X \hookrightarrow \PP,+} (G)$
est surholonome.
\end{prop}
\begin{proof}
Comme $X$ est somme de ses composantes irréductibles, il ne coûte rien de
supposer $X$ irréductible.
Par récurrence sur $n \geq 0$,
prouvons la propriété $(P_n)$ suivante :

$(P _n )\, :\,$ {\og Pour tout isocristal convergent $G$ sur $X$
($X \hookrightarrow \PP$ pouvant naturellement varier),
le $\D ^\dag _{\PP ,\Q}$-module cohérent $\sp _{X \hookrightarrow \PP,+} (G)$ associé
est $n$-surholonome\fg}.

L'assertion $P _0$ est un cas particulier de \cite[3.1.2.2]{caro_surcoherent} (ou alors,
on remarque qu'il suffit d'adapter la preuve de $P _n\Rightarrow P_{n+1}$, i.e. on enlève les foncteurs $\DD$,
pour établir $P_0$).
Pour $n \geq 0$, en supposant $P _n$ vraie, démontrons $P_{n+1}$.

Par \cite[4.1.8]{caro-construction},
si $g$ : $\PP ' \rightarrow \PP$ est un morphisme lisse de $\V$-schémas
formels lisses, $X ':= g _0 ^{-1}(X)$,
$f$ : $ X' \rightarrow X$ est le morphisme de $k$-schémas lisses
induit par $g$ et $d _{X'/X}$ est la dimension relative de $X'$ sur $X$,
alors, pour tout isocristal convergent $G$ sur $X$,
$\sp _{X '\hookrightarrow \PP',+} (f ^*  (G))[d _{X'/X}]
\riso \R \underline{\Gamma} ^\dag _{X '} \circ g  ^! ( \sp _{X \hookrightarrow \PP,+} (G))$.
Or, comme $g  ^! ( \sp _{X \hookrightarrow \PP,+} (G))$ est à support dans $X'$, le morphisme
canonique
$\R \underline{\Gamma} ^\dag _{X '} \circ g  ^! ( \sp _{X \hookrightarrow \PP,+} (G))
\rightarrow
g  ^! ( \sp _{X \hookrightarrow \PP,+} (G))$ est un isomorphisme. D'où,
$\sp _{X' \hookrightarrow \PP',+} (f ^*  (G))[d _{X'/X}] \riso g  ^! ( \sp _{X \hookrightarrow \PP,+} (G))$.
Il en résulte qu'il suffit de prouver que pour tout isocristal convergent $G$ sur $X$,
pour tout diviseur $T$ de $P$,
$\DD \circ (\hdag T) \circ \sp _{X \hookrightarrow \PP,+} (G)$ est $n$-surholonome.
De plus, puisque $X$ est irréductible, soit $T \cap X =X$ soit $T \cap X$ est un diviseur de $X$.
Le cas où $T \cap X =X$ implique alors $\DD \circ (\hdag T) \circ \sp _{X \hookrightarrow \PP,+} (G)=0$,
qui $n$-surholonome. On se ramène donc au cas où $T \cap X$ est un diviseur de $X$.

Comme le fait que, pour tout isocristal convergent $G$ sur $X$,
pour tout diviseur $T$ de $P$,
$\DD \circ (\hdag T) \circ \sp _{X \hookrightarrow \PP,+} (G)$ soit $n$-surholonome est locale en $\PP$,
on se ramène au cas où $\PP$ est affine.
L'immersion fermée $X \hookrightarrow P$ se relève alors en
un morphisme $u$ : $\X \hookrightarrow \PP$ de $\V$-schémas formels lisses.
Dans ce cas, $\sp _{X \hookrightarrow \PP,+} (G) \riso u _+ \circ \sp _* (G)$,
où $\sp $ : $ \X _K \rightarrow \X$ est le morphisme de spécialisation.
Comme le foncteur $u_+$ commute aux foncteurs de localisation (voir \ref{commutfonctcohlocal2+})
et au foncteur dual (voir \ref{iso-dualrel}), on obtient donc
$\DD \circ (\hdag T) \circ \sp _{X \hookrightarrow \PP,+} (G)
\riso
u _+ \circ \DD \circ (\hdag T \cap X) \circ \sp _* (G)$.

La $n$-surholonomie étant stable
par image directe par un morphisme propre de $\V$-schémas formels lisses,
on se ramène ainsi au cas où $X =P$.
On pose alors $\X :=\PP$ et $\G= \sp _* (G)$.
Il s'agit de prouver que $\DD \circ (\hdag T) (\G)$ est $n$-surholonome.

i) Supposons dans un premier temps que $T =\emptyset$.
D'après \cite{caro_comparaison}, le foncteur $\sp _*$ commute aux foncteurs duaux, i.e.,
$\DD \circ \sp _* (G) \riso
 \sp _* (G ^\vee )$, où $G ^\vee $ désigne l'isocristal convergent sur $X$ dual de $G$.
 On conclut alors par hypothèse de récurrence.
 \medskip

ii) Traitons maintenant le cas où $T$ est un diviseur lisse. Comme $T$
est affine et lisse, il existe un relèvement de $\V$-schémas
formels affines et lisses $i$ : $\mathfrak{T} \hookrightarrow \X$
de l'immersion fermée $T \hookrightarrow X$. Or, en appliquant le foncteur dual $\DD$
au triangle de localisation de $\G$ en $T$, on obtient le triangle distingué
\begin{equation}\label{tria-loc-proptheo}
\DD (\hdag T) (\G) \rightarrow \DD (\G )\rightarrow   \DD \circ \R \underline{\Gamma} ^\dag _T  (\G)
\rightarrow \DD (\hdag T) (\G)  [1].
\end{equation}
On sait que $i _ + \circ i ^! (\G) \riso \R \underline{\Gamma} ^\dag _T  (\G) $ (voir \ref{445Beintro}
ou \cite[4.4.5.2]{Beintro2}).
Comme $i ^! (\G)$ est cohérent (car isomorphe à $\sp _* ( i ^* G) [-1] $),
en utilisant le théorème de dualité relative (voir \ref{iso-dualrel}), il en découle
$\DD \circ  \R \underline{\Gamma} ^\dag _T  (\G) \riso i _ + \circ \DD \circ i ^! (\G) $.
Or, d'après respectivement \cite[4.3.5]{caro-construction} et \cite[4.1.8]{caro-construction},
le foncteur $\sp _*$ commute aux images inverses (extraordinaire)
et aux foncteurs duaux. Ainsi,
$\DD \circ i ^! (\G) \riso  \sp _* (i ^*  (G ^\vee ) )[1]$.
D'où :
$\DD\circ  \R \underline{\Gamma} ^\dag _T  (\G) \riso
i _ + ( \sp _*  i ^* (G ^\vee ) [1])$.
Comme l'hypothèse de récurrence implique que $\sp _* (i ^*  (G ^\vee ) )$ est $n$-surholonome,
il découle alors du théorème \ref{stabinmdir} que
$i _ + ( \sp _*  i ^* (G ^\vee ) [1])$ est $n$-surholonome et donc
que $\DD\circ  \R \underline{\Gamma} ^\dag _T  (\G)$ est $n$-surholonome.

Il résulte alors du triangle distingué
\ref{tria-loc-proptheo} et du cas i), que le faisceau
$\DD (\hdag T) (\G)$ est $n$-surholonome.
\medskip

iii) Supposons à présent que $T$ soit un diviseur à croisements
normaux stricts, i.e. (voir \cite[2.4]{dejong}), un diviseur tel que
$T$ soit non vide, pour tout $x \in T$, l'anneau local $\O _{T,\,x}$ soit régulier,
$T$ soit un schéma réduit, autrement dit $T = \cup _{i=1} ^r T _i$
avec $T _i$ composantes irréductibles et distinctes de $T$, et
pour tout sous-ensemble non-vide $J \subset I$, le sous-schéma fermé
$T _J := \cap _{j\in J} T _j$
soit régulier et de codimension dans $X$ égale au cardinal de $J$.

Effectuons alors une (deuxième) récurrence sur $r\geq 1$, où $r$ est l'entier qui
apparaît dans l'expression $T = \cup _{i=1} ^r T _i$.
Le cas où $r =1$ a été traité dans ii). Considérons le triangle distingué
de Mayer-Vietoris (\ref{eq1mayer-vietoris}) :
\begin{equation}\label{mayer-vietorisproptheo1}
  \R \underline{\Gamma} ^\dag _{ \cup _{i=2} ^r T _1 \cap T _i} \ (\G)
  \rightarrow
\R \underline{\Gamma} ^\dag _{ T _1  } (\G)
\oplus
\R \underline{\Gamma} ^\dag _{ \cup _{i=2} ^r T _i} \ (\G)
  \rightarrow
\R \underline{\Gamma} ^\dag _{ \cup _{i=1} ^r T _i} \ (\G)
\rightarrow
\R \underline{\Gamma} ^\dag _{ \cup _{i=2} ^r T _1 \cap T _i}  (\G) [1].
\end{equation}
Notons $i _1$ : $\mathfrak{T} _1 \hookrightarrow \X$,
un relèvement de l'immersion fermée $T _1 \hookrightarrow X$.
Comme $\R \underline{\Gamma} ^\dag _{ T _1 } \riso i_{1+} i_1 ^!$ (voir \ref{445Beintro}),
par \ref{GammaZGammaZ'},
on obtient l'isomorphisme
$\R \underline{\Gamma} ^\dag _{ \cup _{i=2} ^r T _1 \cap T _i} \ (\G)
\riso i _{1+}  i _1 ^!  \R \underline{\Gamma} ^\dag _{ \cup _{i=2} ^r T _1 \cap T _i} \ (\G) $.
Comme l'image inverse extraordinaire commute au foncteur cohomologique local (voir \ref{commutfonctcohlocal2}),
cela implique l'isomorphisme :
$\R \underline{\Gamma} ^\dag _{ \cup _{i=2} ^r T _1 \cap T _i} \ (\G)
\riso i _{1+} \circ \R \underline{\Gamma} ^\dag _{ \cup _{i=2} ^r T _1 \cap T _i} \ ( i _1 ^! \G)$.
Comme $\cup _{i=2} ^r T _1 \cap T _i$ est un diviseur à croisements normaux stricts de
$T _1$ et puisque que $i _1 ^! (\G) \riso \sp _* (i _1 ^* (G))[-1]$,
l'hypothèse de récurrence sur $r$
implique que
$\DD \circ ( \hdag \cup _{i=2} ^r T _1 \cap T _i)  \ (i _1 ^! (\G))$
est $n$-surholonome.
Il résulte alors du cas $i)$ (via un triangle de localisation de
$i _1 ^! (\G)$ en $\cup _{i=2} ^r T _1 \cap T _i$) que
$\DD \circ \R \underline{\Gamma} ^\dag _{ \cup _{i=2} ^r T _1 \cap T _i} \ (i _1 ^! (\G))$
est $n$-surholonome.
Grâce au théorème de dualité relative (voir \ref{iso-dualrel}) et à la préservation de la $n$-surholonomie
par image directe par un morphisme propre (voir \ref{stabinmdir}),
il s'en suit que
$\DD \circ  \R \underline{\Gamma} ^\dag _{ \cup _{i=2} ^r T _1 \cap T _i} \ (\G)
\riso i _{1+} \circ \DD  \circ \R \underline{\Gamma} ^\dag _{ \cup _{i=2} ^r T _1 \cap T _i} \ ( i _1 ^! \G)$
est $n$-surholonome.
On obtient en outre (via aussi un triangle de localisation) par hypothèse de récurrence sur $r$ que
$\DD \circ \R \underline{\Gamma} ^\dag _{ T _1  } (\G)
\oplus
\DD \circ \R \underline{\Gamma} ^\dag _{ \cup _{i=2} ^r T _i} \ (\G)$ est $n$-surholonome.
En appliquant le foncteur $\DD$ au triangle distingué \ref{mayer-vietorisproptheo1},
il en dérive que le complexe $\DD \circ \R \underline{\Gamma} ^\dag _T (\G)$ est $n$-surholonome.
Puis, en appliquant $\DD$ au triangle distingué de localisation en $T$ de $\G$, via le cas i),
il en résulte que
$\DD \circ (\hdag T) (\G)$ est $n$-surholonome.
\medskip

iv) Enfin, passons au cas où $T $ est un diviseur quelconque mais différent de l'ensemble vide.
Grâce au théorème de de Jong
sur les altérations de variétés algébriques \cite{dejong}, il
existe un $k$-schéma lisse $X '$ et un morphisme $f _0$ : $X '
\rightarrow X$ projectif, surjectif et génériquement fini et étale,
tels que $f _0 ^{-1} (T)$ soit un diviseur à croisements normaux stricts
de $X '$. Désignons par $j^\dag G$ l'isocristal sur $Y$ surconvergent le long de $T _X$
canoniquement déduit de $G$. Comme
$\DD  \sp _{X \hookrightarrow \PP,T, +} (j ^\dag G)\riso \DD (\hdag T) (\G)$,
il suffit de prouver que
$\DD  \sp _{X \hookrightarrow \PP,T, +} (j ^\dag G)$
est $n$-surholonome.

D'après \ref{prop-desc}.\ref{prop-desc(1)}, il existe un morphisme propre et lisse $g$ : $\PP' \rightarrow \PP$ et
  une immersion fermée $u' _0$ : $X ' \hookrightarrow P'$ telle que $g _0 \circ u' _0 = u _0 \circ f_0$.
Or, en notant $T' := g _0 ^{-1} (T)$,
 $f _0 ^* j ^\dag G$ est un isocristal sur $X ' \setminus T'$ surconvergent le long du
diviseur à croisements normaux stricts
$T' \cap X'$ de $X '$
et
$\sp _{X '\hookrightarrow \PP',T', +} (f _0 ^* j ^\dag G)
\riso
(\hdag T')\sp _{X '\hookrightarrow \PP',T', +} (f _0 ^* G)$
Le cas iii) implique alors que
$\sp _{X '\hookrightarrow \PP',T', +} (f _0 ^* j ^\dag G)$
et
$\DD \sp _{X '\hookrightarrow \PP',T', +} (f _0 ^* j ^\dag G)$
sont $n$-surholonomes.
Il découle alors de \ref{prop-desc}.\ref{prop-desc(2)} que
$\DD \sp _{X \hookrightarrow \PP,T, +} (j ^\dag G)$
est $n$-surholonome.

\end{proof}

\begin{theo}
\label{unitesurho}
Soit $E$ un $F$-isocristal unité sur $Y$ surconvergent le long de $T _X$.
Le faisceau $\sp _{X \hookrightarrow \PP, T,+} (E)$ est surholonome.
\end{theo}
\begin{proof}
On peut supposer $X$ intègre.
  D'après le théorème de monodromie génériquement finie de
Tsuzuki \cite[1.3.1]{tsumono},
il existe un $k$-schéma lisse $X '$ et un morphisme $f$ :
$X' \rightarrow X$ propre, surjectif et génériquement fini et étale, tels que,
si on pose
$Y '= f ^{-1} ( Y)$, $T ' = X ' \setminus Y '$ et $j'$
l'immersion ouverte $Y '\hookrightarrow X'$, alors il existe $G'$,
un (unique) $F$-isocristal convergent sur $X'$ vérifiant
l'isomorphisme
$f ^*  (E) \riso (j ') ^\dag G '$
de $F$-isocristaux sur $Y'$ surconvergents le long de $T'$.
En utilisant le théorème de désingularisation de de Jong \cite{dejong},
on peut en outre supposer que le morphisme $f$ est projectif et que $T '$ est un
diviseur de $X '$.

D'après \ref{prop-desc}, il existe alors un morphisme propre et lisse $g$ : $\PP' \rightarrow \PP$ et
  une immersion fermée $u' _0$ : $X ' \hookrightarrow P'$ telle que $g _0 \circ u' _0 = u _0 \circ f$.
En notant $T '' := g _0 ^{-1} (T)$,
comme $\sp _{X' \hookrightarrow \PP', T'',+} (j ') ^\dag (G ')
\riso
(\hdag T'') \sp _{X' \hookrightarrow \PP', T'',+} (G ')$,
il résulte de \ref{proptheo} et de \ref{stab-surhol-Dual-hdag}.\ref{stab-surhol-Dual-hdag1} que
$\sp _{X' \hookrightarrow \PP', T'',+} (j ') ^\dag (G ')$
est surholonome.
Puisque $f ^*  (E) \riso (j ') ^\dag G '$, il dérive de \ref{prop-desc} que
$\sp _{X \hookrightarrow \PP, T,+} (E)$ est surholonome.

\end{proof}

\begin{theo}
  \label{unitesurho2}
Soient $U $ un $k$-schéma séparé et lisse, $E$ un $F$-isocristal unité sur $U$ et $\E:= \sp _{U,+} (E)$
l'objet de $F\text{-}\mathrm{Isoc} ^{\dag \dag}( U/K )$ associé (voir \ref{eqcat-gen}).
Alors $\E \in F\text{-}\mathfrak{M} ^+ _{U }$, i.e., $\E$ est un
$F\text{-}\D_U$-module arithmétique surholonome.
\end{theo}

\begin{proof}
  Par construction de $\sp _{U,+}$, de $F\text{-}\mathrm{Isoc} ^{\dag \dag}( U/K )$ et de $F\text{-}\mathfrak{M} ^+ _{U }$
(constructions qui s'établissent par recollement via un recouvrement de $U$ par des ouverts affines), on peut supposer $U$ affine et intègre.
Il existe donc un plongement de
$U$ dans un $\V$-schéma formel $\X$ propre et lisse,
un diviseur $T$ de $X$ tels que $U = Z \setminus T$, où $Z$ est l'adhérence de $U$ dans $P$
(voir par exemple le début de \ref{nota-rem-sch-sep}).
D'après le théorème de désingularisation de de Jong de \cite{dejong}, il
existe un $k$-schéma lisse $Z '$ et un morphisme $f _0$ : $Z '
\rightarrow Z$ projectif, surjectif et génériquement fini et étale,
tels que $f _0 ^{-1} (T \cap Z)$ soit un diviseur à croisements normaux stricts
dans $Z '$. De manière analogue à \ref{prop-desc}.\ref{prop-desc(1)} , il existe un
 morphisme propre et lisse $g$ : $\X' \rightarrow \X$ et
  une immersion fermée $u' _0$ : $Z ' \hookrightarrow P'$ telle que $g _0 \circ u' _0 = u _0 \circ f_0$.
  Notons $U' := Z ' \setminus f _0 ^{-1} (T \cap Z) $, $\Y:= \X \setminus T$, $\Y':= \X ' \setminus g _0^{-1} (T)$,
 $h\,:\, \Y' \rightarrow \Y$ le morphisme induit par $g$,
  $E'$ le $F$-isocristal surconvergent unité sur $U'$ déduit de $E$,
  $\widehat{E}$ (resp. $\widehat{E'}$) le $F$-isocristal convergent sur $U$ (resp. $U'$)
  déduit de $E$ (resp. $E'$).

D'après \cite[6.3.1]{caro_devissge_surcoh},
$\R \underline{\Gamma} ^\dag _{Z' } g ^! (\E)$ est dans l'image essentielle de
$\sp _{Z '\hookrightarrow \X', g _0 ^{-1} (T ),+}$ (pour plus de précisions, voir la remarque \ref{rema-unitesurho2}).
De plus, on dispose des isomorphismes :
\begin{equation}
\label{unitesurho2-iso}
  \R \underline{\Gamma} ^\dag _{Z' } g ^! (\E) |\Y' \riso \R \underline{\Gamma} ^\dag _{U'}  h ^! (\E|\Y) \riso
  \R \underline{\Gamma} ^\dag _{U'}  h ^!
(  \sp _{U \hookrightarrow \Y , +} ( \widehat{E}))
\riso \sp _{U '\hookrightarrow \Y ', +} ( \widehat{E'}),
\end{equation}
le dernier isomorphisme provenant
de \cite{caro_unite} ou de \cite[4.1.8]{caro-construction}.

  Grâce au théorème de pleine fidélité de Kedlaya du foncteur restriction en dehors du lieu de surconvergence
  (voir \cite{kedlaya_full_faithfull})
  et à la pleine fidélité du foncteur $\sp _{Z '\hookrightarrow \X', g _0 ^{-1} (T ),+}$,
  pour vérifier que l'on dispose d'un isomorphisme de la forme
  $\R \underline{\Gamma} ^\dag _{Z' } g ^! (\E) \riso \sp _{Z '\hookrightarrow \X', g _0 ^{-1} (T ),+} (E')$,
  il suffit alors de le voir en dehors de $g _0 ^{-1} (T )$, ce qui découle de \ref{unitesurho2-iso}.
Cela implique, d'après \cite[6.3.1]{caro_devissge_surcoh}, que $\E$ est un facteur direct de
$g _+ \sp _{Z '\hookrightarrow \X', g _0 ^{-1} (T ),+} (E')$.
Via le théorème \ref{unitesurho}, $\sp _{Z '\hookrightarrow \X', g _0 ^{-1} (T ),+} (E')$ est
 surholonome.
 Les propriétés de stabilité de la surholonomie (voir \ref{stab-surhol-Dual-hdag}.\ref{stab-surhol-Dual-hdag2} et \ref{stabinmdir})
 nous permettent de conclure.
\end{proof}

\begin{rema}
\label{rema-unitesurho2}
Avec les notations de la preuve de \ref{unitesurho2},
il résulte de la description de l'image essentielle du foncteur
  $\sp _{Z '\hookrightarrow \X', g _0 ^{-1} (T ),+}$ de \cite[4.1.9]{caro-construction}
  et de celle de la dernière phrase du théorème de Berthelot de \ref{be4eqiso},
  qu'un $\D ^\dag _{\X'} ( g _0 ^{-1} (T )) _{\Q}$-cohérent et à support dans $Z'$
  est dans l'image essentielle de
  $\sp _{Z '\hookrightarrow \X', g _0 ^{-1} (T ),+}$ si et seulement si sa restriction en dehors
  de $g _0 ^{-1} (T )$ est dans l'image essentielle de $\sp _{U '\hookrightarrow \Y ', +}$.

Ainsi, (si ce fait peut sembler quelque peu abrupte dans la preuve de \cite[6.3.1]{caro_devissge_surcoh}),
pour voir que $\R \underline{\Gamma} ^\dag _{Z' } g ^! (\E)$ est dans l'image
essentielle de $\sp _{Z '\hookrightarrow \X', g _0 ^{-1} (T ),+}$, il suffit alors d'invoquer l'isomorphisme :
$\R \underline{\Gamma} ^\dag _{Z' } g ^! (\E) |\Y'
\riso \sp _{U '\hookrightarrow \Y ', +} ( \widehat{E'}).$

\end{rema}

\section{Sur les conjectures de la stabilité de l'holonomie}
Soit $f$ : $\X \rightarrow \Y$ un morphisme de $\V$-schémas formels lisses.
Considérons les conjectures suivantes :
\begin{enumerate}
\item[A)] On suppose que $\X$ se plonge dans un $\V$-schéma formel propre et lisse.
Si $\E$ est un objet de $F\text{-}D ^\mathrm{b} _\mathrm{hol} ( \D ^\dag _{\X,\,\Q})$
à support propre sur $\Y$,
alors l'image directe $f _+ (\E)$ appartient à $F\text{-}D ^\mathrm{b} _\mathrm{hol} ( \D ^\dag _{\Y,\,\Q})$.
\item[A')] Si $f$ est propre et
$\E$ est un objet de $F\text{-}D ^\mathrm{b} _\mathrm{hol} ( \D ^\dag _{\X,\,\Q})$,
alors $f _+ (\E)$ appartient à $F\text{-}D ^\mathrm{b} _\mathrm{hol} ( \D ^\dag _{\Y,\,\Q})$.
\item[B)] Si $\FF$ est un objet de $F\text{-}D ^\mathrm{b} _\mathrm{hol} ( \D ^\dag _{\Y,\,\Q})$,
alors l'image inverse extraordinaire $f ^! (\FF)$
appartient à $F\text{-}D ^\mathrm{b} _\mathrm{hol} ( \D ^\dag _{\X,\,\Q})$.
\item[C)] Si $Z$ est un fermé de $X$ et si $\E$ appartient $F\text{-}D ^\mathrm{b} _\mathrm{hol} ( \D ^\dag _{\Y,\,\Q})$,
alors $\R \underline{\Gamma} ^\dag _Z (\E)$ appartient à $F\text{-}D ^\mathrm{b} _\mathrm{hol} ( \D ^\dag _{\Y,\,\Q})$.
\item[C')] Si $Z$ est un fermé de $X$ et si $\E$ appartient $F\text{-}D ^\mathrm{b} _\mathrm{hol} ( \D ^\dag _{\Y,\,\Q})$,
alors $(\hdag Z ) (\E)$ appartient à $F\text{-}D ^\mathrm{b} _\mathrm{hol} ( \D ^\dag _{\Y,\,\Q})$.
\end{enumerate}
Rappelons que dans \cite[5.3.6]{Beintro2}, Berthelot a énoncé les conjectures
$A$ (sans l'hypothèse {\og $\X$ se plonge dans un $\V$-schéma formel propre et lisse\fg}),
$B$ et $C$ ci-dessus. D'après \cite[5.3.6]{Beintro2}, les conjectures $B$ et $C$ sont équivalentes.
De plus, via le triangle de localisation en $Z$, on remarque que les conjecture $C$ et $C'$ sont équivalentes.

\begin{lemm}
  \label{lemm-conjB-eq}
  On suppose la conjecture $B$ vérifiée.
  Alors, pour tout $\V$-schéma formel lisse $\X$,
  un $F\text{-}\D ^\dag _{\X,\,\Q}$-module cohérent est holonome si et seulement s'il est à fibres extraordinaires finies (voir \ref{defi-fef}).
\end{lemm}

\begin{proof}
  Comme la conjecture $B$ est supposée vérifiée, il est clair qu'un
  $F\text{-}\D ^\dag _{\X,\,\Q}$-module holonome est à fibres extraordinaires finies.
Réciproquement, soit $\E$ un $F\text{-}\D ^\dag _{\X,\,\Q}$-module cohérent à fibres extraordinaires finies.
Prouvons par récurrence sur la dimension de son support, noté $Z$, que celui-ci est holonome.
  D'après \ref{lemmsurcoh=>holrem},
  $\E' :=\mathcal{H} ^0 \DD \circ \mathcal{H} ^0 \DD(\E)$ est le plus grand sous-$F\text{-}\D ^\dag _{\X,\,\Q}$-module
  holonome de $\E$.
  D'après \cite[5.3.5(ii)]{Beintro2}, le terme du milieu d'une suite exacte courte de
$F\text{-}\D ^\dag _{\X,\,\Q}$-modules cohérents est holonome si et seulement
si les deux autres termes le sont.
Pour obtenir l'holonomie de $\E$, il suffit alors de prouver celle de $\E / \E '$.
 Or, la conjecture $B$ étant supposée vraie, $\E'$ est à fibres extraordinaires finies.
Comme $\E'$ et $\E$ sont alors à fibres extraordinaires finies,
il en est donc de même de $\E / \E '$.
 Par hypothèse de récurrence, il suffit alors vérifier que la dimension du support de
 $\E /\E '$ est strictement inférieure à celle de $Z$. Or, on déduit de \ref{lemmsurcoh=>hol} que l'inclusion
 $\E '\subset \E$ est un isomorphisme au dessus d'un ouvert $\U$ de $\X$ induisant un ouvert dense de $Z$.
 Par conséquent, $\E / \E'$ est nulle au dessus de $\U$ et son support est de dimension strictement inférieure à celle de $Z$.
D'où le résultat.
\end{proof}

\begin{prop}\label{prop-conjB-eq}
  Les assertions suivantes sont équivalentes :
  \begin{enumerate}
    \item \label{prop-conjB-eqi}La conjecture $B$ est exacte ;
    \item \label{prop-conjB-eqii} Pour tout $\V$-schéma formel lisse $\X$,
    un $F$-complexe $\E$ de $F\text{-}D ^\mathrm{b} _\mathrm{coh} ( \D ^\dag _{\X,\,\Q})$ est holonome si et seulement s'il est surcohérent ;
    \item \label{prop-conjB-eqiii} Pour tout $\V$-schéma formel lisse $\X$,
    un $F$-complexe $\E$ de $F\text{-}D ^\mathrm{b} _\mathrm{coh} ( \D ^\dag _{\X,\,\Q})$
    est holonome si et seulement s'il est surholonome.
  \end{enumerate}
\end{prop}
\begin{proof}
Prouvons d'abord $\ref{prop-conjB-eqi}) \Rightarrow \ref{prop-conjB-eqii}) $.
Comme la conjecture $B$ est supposée vraie, il est clair qu'un
$F$-complexe de $F\text{-}D ^\mathrm{b} _\mathrm{hol} ( \D ^\dag _{\X,\,\Q})$
soit surcohérent.
Réciproquement, soit $\E \in F\text{-}D ^\mathrm{b} _\mathrm{surcoh} ( \D ^\dag _{\X,\,\Q})$.
Pour tout entier $l$, $\mathcal{H} ^l (\E)$ sont alors des
$F\text{-}\D ^\dag _{\X,\,\Q}$-modules surcohérents (voir \ref{defi-surcoh}).
En particulier, ils sont alors à fibres extraordinaires finies.
Il résulte alors du lemme \ref{lemm-conjB-eq} que
pour tout entier $l$, $\mathcal{H} ^l (\E)$ sont des
$F\text{-}\D ^\dag _{\X,\,\Q}$-modules holonomes. Ainsi, $\E$ est holonome.

Vérifions $\ref{prop-conjB-eqii}) \Rightarrow \ref{prop-conjB-eqiii}) $.
Comme l'holonomie est stable par foncteur dual (voir \cite[III.4.4]{virrion}), l'assertion \ref{prop-conjB-eqii})
implique que la surcohérence est stable par foncteur dual.
Il en résulte que les notions de surcohérence et de surholonomie sont équivalentes (voir la deuxième remarque \ref{rema-hol}).
D'où : $\ref{prop-conjB-eqii}) \Rightarrow \ref{prop-conjB-eqiii}) $.

Enfin, l'implication $\ref{prop-conjB-eqiii}) \Rightarrow \ref{prop-conjB-eqi}) $
  découle de la stabilité par image inverse extraordinaire de la surholonomie (\ref{stabinminv}).

\end{proof}

\begin{rema}
\label{rema-prop-conjB-eq}
  Il résulte de \ref{prop-conjB-eq} que si la conjecture $B$ est exacte alors les notions
  d'holonomie, de surcohérence et de surholonomie sont équivalentes.
\end{rema}

\begin{prop}\label{stabholprodint}
  Soient $T$ diviseur de $X$,
$\E $ et $\FF \in F\text{-}D ^\mathrm{b} _\mathrm{coh} ( \D ^\dag _{\X} (\hdag T) _{\Q} )
\cap
F\text{-}D ^\mathrm{b} _\mathrm{hol} ( \D ^\dag _{\X,\Q} )$.
Si la conjecture $B$ est validée, alors
\begin{enumerate}
\item \label{stabholprodint.i)} le $F$-complexe $\DD _T (\E)$ appartient à
$F\text{-}D ^\mathrm{b} _\mathrm{coh} ( \D ^\dag _{\X} (\hdag T) _{\Q} )
\cap
F\text{-}D ^\mathrm{b} _\mathrm{hol} ( \D ^\dag _{\X,\Q} )$,
\item
\label{stabholprodint.ii)}
 le $F$-complexe  $\E
 \smash{\overset{\L}{\otimes}}^{\dag} _{\O _{\X } ( \hdag T ) _{\Q}}
  \FF$ (voir \ref{otimes-coh})
  appartient à $F\text{-}D ^\mathrm{b} _\mathrm{coh} ( \D ^\dag _{\X} (\hdag T) _{\Q} )
\cap
F\text{-}D ^\mathrm{b} _\mathrm{hol} ( \D ^\dag _{\X,\Q} )$.
\end{enumerate}
\end{prop}
\begin{proof}

D'après \ref{4312Be1}, comme $\E$ et $(\hdag T) (\E)$ sont des complexes à cohomologie
$\D ^\dag _{\X} (\hdag T) _{\Q}$-cohérente,
le morphisme canonique $\E \rightarrow (\hdag T) (\E)$ est un isomorphisme
car il l'est en dehors de $T$.
Cela implique l'isomorphisme :
$\DD _T   (\E) \liso \DD _T (\hdag T)  (\E)$.
Le foncteur dual commutant à l'extension des scalaires,
il en résulte  $\DD _T   (\E) \riso  (\hdag T) \DD  (\E) $.
Le foncteur $\DD$ stabilisant l'holonomie, $\DD  (\E)$ est holonome.
Comme la conjecture $B$ équivaut à la conjecture $C$,
elle implique la stabilité de l'holonomie par le foncteur $(\hdag T)$.
Ainsi, $(\hdag T) \DD  (\E) $ est holonome.
D'où \ref{stabholprodint.i)}).

Passons à la deuxième propriété.
Soient $p _1$, $p _2$ : $\X \times _\S \X \rightarrow \X$ les projections respectives à gauche et à droite,
  $\delta $ : $ \X \hookrightarrow \X \times _\S \X$ l'immersion diagonale.
Comme en particulier $\E,\,\FF \in F\text{-}D ^\mathrm{b} _\mathrm{coh} ( \D ^\dag _{\X,\Q} )$,
il existe
  $\E ^{(\bullet)},\,\FF ^{(\bullet)}
\in \smash[b]{\underset{^{\longrightarrow }}{LD }}  ^\mathrm{b} _{\Q, \mathrm{coh}}
( \smash{\widehat{\D}} _{\X} ^{(\bullet)})$ tels que
$\E  \riso \underset{\longrightarrow}{\lim} \E ^{(\bullet)} $,
$\FF  \riso \underset{\longrightarrow}{\lim} \FF ^{(\bullet)} $ (voir \ref{lim-coh}).
Par commutation (à un décalage près) du foncteur image inverse extraordinaire au
  produit tensoriel interne (cela résulte de \cite[3.3.1]{Be2} ou de \cite[1.2.22]{caro_surcoherent}),
  on obtient :
\begin{equation}
  \label{stabholprodint-iso}
  \delta ^! ( \E ^{(\bullet)} \smash{\overset{\L}{\boxtimes}} _{\O _\S}^{\dag} \FF ^{(\bullet)})
=
  \delta ^! ( p _1 ^! (\E ^{(\bullet)} ) \smash{\overset{\L}{\otimes}}^{\dag}_{\O _{\X \times _\S \X,\, \Q}} p _2^!(\FF ^{(\bullet)})
  [-2 d _{\X}])
  \riso
  \E ^{(\bullet)} \smash{\overset{\L}{\otimes}}^{\dag}_{\O _{\X ,\, \Q}} \FF ^{(\bullet)} [- d _{\X}] .
\end{equation}
Or, d'après \ref{courbe-1171}, le foncteur
$(\hdag T)$ est canoniquement isomorphe au foncteur
$\O _{\X } ( \hdag T ) _{\Q}
 \smash{\overset{\L}{\otimes}}^{\dag} _{\O _{\X ,\Q}}
-$.
On en déduit l'isomorphisme :
\begin{equation}
  \label{stabholprodint-iso2} (\hdag T) (\E ^{(\bullet)} \smash{\overset{\L}{\otimes}}^{\dag}_{\O _{\X ,\, \Q}} \FF ^{(\bullet)} )
\riso
(\hdag T) (\E ^{(\bullet)}) \smash{\overset{\L}{\otimes}}^{\dag}_{\O _{\X } (\hdag T) _{\Q}} (\hdag T) (\FF ^{(\bullet)} ).
\end{equation}

Comme l'holonomie est stable par produit tensoriel externe (voir \cite[5.3.5.(v)]{Beintro2}),
par image inverse extraordinaire et par foncteur de localisation,
  il résulte de \ref{stabholprodint-iso} et de \ref{stabholprodint-iso2},
  que
  $\underset{\longrightarrow}{\lim}\,( (\hdag T) (\E ^{(\bullet)})
  \smash{\overset{\L}{\otimes}}^{\dag}_{\O _{\X } (\hdag T) _{\Q}} (\hdag T) (\FF ^{(\bullet)} ))$
  est holonome.
Comme $\E \riso (\hdag T) (\E) \riso \underset{\longrightarrow}{\lim}\, (\hdag T) (\E ^{(\bullet)})$
et aussi
$\FF \riso \underset{\longrightarrow}{\lim}\, (\hdag T) (\FF ^{(\bullet)})$, on obtient :
  $$\underset{\longrightarrow}{\lim}\, (
  (\hdag T) (\E ^{(\bullet)})
  \smash{\overset{\L}{\otimes}}^{\dag}_{\O _{\X } (\hdag T) _{\Q}} (\hdag T) (\FF ^{(\bullet)} )
  )
  \riso
\underset{\longrightarrow}{\lim}\, (\hdag T) (\E ^{(\bullet)})
  \smash{\overset{\L}{\otimes}}^{\dag}_{\O _{\X } (\hdag T) _{\Q}}
  \underset{\longrightarrow}{\lim}\, (\hdag T) (\FF ^{(\bullet)})
  \riso
  \E
 \smash{\overset{\L}{\otimes}}^{\dag} _{\O _{\X } ( \hdag T ) _{\Q}}
  \FF.$$
Ainsi, $\E
 \smash{\overset{\L}{\otimes}}^{\dag} _{\O _{\X } ( \hdag T ) _{\Q}}
  \FF$ est holonome.
\end{proof}

\begin{lemm}\label{lemmxzz'}
Soient $Z$ un sous-schéma fermé intègre d'un $k$-schéma lisse $X$ et $Z'$ un sous-schéma fermé de $Z$.
On suppose que $\dim Z' < \dim Z < \dim X$.
Il existe alors un diviseur $T$ de $X$ tel que $T \supset Z' $ et $T \not \supset Z$.
\end{lemm}
\begin{proof}
  Soit $(X _i ) $ un recouvrement fini d'ouverts affines de $X$.
  En choisissant un élément non nul de l'idéal de définition de $Z'\cap X _i$ dans $X _i$ qui n'appartient pas
  à celui de $Z\cap X _i$ dans $X _i$, on lui associe un diviseur $T _i$ de $X _i$ tel que
  $T _i \supset Z'\cap X _i$ et $T _i \not \supset Z\cap X _i$. En notant $\overline{T} _i$ l'adhérence schématique
  de $T _i$ dans $X$, on pose $T := \cup _i \overline{T} _i$. On a $T \supset Z'$. Supposons par l'absurde
  que $T \supset Z$. Comme $Z$ est intègre, il existe un $i$ tel que $\overline{T} _i \supset Z$.
  D'où $T _i = \overline{T} _i \cap X _i \supset Z \cap X _i$. Contradiction.
\end{proof}

\begin{lemm}\label{bidualetale}
  Soient $b $ : $ \mathfrak{V} '\rightarrow \mathfrak{V}$ un morphisme étale de
  $\V$-schémas formels lisses de dimension relative nulle et
  $\E \in D ^\mathrm{b} _\mathrm{coh} (\D ^\dag _{\mathfrak{V},\Q} )$.
  On dispose de l'isomorphisme canonique
  $ b ^! \DD _{\mathfrak{V}}(\E) \riso   \DD _{\mathfrak{V}'}b ^! (\E) $ compatible à Frobenius.
\end{lemm}
\begin{proof}
  Comme $b$ est étale, $b ^! (\E) =
  \D ^\dag _{\mathfrak{V}',\Q} \otimes _{b ^{-1} \D ^\dag _{\mathfrak{V},\Q}} b ^{-1} \E $.
  Notons $(\D ^\dag _{\mathfrak{V},\Q} \otimes _{\O _{\mathfrak{V},\Q}} \omega ^{-1} _{\mathfrak{V},\Q})  _{\mathrm{t}} $
  le $\D ^\dag _{\mathfrak{V},\Q}$-bimodule à gauche égal à
$\D ^\dag _{\mathfrak{V},\Q} \otimes _{\O _{\mathfrak{V},\Q}} \omega ^{-1} _{\mathfrak{V},\Q}$ mais dont on a
inversé les structures droite et gauche. De même avec des primes.
Avec ces notations, on obtient alors l'isomorphisme canonique :
$b ^! _\mathrm{d}
  ((\D ^\dag _{\mathfrak{V},\Q} \otimes _{\O _{\mathfrak{V},\Q}} \omega ^{-1} _{\mathfrak{V},\Q})  _{\mathrm{t}})
  \riso
  (\D ^\dag _{\mathfrak{V}',\Q} \otimes _{\O _{\mathfrak{V}',\Q}} \omega ^{-1} _{\mathfrak{V}',\Q})  _{\mathrm{t}}$
où l'indice {\og d\fg} signifie que nous choisissons la structure {\it droite} de $\D ^\dag _{\mathfrak{V},\Q}$-module à gauche
sur $(\D ^\dag _{\mathfrak{V},\Q} \otimes _{\O _{\mathfrak{V},\Q}} \omega ^{-1} _{\mathfrak{V},\Q})  _{\mathrm{t}}$
pour calculer le foncteur $b ^!$.
On bénéficie des isomorphismes :
\begin{gather}
\notag
  b ^! \R \mathcal{H} om _{\D ^\dag _{\mathfrak{V},\Q}} (\E,
  (\D ^\dag _{\mathfrak{V},\Q} \otimes _{\O _{\mathfrak{V},\Q}} \omega ^{-1} _{\mathfrak{V},\Q})  _{\mathrm{t}}) [d _V]
  \riso
\R \mathcal{H} om _{b ^{-1}\D ^\dag _{\mathfrak{V},\Q}} ( b ^{-1} (\E),
  b ^! _\mathrm{d}
  ((\D ^\dag _{\mathfrak{V},\Q} \otimes _{\O _{\mathfrak{V},\Q}} \omega ^{-1} _{\mathfrak{V},\Q})  _{\mathrm{t}}))
  [d _V]
  \riso
  \\
  \notag
  \riso
  \R \mathcal{H} om _{b ^{-1}\D ^\dag _{\mathfrak{V},\Q}} ( b ^{-1} (\E),
  (\D ^\dag _{\mathfrak{V}',\Q} \otimes _{\O _{\mathfrak{V}',\Q}} \omega ^{-1} _{\mathfrak{V}',\Q})  _{\mathrm{t}})[d _V]
  \riso
\R \mathcal{H} om _{\D ^\dag _{\mathfrak{V}',\Q}} ( b ^! (\E),
  (\D ^\dag _{\mathfrak{V}',\Q} \otimes _{\O _{\mathfrak{V}',\Q}} \omega ^{-1} _{\mathfrak{V}',\Q})  _{\mathrm{t}})[d _V],
\end{gather}
où, pour calculer les foncteurs
$\R \mathcal{H} om $, nous prenons les structures gauches respectives (des bimodules à gauche)
et où le premier isomorphisme se déduit par exemple de \cite[2.1.12.(i)]{caro_comparaison}.
On conclut ensuite via les deux isomorphismes de transposition compatibles à Frobenius
$(\D ^\dag _{\mathfrak{V},\Q} \otimes _{\O _{\mathfrak{V},\Q}} \omega ^{-1} _{\mathfrak{V},\Q})  _{\mathrm{t}}
\riso \D ^\dag _{\mathfrak{V},\Q} \otimes _{\O _{\mathfrak{V},\Q}} \omega ^{-1} _{\mathfrak{V},\Q}$ et
$(\D ^\dag _{\mathfrak{V}',\Q} \otimes _{\O _{\mathfrak{V}',\Q}} \omega ^{-1} _{\mathfrak{V}',\Q})  _{\mathrm{t}}
\riso \D ^\dag _{\mathfrak{V}',\Q} \otimes _{\O _{\mathfrak{V}',\Q}} \omega ^{-1} _{\mathfrak{V}',\Q}$
(voir \cite[2.1.7]{caro_comparaison}).
\end{proof}

\begin{vide}
  \label{bidualetale2}
Soient $b $ : $ \mathfrak{V} '\rightarrow \mathfrak{V}$ un morphisme fini, étale et surjectif de
  $\V$-schémas formels lisses et
  $\E \in D ^\mathrm{b} _\mathrm{coh} (\D ^\dag _{\mathfrak{V},\Q} ) \cap
  D ^\mathrm{b} _\mathrm{coh} (\O _{\mathfrak{V},\Q} )$.
Le foncteur $b _+$ :
$D ^\mathrm{b} _\mathrm{coh} (\D ^\dag _{\mathfrak{V}',\Q} ) \rightarrow
D ^\mathrm{b} _\mathrm{coh} (\D ^\dag _{\mathfrak{V},\Q} )$
est adjoint à gauche (resp. à droite) de $b ^!$ (resp. $b ^+$) :
$D ^\mathrm{b} _\mathrm{coh} (\D ^\dag _{\mathfrak{V},\Q} ) \rightarrow
D ^\mathrm{b} _\mathrm{coh} (\D ^\dag _{\mathfrak{V}',\Q} )$.
Or, d'après \ref{bidualetale}, $b ^+ \riso b ^!$.
On obtient donc les morphismes d'adjonction
$\E \rightarrow b _+ b ^! (\E) \rightarrow \E$.
Le composé $\E \rightarrow \E$ est la multiplication par le degré de $b$
et est en particulier un isomorphisme.

En effet, avec le théorème de dualité de Grothendieck, le foncteur $b _*$ :
$D ^\mathrm{b} _\mathrm{coh} (\D ^\dag _{\mathfrak{V}',\Q} ) \cap
  D ^\mathrm{b} _\mathrm{coh} (\O _{\mathfrak{V}',\Q} ) \rightarrow
D ^\mathrm{b} _\mathrm{coh} (\D ^\dag _{\mathfrak{V},\Q} ) \cap
  D ^\mathrm{b} _\mathrm{coh} (\O _{\mathfrak{V},\Q} )$
est adjoint à droite et à gauche de
$b ^*$ : $D ^\mathrm{b} _\mathrm{coh} (\D ^\dag _{\mathfrak{V},\Q} ) \cap
  D ^\mathrm{b} _\mathrm{coh} (\O _{\mathfrak{V},\Q} )
  \rightarrow
D ^\mathrm{b} _\mathrm{coh} (\D ^\dag _{\mathfrak{V}',\Q} ) \cap
  D ^\mathrm{b} _\mathrm{coh} (\O _{\mathfrak{V}',\Q} )$.
De plus, le composé
$\E \rightarrow b_* b ^* (\E) \rightarrow \E$
induit par adjonction, est la multiplication par le degré de $b$.
Comme $b_* = b_+$, par unicité des foncteurs adjoints, il en découle que
$\E \rightarrow b _+ b ^! (\E) \rightarrow \E$ est la multiplication par le degré de $b$.

\end{vide}

\begin{prop}
\label{conja'b->a}
  Les conjectures $A'$ et $B$ impliquent la conjecture $A$.
\end{prop}
\begin{proof}
 Soit $\E$ un objet de $F\text{-}D ^\mathrm{b} _\mathrm{hol} ( \D ^\dag _{\X,\,\Q})$ tel que son support noté $Z$ soit propre sur $\Y$.
 Notons $u_0$ : $Z \hookrightarrow X$ l'immersion fermée canonique.

{\it \'Etape $1.$} Dans un premier temps, supposons que le support $Z$ de $\E$ soit lisse.
Par hypothèse,
il existe une immersion $\rho$ : $\X \hookrightarrow \PP$ avec $\PP$ un $\V$-schéma formel propre et lisse.
On pose $\delta=(\rho, f)$ : $\X \hookrightarrow \PP \times \Y$.
Comme la projection $\PP \times \Y \rightarrow \Y$ est propre, via la conjecture $A'$ que nous supposons exacte,
pour prouver que $f _+ (\E)$
est holonome, il suffit de vérifier que $\delta _+ (\E)$ est holonome.
Comme cela est local en $\PP \times \Y$, on peut supposer $\PP \times \Y$ affine.
Puisque $Z$ et $P \times Y$ sont propres sur $Y$, la $Y$-immersion
$\delta _0\circ u _0= (\rho _0 \circ u _0,f_0\circ u _0)\,:\,Z \hookrightarrow P \times Y$
est propre, donc fermée.
Il en résulte que $Z$ est aussi affine et donc que
$Z\hookrightarrow X$ se relève en un morphisme $u$ : $\ZZ \hookrightarrow \X$
de $\V$-schémas formels lisses. D'après l'analogue $p$-adique de Berthelot
du théorème de Kashiwara (voir \ref{kashiwara}),
$\E$ est de la forme $u _+ (\FF)$, avec $\FF \in F\text{-}D ^\mathrm{b} _\mathrm{hol} ( \D ^\dag _{\ZZ,\,\Q})$.
D'où: $\delta _+ (\E) \riso \delta _+ u _+ (\FF) \riso (\delta \circ u) _+ (\FF)$.
Comme $\delta \circ u$ est une immersion fermée,
on conclut via le théorème de Kashiwara (voir \ref{kashiwara}) que $\delta _+ (\E) $ est holonome.
\\

{\it \'Etape $2$.} Prouvons maintenant, par récurrence sur la dimension de $Z$, que
 $f _+ (\E)$ appartient à $F\text{-}D ^\mathrm{b} _\mathrm{hol} ( \D ^\dag _{\Y,\,\Q})$.
 Lorsque $\dim Z =0$, $Z$ est lisse et cela résulte alors du premier cas.
 Supposons donc $\dim Z > 0$ et le théorème vrai lorsque $\dim Z$ est strictement inférieure.
\\

{\it \'Etape $3$ (Réduction au cas où $\E$ est {\og associé\fg} à un $F$-isocristal surconvergent).}
 En considérant la deuxième suite spectrale d'hypercohomologie du foncteur
 $f _+$, on se ramène au cas où $\E$ est réduit à un terme.
Soient $Z _1, \dots, Z _n$ les composantes irréductible de $Z$.
Quitte à les réordonner, supposons $\dim Z _1 = \dim Z$.
Posons $Z ' = \cup _{i=2,\dots , n} Z _i$.
Comme on suppose la conjecture $B$ vraie,
les complexes du triangle distingué de Mayer-Vietoris (voir \ref{eq1mayer-vietoris}) :
\begin{equation}
\label{conja'b->a-MY}
  \R \underline{\Gamma} ^\dag _{Z _1\cap Z'}(\E ) \rightarrow
  \R \underline{\Gamma} ^\dag _{Z _1}(\E ) \oplus
\R \underline{\Gamma} ^\dag _{Z ' }(\E )  \rightarrow
\R \underline{\Gamma} ^\dag _{Z }(\E ) \rightarrow
\R \underline{\Gamma} ^\dag _{Z _1\cap Z' }(\E )[1]
\end{equation}
sont holonomes et à support propre sur $Y$.
De plus, comme $\E$ est à support dans $Z$,
le morphisme canonique
$\R \underline{\Gamma} ^\dag _{Z }(\E ) \rightarrow \E$ est un isomorphisme (voir \cite[2.2.9]{caro_surcoherent}).
En appliquant $f _+$ à \ref{conja'b->a-MY},
quitte à procéder à une seconde récurrence sur le nombre de composantes irréductibles de dimension
$\dim Z$,
on se ramène alors au cas où $Z$ est intègre.

D'après le théorème de désingularisation de de Jong (\cite{dejong}),
 il existe un morphisme
  $h_0$ : $Z'\rightarrow Z$ surjectif, projectif, génériquement fini et étale
  tel que $Z'$ soit intègre et lisse.
  Le morphisme $h_0$ se décompose en une immersion fermée $Z ' \hookrightarrow \P ^r _Z $ suivi de la projection
  $\P ^r _Z \rightarrow Z$. En notant $\X' $ le complété $p$-adique de $\P ^r _{\X}$,
$g$ : $\X '\rightarrow \X$ la projection canonique,
$u' _0$ : $Z '\hookrightarrow X'$ l'immersion fermée induite,
on obtient l'égalité
 $g _0 \circ u ' _0= u _0\circ h_0$.

Il existe un ouvert $\U$ de $\X$
tel que $V:=Z \cap U$ soit un ouvert lisse non vide de $Z$ et
tel que $V$ soit inclus dans un ouvert affine et lisse $W$ de $Z$
tel que $h _0^{-1} (W) \rightarrow W $ soit
fini et étale. Notons $\U' := g ^{-1} (\U)$, $a$ : $\U' \rightarrow \U$ le morphisme induit par $g$.
Il existe alors des morphismes
$b $ : $\mathfrak{V} ' \rightarrow \mathfrak{V}$,
$v$ : $\mathfrak{V} \hookrightarrow \U$ et $v '$ :
$\mathfrak{V} '\hookrightarrow \U'$ de $\V$-schémas formels lisses relevant de respectivement
$h _0^{-1} (V) \rightarrow V$, $V \hookrightarrow U$ et $h_0 ^{-1} (V) \hookrightarrow U' _0$
tels que $a \circ v' = v \circ b$ (en effet, $W$ se relève, donc $V$ et $V'$ aussi, puis on choisit un relèvement $v$, puis
un autre de la forme $\mathfrak{V}' \rightarrow \U' \times _{\U} \mathfrak{V}$).

Grâce à \cite[5.3.5.(i)]{Beintro2} (ou à \ref{lemmsurcoh=>hol} et \ref{prop-conjB-eq}),
quitte à rétrécir $\U$, on peut en outre supposer
qu'il existe un isocristal convergent $G$ sur $V$ vérifiant
$\E |\U \riso v _+ \sp _{*} (G)$.
Grâce au lemme \ref{lemmxzz'} (appliqué à $X \supset Z \supset Z \setminus V$),
quitte à rétrécir $\U$,
on se ramène au cas où $T:= X \setminus U$ est le support d'un diviseur de $X$.
On pose $T' := g _0 ^{-1} (T)$.

Comme $Z'$ est intègre et lisse, que $Z' \cap T'$ ne contient pas $Z'$ (car $b_0$ est fini et étale
 et $\dim Z \cap T < \dim Z$), $Z' \cap T'$ est un diviseur de $Z'$.
 On dispose donc du foncteur pleinement fidèle  $\sp _{Z'\hookrightarrow \X',T'+}$.

 Comme $\R \underline{\Gamma} ^\dag _{T} \E$ est à support
 dans $T \cap Z$, son support est donc propre sur $Y$ et de dimension plus petite que celle de $Z$.
  Par hypothèse de récurrence, $f _+ (\R \underline{\Gamma} ^\dag _{T} \E)$ est donc holonome.
 En appliquant $f _+ $ au le triangle de localisation de $\E$ en $T$,
 on se ramène ainsi à établir l'holonomie de $f _+ ((\hdag T) \E)$.
\\

{\it \'Etape $4$} :
avec les notations de l'étape $3$, on vérifie que
$ ((\hdag T) \E)$ est un facteur direct de
 $ g _{T,+} \R \underline{\Gamma} ^\dag _{Z'} g _T ^! ( (\hdag T) \E)$.

Via le théorème de Kashiwara (voir \ref{kashiwara}),
on obtient le dernier isomorphisme de la
suite d'isomorphismes compatibles à Frobenius :
\begin{equation}
  \label{conja'b->a-iso1}
(\R \underline{\Gamma} ^\dag _{Z'} g _T ^! ( (\hdag T) \E) )|\U
\riso
v' _+ v ^{\prime!} a ^! ( \E |\U)
\riso
v' _+ v ^{\prime!} a ^! v _+ \sp _{*} (G)
\riso
v' _+ b ^! \sp _{*} (G)
\riso
v' _+  \sp _{*} (b ^* (G)).
\end{equation}
De même,
grâce au théorème de dualité relative compatible à Frobenius appliqué aux immersions fermées
$v$ et $v'$ (voir \cite[2.4.3]{caro-frobdualrel} ou son résumé),
au théorème de Kashiwara (voir \ref{kashiwara}), à \ref{bidualetale}
 et au théorème de bidualité (\cite[II.3.5]{virrion}),
on obtient respectivement le deuxième, troisième, quatrième et cinquième isomorphisme compatible à Frobenius
de la suite :
\begin{gather}
\notag
  (\DD  _{\X',T'} \circ \R \underline{\Gamma} ^\dag _{Z'} g _T ^! \DD  _{\X,T} ( (\hdag T) \E)) |\U
\riso
\DD  _{\U'} \circ v' _+ v ^{\prime!} a ^! \DD  _{\U} v _+ \sp _{*} (G)
\\
\notag
\riso
v' _+ \DD  _{\mathfrak{V}'} v ^{\prime!} a ^! v _+ \DD  _{\mathfrak{V}}  \sp _{*} (G)
\underset{\ref{kashiwara}}{\riso}
v' _+ \DD  _{\mathfrak{V}'} b ^! \DD  _{\mathfrak{V}}  \sp _{*} (G)
\\
\label{conja'b->a-iso2}
\underset{\ref{bidualetale}}{\riso}
v' _+ \DD  _{\mathfrak{V}'}  \DD  _{\mathfrak{V}'}  b ^! \sp _{*} ( G )
\riso
v' _+ b ^! \sp _{*} (G)
\riso
v' _+  \sp _{*} (b ^* (G)).
\end{gather}

On déduit de \ref{conja'b->a-iso1} et \ref{conja'b->a-iso2} que
 $\R \underline{\Gamma} ^\dag _{Z'} g _T ^! ( (\hdag T) \E)$ et
 $\DD  _{\X',T'} \circ \R \underline{\Gamma} ^\dag _{Z'} g _T ^! \DD _{\X,T} ( (\hdag T) \E)$
 sont dans l'image essentielle de
 $\sp _{Z'\hookrightarrow \X',T'+}$ (voir la remarque \ref{rema-unitesurho2}).
D'après \cite[5.3.1]{caro_devissge_surcoh},
comme d'après Kedlaya le foncteur restriction qui à un $F$-isocristal sur $Z '\setminus T'$
 surconvergent le long de $T' \cap Z'$
 associe le $F$-isocristal convergent sur $Z '\setminus T'$
 correspondant est pleinement fidèle (voir par exemple \cite[4.2.1]{kedlaya-semistableII}),
on en déduit que pour vérifier que
 $\R \underline{\Gamma} ^\dag _{Z'} g  _T ^! ( (\hdag T) \E)$ et
 $\DD _{\X',T'} \circ  \R \underline{\Gamma} ^\dag _{Z'} g _T ^! \DD _{\X,T} ( (\hdag T) \E)$ sont isomorphes,
 il suffit de le voir en dehors de $T'$, ce qui est bien le cas d'après
 \ref{conja'b->a-iso1} et \ref{conja'b->a-iso2}.

 Or, d'après \ref{f+adjf!},
 on dispose des morphismes d'adjonction
 $g _{T+} \R \underline{\Gamma} ^\dag _{Z'} g _T^! ( (\hdag T) \E) \rightarrow  (\hdag T) \E$ et
 $(\hdag T) \E \rightarrow
 g _{T,+} \DD  _{\X',T'} \circ \R \underline{\Gamma} ^\dag _{Z'} g _T ^! \DD  _{\X,T} ( (\hdag T) \E)$
 (pour le deuxième, on applique le foncteur dual et on utilise l'isomorphisme
 de dualité relative).
On obtient ainsi la suite de morphismes :
\begin{equation}
\label{conja'b->aadjadj}
  (\hdag T) \E \rightarrow  g _{T,+} \R \underline{\Gamma} ^\dag _{Z'} g _T ^! ( (\hdag T) \E)
  \rightarrow  (\hdag T) \E.
\end{equation}
 Par \ref{4312Be1}, pour prouver que ce composé est un isomorphisme, il suffit de le vérifier
 au dessus de $\U$.
 La flèche de droite de \ref{conja'b->aadjadj} restreinte à $\U$, est canoniquement isomorphe au morphisme induit par adjonction
 $a _+ v ' _+ v ^{\prime !} a ^! (\E |\U)  \rightarrow \E |\U$ (voir \ref{+!=gamma} et \ref{rema-prop-conjB-eq}).
 Par transitivité des morphismes d'adjonction (voir \cite[1.2.11]{caro-construction} ou \ref{commFrobf+}),
 il est alors canoniquement isomorphe au morphisme construit par adjonction :
$v _+ b _+ b ^{!} v ^! (\E |\U)  \rightarrow \E |\U$.
Avec le théorème de Kashiwara (voir \ref{kashiwara}), en lui appliquant
$v ^!$, on obtient le morphisme : $b _+ b ^{!} v ^! (\E |\U)  \rightarrow v ^!  (\E |\U)$.
D'un autre côté, la flèche de gauche de \ref{conja'b->aadjadj} restreinte à $\U$
est canoniquement isomorphe
au morphisme induit par adjonction
$\E |\U \rightarrow a _+ v ' _+ v ^{\prime +} a ^+ (\E |\U)$.
Par transitivité, il correspond au morphisme d'adjonction
$\E |\U \rightarrow v _+ b _+ b ^{+} v ^+ (\E |\U)$.
En lui appliquant $v ^! $, on obtient donc (modulo $v ^! v_= \riso Id $) le morphisme
$v ^!  (\E |\U) \rightarrow b _+ b ^{+} v ^+ (\E |\U)$.

Or, comme $\E \riso v _+ \sp _{*} (G)$, le théorème de Kashiwara et l'isomorphisme de dualité relative
donnent $v ^+  (\E |\U) \riso v ^!  (\E |\U) \riso \sp _{*} (G)$.
De plus, d'après \ref{bidualetale}, $b ^+ \riso b^!$.
En appliquant $v ^!$ à la restriction sur $\U$ de \ref{conja'b->aadjadj},
on obtient donc la suite de morphismes construits par adjonction :
$\sp _{*} (G)\rightarrow b _+ b ^{!} (\sp _{*} (G))\rightarrow \sp _{*} (G)$.
D'après \ref{bidualetale2}, ce composé est un isomorphisme.
On a ainsi établi que
le composé de \ref{conja'b->aadjadj} est un isomorphisme, ce qui implique que
$ ((\hdag T) \E)$ est un facteur direct de
 $ g _{T,+} \R \underline{\Gamma} ^\dag _{Z'} g _T ^! ( (\hdag T) \E)$.
\\

{\it \'Etape 5 (Conclusion)}. On déduit de l'étape $4$ que $f _{+} ((\hdag T) \E)$ est un facteur direct de
 $f _+ ( g _{T,+} \R \underline{\Gamma} ^\dag _{Z'} g _T ^! ( (\hdag T) \E))$.
Comme $f _+ g _{T,+}\riso f _+  g _{+} \riso f \circ g _+$ et
comme $\R \underline{\Gamma} ^\dag _{Z'} g _T ^! ( (\hdag T) \E))$
est à support dans $Z'$ qui est propre sur $Y$ et lisse,
d'après le cas lisse (i.e., l'étape $1$), le $F$-complexe $f _+ ( g _{T,+} \R \underline{\Gamma} ^\dag _{Z'} g _T ^! ( (\hdag T) \E))$
est holonome. D'où le résultat.

\end{proof}

\begin{theo}\label{B->A}
  La conjecture $B$ implique les conjectures $A$ et $A'$.
\end{theo}
\begin{proof}
  D'après \ref{prop-conjB-eq}, si la conjecture $B$ est exacte alors les notions de surcohérence et d'holonomie
  sont équivalentes.
  Comme la surcohérence
se préserve par l'image directe d'un morphisme propre (\cite[3.1.9]{caro_surcoherent}),
il en serait de même pour l'holonomie, i.e., la conjecture $A'$ est une conséquence de la conjecture $B$.
Par \ref{conja'b->a}, il en résulte que la conjecture $B$ implique la conjecture $A$.
\end{proof}

\bibliographystyle{smfalpha}
\bibliography{bib1}

\newcommand{\etalchar}[1]{$^{#1}$}
\providecommand{\bysame}{\leavevmode ---\ }
\providecommand{\og}{``}
\providecommand{\fg}{''}
\providecommand{\smfandname}{et}
\providecommand{\smfedsname}{\'eds.}
\providecommand{\smfedname}{\'ed.}
\providecommand{\smfmastersthesisname}{M\'emoire}
\providecommand{\smfphdthesisname}{Th\`ese}
\begin{thebibliography}{BGK{\etalchar{+}}87}

\bibitem[Ber]{Be4}
{\scshape P.~Berthelot} -- {\og {$\mathcal{D}$}-modules arithm\'etiques.
  {I}{V}. {V}ari\'et\'e caract\'eristique\fg}, en pr{\'e}paration.

\bibitem[Ber74]{Becris}
\bysame , \emph{Cohomologie cristalline des sch\'emas de caract\'eristique
  $p>0$}, Springer-Verlag, Berlin, 1974, Lecture Notes in Mathematics, Vol.
  407.

\bibitem[Ber90]{Be0}
\bysame , {\og Cohomologie rigide et th\'eorie des ${\mathcal{d}}$-modules\fg},
  $p$-adic analysis (Trento, 1989), Springer, Berlin, 1990, p.~80--124.

\bibitem[Ber96a]{Becohdiff}
\bysame , {\og Coh\'erence diff\'erentielle des alg\`ebres de fonctions
  surconvergentes\fg}, \emph{C. R. Acad. Sci. Paris S\'er. I Math.}
  \textbf{323} (1996), no.~1, p.~35--40.

\bibitem[Ber96b]{Berig}
\bysame , {\og {Cohomologie rigide et cohomologie rigide \`a support propre.
  Première partie}\fg}, {Pr\'epublication IRMAR 96-03, Universit\'e de Rennes},
  1996.

\bibitem[Ber96c]{Be1}
\bysame , {\og ${\mathcal{d}}$-modules arithm\'etiques. {I}. {O}p\'erateurs
  diff\'erentiels de niveau fini\fg}, \emph{Ann. Sci. \'Ecole Norm. Sup. (4)}
  \textbf{29} (1996), no.~2, p.~185--272.

\bibitem[Ber00]{Be2}
\bysame , {\og {$\mathcal{D}$}-modules arithm\'etiques. {I}{I}. {D}escente par
  {F}robenius\fg}, \emph{M\'em. Soc. Math. Fr. (N.S.)} (2000), no.~81,
  p.~vi+136.

\bibitem[Ber02]{Beintro2}
\bysame , {\og {Introduction \`a la th\'eorie arithm\'etique des
  {$\mathcal{D}$}-modules}\fg}, \emph{Ast\'erisque} (2002), no.~279, p.~1--80,
  Cohomologies {$p$}-adiques et applications arithm\'etiques, {II}.

\bibitem[BGK{\etalchar{+}}87]{borel}
{\scshape A.~Borel, P.-P. Grivel, B.~Kaup, A.~Haefliger, B.~Malgrange
  {\normalfont \smfandname} F.~Ehlers} -- \emph{Algebraic ${D}$-modules},
  Academic Press Inc., Boston, MA, 1987.

\bibitem[Car04a]{caro_unite}
{\scshape D.~Caro} -- {\og Coh\'erence diff\'erentielle des {$F$}-isocristaux
  unit\'es\fg}, \emph{C. R. Math. Acad. Sci. Paris} \textbf{338} (2004), no.~2,
  p.~145--150.

\bibitem[Car04b]{caro_surcoherent}
\bysame , {\og {$\mathcal{D}$}-modules arithm{\'e}tiques surcoh{\'e}rents.
  {A}pplication aux fonctions {L}\fg}, \emph{Ann. Inst. Fourier, Grenoble}
  \textbf{54} (2004), no.~6, p.~1943--1996.

\bibitem[Car05a]{caro_comparaison}
\bysame , {\og Comparaison des foncteurs duaux des isocristaux
  surconvergents\fg}, \emph{Rend. Sem. Mat. Univ. Padova} \textbf{114} (2005),
  p.~131--211.

\bibitem[Car05b]{caro-construction}
\bysame , {\og {$\mathcal{D}$-modules arithmétiques associés aux isocristaux
  surconvergents. Cas lisse}\fg}, \emph{ArXiv Mathematics e-prints} (2005).

\bibitem[Car05c]{caro-frobdualrel}
\bysame , {\og {Sur la compatibilité à Frobenius de l'isomorphisme de dualité
  relative}\fg}, \emph{ArXiv Mathematics e-prints} (2005).

\bibitem[Car06a]{caro_devissge_surcoh}
\bysame , {\og {Dévissages des $F$-complexes de $\mathcal{D}$-modules
  arithmétiques en $F$-isocristaux surconvergents}\fg}, \emph{Invent. Math.}
  \textbf{166} (2006), no.~2, p.~397--456.

\bibitem[Car06b]{caro_courbe-nouveau}
\bysame , {\og Fonctions {L} associ{\'e}es aux {$\mathcal{D}$}-modules
  arithm{\'e}tiques. {C}as des courbes\fg}, \emph{Compositio Mathematica}
  \textbf{142} (2006), no.~01, p.~169--206.

\bibitem[Car06c]{caro_log-iso-hol}
\bysame , {\og {Log-isocristaux convergents et holonomie}\fg}, \emph{ArXiv
  Mathematics e-prints} (2006).

\bibitem[Car07]{caro-2006-surcoh-surcv}
\bysame , {\og {$F$-isocristaux surconvergents et surcohérence
  différentielle}\fg}, \emph{Invent. Math.} \textbf{170} (2007), no.~3,
  p.~507--539.

\bibitem[dJ96]{dejong}
{\scshape A.~J. de~Jong} -- {\og Smoothness, semi-stability and
  alterations\fg}, \emph{Inst. Hautes \'Etudes Sci. Publ. Math.} (1996),
  no.~83, p.~51--93.

\bibitem[Elk73]{elkik}
{\scshape R.~Elkik} -- {\og Solutions d'\'equations \`a coefficients dans un
  anneau hens\'elien\fg}, \emph{Ann. Sci. \'Ecole Norm. Sup. (4)} \textbf{6}
  (1973), p.~553--603 (1974).

\bibitem[{\'E}LS93]{E-LS}
{\scshape J.-Y. {\'E}tesse {\normalfont \smfandname} B.~Le~Stum} -- {\og
  Fonctions ${L}$ associ\'ees aux ${F}$-isocristaux surconvergents. {I}.
  {I}nterpr\'etation cohomologique\fg}, \emph{Math. Ann.} \textbf{296} (1993),
  no.~3, p.~557--576.

\bibitem[Har77]{HaAG}
{\scshape R.~Hartshorne} -- \emph{Algebraic geometry}, Springer-Verlag, New
  York, 1977, Graduate Texts in Mathematics, No. 52.

\bibitem[Keda]{kedlaya-padiccohomology}
{\scshape K.~S. Kedlaya} -- {\og {p-adic cohomology}\fg},
  arXiv:math.NT/0601507.

\bibitem[Kedb]{kedlaya-semistableI}
\bysame , {\og {Semistable reduction for overconvergent F-isocrystals, I:
  Unipotence and logarithmic extensions}\fg}, To appear in Compositio
  Mathematica.

\bibitem[Kedc]{kedlaya-semistableII}
\bysame , {\og {Semistable reduction for overconvergent F-isocrystals, II: A
  valuation-theoretic approach}\fg}, To appear in Compositio Mathematica.

\bibitem[Ked04]{kedlaya_full_faithfull}
\bysame , {\og Full faithfulness for overconvergent {$F$}-isocrystals\fg},
  Geometric aspects of Dwork theory. Vol. I, II, Walter de Gruyter GmbH \& Co.
  KG, Berlin, 2004, p.~819--835.

\bibitem[NH]{huyghe_finitude_coho}
{\scshape C.~Noot-Huyghe} -- {\og {Finitude de la dimension homologique
  d'algèbres d'opérateurs différentiels faiblement complètes et à coefficients
  surconvergents}\fg}, {A paraître au journal of Algebra (2005)}.

\bibitem[Tsu02]{tsumono}
{\scshape N.~Tsuzuki} -- {\og Morphisms of ${F}$-isocrystals and the finite
  monodromy theorem for unit-root ${F}$-isocrystals\fg}, \emph{Duke Math. J.}
  \textbf{111} (2002), no.~3, p.~385--418.

\bibitem[Vir00]{virrion}
{\scshape A.~Virrion} -- {\og Dualit\'e locale et holonomie pour les
  {$\mathcal{D}$}-modules arithm\'etiques\fg}, \emph{Bull. Soc. Math. France}
  \textbf{128} (2000), no.~1, p.~1--68.

\bibitem[Vir04]{Vir04}
\bysame , {\og Trace et dualit\'e relative pour les {$\mathcal{D}$}-modules
  arithm\'etiques\fg}, Geometric aspects of Dwork theory. Vol. I, II, Walter de
  Gruyter GmbH \& Co. KG, Berlin, 2004, p.~1039--1112.

\end{thebibliography}

\bigskip
\noindent Daniel Caro\\
Arithmétique et Géométrie algébrique, Bât. 425\\
Université Paris-Sud\\
91405 Orsay Cedex\\
France.\\
email: daniel.caro@math.u-psud.fr

\end{document}